# HIGHER SCHOOL OF ECONOMICS
NATIONAL RESEARCH UNIVERSITY

*Sergey Shvydun*

# NORMATIVE PROPERTIES OF MULTI-CRITERIA CHOICE PROCEDURES AND THEIR SUPERPOSITIONS: I

We consider different choice procedures such as scoring rules, rules, using majority relation, value function and tournament matrix, which are used in social and multi-criteria choice problems. We focus on the study of the properties that show how the final choice is changed due to changes of preferences or a set of feasible alternatives. As a result a theorem is provided showing which normative properties (rationality, monotonicity, non-compensability) are satisfied for the given choice procedures.

*Shvydun Sergey* – National Research University Higher School of Economics, Moscow, Russia; V.A. Trapeznikov Institute of Control Sciences of Russian Academy of Sciences, Moscow, Russia.

The work was partially financed by the International Laboratory of Decision Choice and Analysis (DeCAn Lab) of the National Research University Higher School of Economics and by V.A. Trapeznikov Institute of Control Sciences of Russian Academy of Sciences. The research was partially supported also by the Russian Foundation for Basic Research (RFBR) (grant № 12-00-00226 "Choice models based on superposition").
The author would like to thank Professor F. Aleskerov for the formulation of the research problem. The author also thank Professor V. Volsky for his comments and suggestions, which helped the author to improve the paper.

# Introduction

The interest in problems of general choice theory can be explained by the fact that many problems of applied mathematics and control theory boil down to choosing the «best» alternatives, in some sense, from a given set.

The choice problem can be formally defined as follows [1]. Let $A$ be a finite set of alternatives, $|A|>2$. Any subset $X \epsilon 2^A$ may be presented for choice; it is called a presentation. Denote by $C(\cdot)$ a choice function that performs the mapping $2^A \to 2^A$ with the restriction $C(X) \subseteq X$ for any $X \epsilon 2^A$. A choice consists in the selection according to some rule from some presentation $X$ of the non-empty subset of alternatives $Y \subseteq X$ (non-empty subset of «best» alternatives).

The simplest example of the choice problem is the choice on one criterion. In that case, some constraints are set and the choice consists in minimization or maximization of a given criterion under these constraints. Thus, the choice of alternatives is accomplished on pre-defined optimality criterion with the use of some extremization procedure.

However, most real-life choice problems usually deal with multiple criteria. Currently, there are many methods of transforming multi-criteria optimization problem into a single criterion optimization problem but, unfortunately, these methods cannot always be applied. For instance, some criteria may conflict with others, and in that case, the choice problem lies in providing a compromise set of alternatives, so the end-user can choose the most appropriate of them.

There are a lot of different multi-criteria choice procedures that allow to choose and rank alternatives from the initial set. All these choice procedures can be divided into 5 following groups [1-9]:
1. Scoring rules;
2. Rules, using majority relation;
3. Rules, using value function;
4. Rules, using tournament matrix;
5. q-Paretian rules,

from which we consider only first 4 groups.

The change of presentation, a set of criteria or criterial values of some alternatives can affect the final choice. Consequently, there is a need for more detailed study of existing choice procedures and for understanding which of them can be used in a particular case.

In this paper we consider the situation when criterial estimates of alternatives are presented in the form of rankings. We focus on the study of the properties of different multi-criteria choice procedures in order to define which normative conditions are satisfied for the given choice procedures. The paper is divided into several sections. First, the description of different choice procedures is given. Second, normative conditions are presented. Finally, we study the properties of choice procedures and give a theorem that shows how stable and sensible is a set of alternatives obtained after applying some choice procedure.

# A description of choice procedures

In Table 1 we provide a list of multi-criteria choice procedures which are studied in this work.

**Table 1. A list of studied choice procedures**

| № | Choice procedure | Type of choice procedure |
|---|---|---|
| 1 | The simple majority rule | Scoring rules |
| 2 | The plurality rule | |
| 3 | The inverse plurality rule | |
| 4 | The q-Approval rule | |
| 5 | The run-off procedure | |
| 6 | The Hare rule (the Ware procedure) | |
| 7 | The Borda rule | |
| 8 | The Black procedure | |
| 9 | The inverse Borda rule | |
| 10 | The Nanson rule | |
| 11 | The Coombs procedure | |
| 12 | Minimal dominant set | Rules, using majority relation |
| 13 | Minimal undominated set | |
| 14 | Minimal weakly stable set | |
| 15 | The Fishburn rule | |
| 16 | Uncovered set I | |
| 17 | Uncovered set II | |
| 18 | The Richelson rule | |
| 19 | The Condorcet winner | |
| 20 | The core | |
| 21 | k-stable set | |
| 22 | The threshold rule | Rules, using value function |
| 23 | The Copeland rule 1 | |
| 24 | The Copeland rule 2 | |
| 25 | The Copeland rule 3 | |
| 26 | The super-threshold rule | |
| 27 | The minimax procedure | Rules, using tournament matrix |
| 28 | The Simpson procedure | |

Let us give a definition of the rules provided in Table 1.

*Definition 1 [10].* The relation $P$ is called
- the partial order if $P$ is irreflexive ($\forall x \in A\ (x,x) \notin P$) and transitive ($\forall x, y, z \in A\ xPy$ and $yPz \Rightarrow xPz$);
- the weak order if $P$ is a partial order satisfying the condition for negative transitivity ($\forall x, y, z \in A\ (x, y) \notin P$ and $(y, z) \notin P \Rightarrow (x, z) \notin P$); and
- the linear order if $P$ is a connected ($\forall x, y \in A\ x \neq y \Rightarrow xPy$ or $yPx$) weak order.

*Definition 2 [11].* A set of binary relations $(P_1, \ldots, P_n)$ is called a preference profile and is denoted as $\vec{P}$ or $\{P\}_1^n$.

Let us denote a set of profiles consisting of linear orders as $\mathcal{P}$. Obviously, if the number of alternatives of the set $A$ is equal to $m$, i.e., $|A| = m$, then $|\mathcal{P}| = (m!)^n$. We also introduce the notion of the contraction of the profile $\vec{P}$ onto a set $X \subseteq A$.

*Definition 3 [10].* A contraction of the profile $\vec{P}$ onto a set $X \subseteq A$, $X \neq \emptyset$, is a profile $\vec{P}_X = (P_1/X, \ldots, P_n/X)$, where $P_i/X = P_i \cap (X \times X)$.

Denote by $C(\vec{P}_X, X)$ a choice function which is a mapping of preference profile $\vec{P}_X$ on the set of alternatives, i.e., $C: \mathcal{P} \times 2^A \to 2^A$, that satisfies the condition $C(\vec{P}_X, X) \subseteq X$ for any $X \in 2^A$. The q-Approval rule and the k-stable set which are described below also depend on parameters $q$ and $k$. We denote these choice procedures by $C(\vec{P}_X, X, q)$ and $C(\vec{P}_X, X, k)$, respectively.

*Definition 4 [1].*
a) Majority relation is a binary relation µ which is constructed as
$$x\mu y \Leftrightarrow card\{i \in N | xP_iy\} > card\{i \in N | yP_ix\}.$$
b) A directed graph of majority relation $\mu$ is called a majority graph.

*Definition 5 [10].*
a) Upper contour set of an alternative $x$ for a given profile $\vec{P}$ is the set $D(x, \vec{P})$ such that
$$D(x, \vec{P}) = \{y \in A | y\mu x\}.$$
b) Lower contour set of $x$ for a given profile $\vec{P}$ is the set $L(x, \vec{P})$ such that
$$L(x, \vec{P}) = \{y \in A | x\mu y\}.$$

In this paper we consider criterial scales of a special type according to which all alternatives are arranged in $n$ linear orders $P_i$. Let us now describe different multi-criteria choice procedures which are used in this case [1-9].

**1. The simple majority rule**

Choose alternatives that have been admitted to be the best by 50% + 1 criteria, i.e.,
$$a \in C(\vec{P}_X, X) \Leftrightarrow \left[ n^+(a, \vec{P}_X) > \frac{1}{2} \cdot \sum_{x \in X} n^+(x, \vec{P}_X) = \frac{1}{2} \cdot |\vec{P}_X| \right],$$
where $n^+(a, \vec{P}_X) = card\{i \in N \mid \forall y \in X \setminus \{a\}: aP_iy\}$.

*Example.* Let $X = \{a, b, c, d\}$ and the profile $\vec{P}_X$ is the following

| $P_1$ | $P_2$ | $P_3$ | $P_4$ | $P_5$ | $P_6$ | $P_7$ | $P_8$ |
|---|---|---|---|---|---|---|---|
| a | a | a | b | b | d | d | d |
| b | b | c | c | c | c | b | b |
| c | c | d | a | a | b | c | c |
| d | d | b | d | d | a | a | a |

In this example $n^+(a, \vec{P}_X) = 3$, $n^+(b, \vec{P}_X) = 2$, $n^+(c, \vec{P}_X) = 0$, $n^+(d, \vec{P}_X) = 3$, and $|\vec{P}_X| = \sum_{x \in X} n^+(x, \vec{P}_X) = 8$. Since $\nexists x \in X: n^+(x, \vec{P}_X) > 4$, the choice is empty, i.e., $C(\vec{P}_X, X) = \emptyset$.

## 2. The plurality rule

Choose alternatives that have been admitted to be the best by the maximum number of criteria, i.e.,

$$a \in C(\vec{P}_X, X) \Leftrightarrow [\forall x \in X: n^+(a, \vec{P}_X) \geq n^+(x, \vec{P}_X)],$$

where $n^+(a, \vec{P}_X) = card\{i \in N \mid \forall y \in X\setminus\{a\}: aP_iy\}$.

*Consider the previous example.* $n^+(a, \vec{P}_X) = 3$, $n^+(b, \vec{P}_X) = 2$, $n^+(c, \vec{P}_X) = 0$, $n^+(d, \vec{P}_X) = 3$. According to the plurality rule the alternatives *a* and *d* will be chosen, i.e., $C(\vec{P}_X, X) = \{a, d\}$.

## 3. The inverse plurality rule

The alternative, which is regarded as the worst by the minimum number of criteria, is chosen, i.e.,

$$a \in C(\vec{P}_X, X) \Leftrightarrow [\forall x \in X: n^-(a, \vec{P}_X) \leq n^-(x, \vec{P}_X)],$$

where $n^-(a, \vec{P}_X) = card\{i \in N \mid \forall y \in X\setminus\{a\}: yP_ia\}$.

*Consider the previous example.* $n^-(a, \vec{P}_X) = 3$, $n^-(b, \vec{P}_X) = 1$, $n^-(c, \vec{P}_X) = 0$, $n^-(d, \vec{P}_X) = 4$. According to the inverse plurality rule the alternative *c* will be chosen, i.e., $C(\vec{P}_X, X) = \{c\}$.

## 4. The q-Approval Rule

Let us define

$$n^+(a, \vec{P}_X, q) = card\{i \in N \mid card\{y \in X \mid yP_ia\} \leq q - 1\},$$

where $n^+(a, \vec{P}_X, q)$ is the number of criteria for which the alternative *a* is placed on one of *q* first places in the ordering of alternatives. Thus, if $q = 1$, then *a* is the first best alternative for criterion *i*; if $q = 2$, then *a* is either first best or second best alternative, etc. The integer *q* can be called as the degree of the choice procedure.

Now we can define the q-Approval rule:

$$a \in C(\vec{P}_X, X, q) \Leftrightarrow [\forall x \in X: n^+(a, \vec{P}_X, q) \geq n^+(x, \vec{P}_X, q)],$$

i.e., the alternatives are chosen that have been admitted to be among *q* best alternatives by the maximum number of criteria.

It can be easily seen that the q-Approval rule is a direct generalization of the plurality rule; for the latter $q = 1$.

*Consider the previous example.* Let us calculate the value $n^+(x, \vec{P}_X, q)$ for each $x \in X$ depending on parameter value *q*.

| q | $n^+(a, \vec{P}_X, q)$ | $n^+(b, \vec{P}_X, q)$ | $n^+(c, \vec{P}_X, q)$ | $n^+(d, \vec{P}_X, q)$ | $C(\vec{P}_X, X, q)$ |
|---|---|---|---|---|---|
| q=1 | 3 | 2 | 0 | 3 | {a,d} |
| q=2 | 3 | 6 | 4 | 3 | {b} |
| q=3 | 5 | 7 | 8 | 4 | {c} |
| q=4 | 8 | 8 | 8 | 8 | {a,b,c,d} |

## 5. The run-off procedure

First, the simple majority rule is used. If a winner exists, the procedure stops. Otherwise, two alternatives that have been admitted to be the best by the maximum number of criteria are taken. Assuming that preferences about these two alternatives do not change, again the simple majority rule is applied. Since alternatives are arranged in $n$ linear orders, single winner (if $n$ is odd) always exists.

*Consider the previous example*. First, the simple majority rule is used. As it was mentioned before, according to the simple majority rule the choice is empty. In that case the alternatives $a$ and $d$ are taken for the next stage of the procedure. Consider now the set $X' = \{a, d\}$ and the profile $\vec{P}_{X'}$ that looks as

| $P_1$ | $P_2$ | $P_3$ | $P_4$ | $P_5$ | $P_6$ | $P_7$ | $P_8$ |
|---|---|---|---|---|---|---|---|
| $a$ | $a$ | $a$ | $a$ | $a$ | $d$ | $d$ | $d$ |
| $d$ | $d$ | $d$ | $d$ | $d$ | $a$ | $a$ | $a$ |

According to the simple majority rule $n^+(a, \vec{P}_{X'}) = 5$, $n^+(d, \vec{P}_{X'}) = 3$, $|\vec{P}_{X'}| = 8$. Consequently, the alternative $a$ will be chosen, i.e., $C(\vec{P}_X, X) = C(\vec{P}_{X'}, X') = \{a\}$.

## 6. The Hare procedure (the Ware procedure)

The Hare procedure (also known as the Ware procedure [4]) is a modification of a single transferable vote system for the case when there is a need to choose a single alternative that have been admitted to be the best by 50% + 1 criteria.

First, the simple majority rule is used. If a winner exists, the procedure stops. Otherwise, the alternative $x$ that have been admitted to be the best by the minimum number of criteria is omitted. Then the procedure again applied to the narrowed set $X' = X \setminus \{x\}$ and the profile $\vec{P}_{X'}$. The procedure stops when an alternative that has been admitted to be the best by 50% + 1 criteria is found.

*Consider the previous example*. First, the simple majority rule is used. According to it the choice is empty. In that case, the alternatives $b$ and $c$ are omitted for the next stage of the procedure, as they are admitted to be the best by the minimum number of criteria. Consider now the set $X' = \{a, d\}$ and the profile $\vec{P}_{X'}$ that looks as

| $P_1$ | $P_2$ | $P_3$ | $P_4$ | $P_5$ | $P_6$ | $P_7$ | $P_8$ |
|---|---|---|---|---|---|---|---|
| $a$ | $a$ | $a$ | $a$ | $a$ | $d$ | $d$ | $d$ |
| $d$ | $d$ | $d$ | $d$ | $d$ | $a$ | $a$ | $a$ |

For this case $n^+(a, \vec{P}_{X'}) = 5$, $n^+(d, \vec{P}_{X'}) = 3$ and $|\vec{P}_{X'}| = 8$. According to the simple majority rule the alternative $a$ will be chosen, i.e., $C(\vec{P}_X, X) = \{a\}$.

## 7. The Borda rule

Put to each $x \in X$ into correspondence a number $r_i(x, \vec{P}_X)$ which is equal to the cardinality of the lower contour set of $x$ in $P_i \in \vec{P}_X$, i.e., $r_i(x, \vec{P}_X) = |\{b \in X : xP_ib\}|$. The sum of that numbers over all $i \in N$ is called the Borda count for the alternative $x$,

$$r(x, \vec{P}_X) = \sum_{i=1}^{n} r_i(x, \vec{P}_X).$$

An alternative with the maximum Borda count is chosen, i.e.,
$$a \in C(\vec{P}_X, X) \Leftrightarrow [\forall b \in X: r(a, \vec{P}_X) \geq r(b, \vec{P}_X)].$$

*Consider the previous example.* $r(a, \vec{P}_X) = 11$, $r(b, \vec{P}_X) = 15$, $r(c, \vec{P}_X) = 12$, $r(d, \vec{P}_X) = 10$. Thus, $C(\vec{P}_X, X) = \{b\}$.

### 8. The Condorcet winner

The Condorcet winner $C(\vec{P}_X, X)$ in the profile $\vec{P}_X$ is an alternative that dominates all other alternatives via majority relation $\mu$ (constructed according to the profile), i.e.,
$$C(\vec{P}_X, X) = [a \in X | \forall x \in X \setminus \{a\}, a\mu x].$$

*Consider the previous example.* Since $a\mu d, b\mu a, b\mu c, c\mu a, c\mu d$, i.e., there is no alternative that dominates all other alternatives via majority relation $\mu$, the Condorcet winner does not exist. Thus, $C(\vec{P}_X, X) = \emptyset$.

### 9. The Black procedure

If the Condorcet winner exists, it is to be chosen. Otherwise, the Borda rule is applied.

*Consider the previous example.* Since there is no Condorcet winner, the Borda rule is applied. As it was shown before, the alternative $b$ will be chosen, i.e., $C(\vec{P}_X, X) = \{b\}$.

### 10. The inverse Borda procedure

The inverse Borda procedure was proposed by J.M. Baldwin [4], so it is also known as the Baldwin procedure.

For each alternative the Borda count is calculated and the alternative with the minimum Borda count is omitted. Then Borda counts are re-calculated for the profile $\vec{P}_{X'}$, $X' = X \setminus \{x\}$, and the procedure is repeated until the choice is found.

*Consider the previous example.* $r(a, \vec{P}_X) = 11$, $r(b, \vec{P}_X) = 15$, $r(c, \vec{P}_X) = 12$, $r(d, \vec{P}_X) = 10$. According to the inverse Borda procedure the alternative $d$ is omitted on the first stage. Consider now the set $X' = \{a, b, c\}$ and the profile $\vec{P}_{X'}$ that looks as

| $P_1$ | $P_2$ | $P_3$ | $P_4$ | $P_5$ | $P_6$ | $P_7$ | $P_8$ |
|---|---|---|---|---|---|---|---|
| a | a | a | b | b | c | b | b |
| b | b | c | c | c | b | c | c |
| c | c | b | a | a | a | a | a |

Let us calculate the Borda count for each alternative $x \in X'$: $r(a, \vec{P}_{X'}) = 6$, $r(b, \vec{P}_{X'}) = 11$, $r(c, \vec{P}_{X'}) = 7$. According to the inverse Borda procedure the alternative $a$ is omitted on the second stage. Consider now the set $X'' = \{b, c\}$ and the profile $\vec{P}_{X''}$ that looks as

| $P_1$ | $P_2$ | $P_3$ | $P_4$ | $P_5$ | $P_6$ | $P_7$ | $P_8$ |
|---|---|---|---|---|---|---|---|
| b | b | c | b | b | c | b | b |
| c | c | b | c | c | b | c | c |

The Borda count for each alternative $x \in X''$ is the following: $r(b, \vec{P}_{X''}) = 6, r(c, \vec{P}_{X''}) = 2$. According to the inverse Borda procedure the alternative $c$ is omitted and the alternative $b$ will be chosen. Thus, $C(\vec{P}_X, X) = \{b\}$.

**11. The Nanson procedure**

For each alternative the Borda count is calculated. Then the average Borda count is calculated,

$$\bar{r}(\vec{P}_X, X) = \frac{\sum_{x \in X} r(x, \vec{P}_X)}{|X|},$$

and alternatives are omitted for which the Borda count is less or equal than the mean value. Then the set $X' = X \setminus \{Y\}$ is considered, where $Y \subset X$ is a set of alternatives with the Borda count that is less or equal than the average value, and the procedure is applied to the profile $\vec{P}_{X'}$. Such procedure is repeated until all remaining alternatives will have the same Borda count.

*Consider the previous example.* $r(a, \vec{P}_X) = 11, r(b, \vec{P}_X) = 15, r(c, \vec{P}_X) = 13, r(d, \vec{P}_X) = 10, \bar{r}(\vec{P}_X, X) = 12$. According to the Nanson procedure the alternatives $a$ and $d$ are omitted. Consider now the set $X' = \{b, c\}$ and the profile $\vec{P}_{X'}$ that looks as

| P₁ | P₂ | P₃ | P₄ | P₅ | P₆ | P₇ | P₈ |
|----|----|----|----|----|----|----|----|
| b  | b  | c  | b  | b  | c  | b  | b  |
| c  | c  | b  | c  | c  | b  | c  | c  |

Let us calculate the Borda count for each alternative $x \in X'$: $r(b, \vec{P}_{X'}) = 6, r(c, \vec{P}_{X'}) = 2$, $\bar{r}(\vec{P}_{X'}, X') = 4$. According to the Nanson procedure the alternative $c$ is omitted and, consequently, the alternative $b$ will be chosen. Thus, $C(\vec{P}_X, X) = \{b\}$.

**12. The Coombs procedure**

First, the simple majority rule is used. If a winner exists, the procedure stops. Otherwise, the alternative $x$ that have been admitted to be the worst by the maximum number of criteria is omitted. Then the procedure again applied to the narrowed set $X' = X \setminus \{x\}$ and the profile $\vec{P}_{X'}$. The procedure stops when an alternative that have been admitted to be the best by 50% + 1 criteria is found.

It is necessary to note the difference between the Coombs procedure and the Hare procedure. In the Coombs procedure we eliminate alternatives that have been admitted to be the worst by the maximum number of criteria while in the Hare procedure we eliminate alternatives that have been admitted to be the best by the minimum number of criteria.

*Consider the previous example.* As it was mentioned before, according to the simple majority rule the choice is empty. In that case, we eliminate alternatives that have been admitted to be the worst by the maximum number of criteria. Since $n^-(a, \vec{P}_X) = 3, n^-(b, \vec{P}_X) = 1$, $n^-(c, \vec{P}_X) = 0, n^-(d, \vec{P}_X) = 4$, the alternative $d$ is omitted on the first stage of the rule. Consider now the set $X' = \{a, b, c\}$ and the profile $\vec{P}_{X'}$ that looks as

|     | $P_1$ | $P_2$ | $P_3$ | $P_4$ | $P_5$ | $P_6$ | $P_7$ | $P_8$ |
|-----|-------|-------|-------|-------|-------|-------|-------|-------|
|     | a     | a     | a     | b     | b     | c     | b     | b     |
|     | b     | b     | c     | c     | c     | b     | c     | c     |
|     | c     | c     | b     | a     | a     | a     | a     | a     |

According to the simple majority rule the choice is empty. Since $n^-(a, \vec{P}_{X'}) = 5$, $n^-(b, \vec{P}_{X'}) = 1$, $n^-(c, \vec{P}_{X'}) = 2$, the alternative $a$ is omitted on the second stage of the rule. Finally, let us consider the set $X'' = \{b, c\}$ and the profile $\vec{P}_{X''}$ that looks as

|     | $P_1$ | $P_2$ | $P_3$ | $P_4$ | $P_5$ | $P_6$ | $P_7$ | $P_8$ |
|-----|-------|-------|-------|-------|-------|-------|-------|-------|
|     | b     | b     | c     | b     | b     | c     | b     | b     |
|     | c     | c     | b     | c     | c     | b     | c     | c     |

The simple majority rule is applied for this case. Since $n^+(b, \vec{P}_{X''}) = 6$, $n^+(c, \vec{P}_{X''}) = 2$, $|\vec{P}_{X''}| = 8$, the alternative $b$ will be chosen ($n^+(b, \vec{P}_{X''}) > \frac{1}{2} \cdot |\vec{P}_{X''}|$). Thus, $C(\vec{P}_X, X) = \{b\}$.

### 13. Minimal dominant set

A set $Q$ is called dominant one if each alternative in $Q$ dominates each alternative outside $Q$ via majority relation $\mu$. Otherwise speaking, $\forall x \in X$
$$x \in Q \Leftrightarrow [\forall y \in X \backslash Q: x\mu y].$$

Then a dominant set $Q$ is called the minimal one if none of its proper subsets is a dominant set. The choice is defined as $C(\vec{P}_X, X) = Q$. If such set is not unique, then the social choice is defined as the union of these sets.

*Consider the previous example.* A majority graph for the profile $\vec{P}_X$ is the following

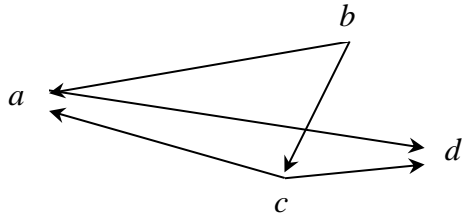

The alternatives $a,b,c,d$ make the minimal dominant set. Thus, $C(\vec{P}_X, X) = Q = X = \{a, b, c, d\}$.

### 14. Minimal undominated set

A set $Q$ is called undominated one if no alternative outside $Q$ dominates some alternative in $Q$ via majority relation $\mu$, i.e.,
$$x \in Q \Leftrightarrow [\nexists y \in X \backslash Q: y\mu x].$$

Undominated set $Q$ is called the minimal one if none of its proper subsets is undominated set. The choice is defined as $C(\vec{P}_X, X) = Q$. If such set is not unique, then the choice is defined as the union of these sets.

*For the previous example* the alternative $b$ makes the minimal undominated set. Thus, $C(\vec{P}_X, X) = Q = \{b\}$.

### 15. Minimal weakly stable set

A set $Q \subseteq A$ is called weakly stable if it has the following property: $\forall x \in X$
$$x \in Q \Leftrightarrow [\exists y \in X \backslash Q: y\mu x \Rightarrow \exists z \in Q: z\mu y],$$
i.e., $x$ belongs to $Q$ if there exists an alternative $y \in X \backslash Q$ which dominates $x$ via majority relation $\mu$, then there exists an alternative $z \in Q$, which dominates $y$, i.e., $z\mu y$. A set $Q \subseteq A$ is called the minimal weakly stable one if none of its proper subsets is a weakly stable set. The choice is defined as $C(\vec{P}_X, X) = Q$. If such set is not unique, then the choice is defined as the union of these sets.

*For the previous example* the alternative $b$ makes the minimal weakly stable set. Thus, $C(\vec{P}_X, X) = Q = \{b\}$.

### 16. The Fishburn rule

Construct the upper contour set $D(x, \vec{P}_X)$ of $x$ in the relation $\mu$, and the binary relation $\gamma$ as follows
$$x\gamma y \Leftrightarrow D(x, \vec{P}_X) \subset D(y, \vec{P}_X).$$
Note that the binary relation $\gamma$ is a partial order.
Then undominated alternatives on $\gamma$ are chosen, i.e.,
$$x \in C(\vec{P}_X, X) \Leftrightarrow [\not\exists y \in X \; y\gamma x].$$

*For the previous example* the upper contour sets for each alternative are $D(a, \vec{P}_X) = \{b, c\}, D(b, \vec{P}_X) = \emptyset, \; D(c, \vec{P}_X) = \{b\}, \; D(d, \vec{P}_X) = \{a, c\}$. Then $b\gamma a$, $b\gamma c$, $b\gamma d$, $c\gamma a$, and $C(\vec{P}_X, X) = \{b\}$.

### 17. Uncovered set I

Construct the lower contour set $L(x, \vec{P}_X)$ of $x$ in the relation $\mu$, and the binary relation $\delta$ as follows
$$x\delta y \Leftrightarrow L(x, \vec{P}_X) \supset L(y, \vec{P}_X).$$
Then undominated alternatives on $\delta$ are chosen, i.e.,
$$x \in C(\vec{P}_X, X) \Leftrightarrow [\not\exists y \in X \; y\delta x].$$

*For the example above* the lower contour sets for each alternative are $L(a, \vec{P}_X) = \{d\}, L(b, \vec{P}_X) = \{a, c\}, \; L(c, \vec{P}_X) = \{a, d\}, \; L(d, \vec{P}_X) = \emptyset$. Then $a\delta d$, $b\delta d$, $c\delta d$, $c\delta a$. Consequently, $C(\vec{P}_X, X) = \{b, c\}$.

### 18. Uncovered set II

An alternative $x$ is said to $B$-dominate an alternative $y$ (denoted as $xBy$) if $x\mu y$ and $D(x, \vec{P}_X) \subset D(y, \vec{P}_X)$, where $D(x, \vec{P}_X)$ is the upper contour set of $x$ in $\mu$. The choice consists of $B$-undominated alternatives, i.e.,
$$x \in C(\vec{P}_X, X) \Leftrightarrow [\not\exists y \in X \; yBx].$$

*For the previous example* the upper contour sets for each alternative are $D(a, \vec{P}_X) = \{b, c\}, D(b, \vec{P}_X) = \emptyset, D(c, \vec{P}_X) = \{b\}, D(d, \vec{P}_X) = \{a, c\}$. Then $bBa$, $bBc$, $cBa$. Consequently, $C(\vec{P}_X, X) = \{b, d\}$.

### 19. The Richelson rule

Construct lower and upper contour sets $D(x, \vec{P}_X)$ and $L(x, \vec{P}_X)$ for each $x \in X$ in the relation $\mu$, and the binary relation $\sigma$ as follows

$$x\sigma y \Leftrightarrow \begin{cases} D(x, \vec{P}_X) \subseteq D(y, \vec{P}_X) \text{ and} \\ L(x, \vec{P}_X) \supseteq L(y, \vec{P}_X) \text{ and} \\ D(x, \vec{P}_X) \subset D(y, \vec{P}_X) \text{ or } L(x, \vec{P}_X) \supset L(y, \vec{P}_X). \end{cases}$$

Then undominated alternatives on $\sigma$ are chosen, i.e.,

$$x \in C(\vec{P}_X, X) \Leftrightarrow [\nexists y \in X \ y\sigma x].$$

*For the example above* the upper and the lower contour sets for each alternative are $D(a, \vec{P}_X) = \{b, c\}, D(b, \vec{P}_X) = \emptyset, D(c, \vec{P}_X) = \{b\}, D(d, \vec{P}_X) = \{a, c\}, L(a, \vec{P}_X) = \{d\}, L(b, \vec{P}_X) = \{a, c\}, L(c, \vec{P}_X) = \{a, d\}, L(d, \vec{P}_X) = \emptyset$. Then $b\sigma d$, $c\sigma a$. Consequently, $C(\vec{P}_X, X) = \{b, c\}$.

### 20. The Core

The choice includes alternatives, which are undominated via majority relation $\mu$, i.e.,

$$C(\vec{P}_X, X) = [a \in X | \nexists x \in X, x\mu a].$$

Obviously, if each alternative is dominated via majority relation $\mu$, the core is empty.

*For the previous example* the alternative $b$ makes the core. Thus, $C(\vec{P}_X, X) = \{b\}$.

### 21. k-stable sets

An alternative $x$ belongs to $Q$ if there exists an alternative $y \in X \backslash Q$ which dominates $x$ via majority relation $\mu$, then there exists a $\mu$-path of length $s$, $s \leq k$, to $y$ from some alternative $z \in Q$. A set $Q \subseteq A$ is called the minimal k-stable one if none of its proper subsets is a k-stable set. The choice is defined as $C(\vec{P}_X, X) = Q$. If such set is not unique, then the choice is defined as the union of these sets. It follows from this definition that a weakly stable set is a 1-stable set.

*For the previous example* let us define k-stable sets for the set $X$ depending on the parameter value $k$.

| k | $C(\vec{P}_X, X, k)$ |
|---|---|
| k=1 | {b} |
| k=2 | {b} |
| k=3 | {b} |

### 22. The threshold rule

Let $v_1(x, \vec{P}_X)$ be the number of criteria for which the alternative $x$ is the worst in their ordering, $v_2(x, \vec{P}_X)$ is the number of criteria placing $x$ the second worst, and so on, $v_m(x, \vec{P}_X)$ is the number of criteria considering the alternative $x$ as their best one. Then we order the alternatives lexicographically. The alternative $x$ is said to $V$-dominate the alternative $y$ if $v_1(x, \vec{P}_X) < v_1(y, \vec{P}_X)$ or, if there exists $k \leq m$, such that $v_i(x, \vec{P}_X) = v_i(y, \vec{P}_X)$, $i = 1, ..., k-1$, and $v_k(x, \vec{P}_X) < v_k(y, \vec{P}_X)$. In other words, first, the number of worst places are compared, if these numbers are equal than the number of second worst places are compared, and so on. The alternatives which are not dominated by other alternatives via $V$ are chosen.

*For the previous example* let us calculate the value $v_i(x, \vec{P}_X)$ for each alternative $x \in X$ depending on the value $i$.

| $i$ | $v_i(a, \vec{P}_X)$ | $v_i(b, \vec{P}_X)$ | $v_i(c, \vec{P}_X)$ | $v_i(d, \vec{P}_X)$ |
|---|---|---|---|---|
| $i=1$ | 3 | 1 | 0 | 4 |
| $i=2$ | 2 | 1 | 4 | 1 |
| $i=3$ | 0 | 4 | 4 | 0 |
| $i=4$ | 3 | 2 | 0 | 3 |

According to the threshold rule the alternative $c$ is the best one as it has the fewest number of the worst places, i.e., $C(\vec{P}_X, X) = c$. As for other alternatives, the alternative $b$ is better than the alternative $a$ which is better than the alternative $d$.

### 23. The Copeland rule 1

Construct a function $u(x, \vec{P}_X)$, which is equal to the difference of cardinalities of lower and upper contour sets of the alternative $x \in X$ via majority relation $\mu$, i.e.,
$$u(x, \vec{P}_X) = |L(x, \vec{P}_X)| - |D(x, \vec{P}_X)|.$$
Then the choice is defined by maximization of $u(x, \vec{P}_X)$, that is
$$x \in C(\vec{P}_X, X) \Leftrightarrow [\nexists y \in X: u(y, \vec{P}_X) > u(x, \vec{P}_X)].$$

*Let us calculate for the previous example* the function value $u(x, \vec{P}_X)$ for each alernative $x \in X$. Then $u(a, \vec{P}_X) = -1$, $u(b, \vec{P}_X) = 2$, $u(c, \vec{P}_X) = 1$, $u(d, \vec{P}_X) = -2$.

According to the Copeland rule 1 the alternative $b$ will be chosen, i.e., $C(\vec{P}_X, X) = \{b\}$.

### 24. The Copeland rule 2

The function $u(x, \vec{P}_X)$ is defined by cardinality of lower contour set of the alternative $x \in X$ via majority relation $\mu$. The choice is defined by maximization of $u(x, \vec{P}_X)$.

*Let us define for the previous example* the lower contour sets for each alternative $x \in X$. $L(a, \vec{P}_X) = \{d\}$, $L(b, \vec{P}_X) = \{a, c\}$, $L(c, \vec{P}_X) = \{a, d\}$, $L(d, \vec{P}_X) = \emptyset$. According to the Copeland rule 2 the alternatives $b$ and $c$ will be chosen, i.e., $C(\vec{P}_X, X) = \{b, c\}$.

### 25. The Copeland rule 3

The function $u(x, \vec{P}_X)$ is defined by cardinality of upper contour set of the alternative $x \in X$ via majority relation $\mu$. The choice is defined by minimization of $u(x, \vec{P}_X)$.

*For the previous example* the upper contour sets for each alternative $x \in X$ are $D(a, \vec{P}_X) = \{b, c\}$, $D(b, \vec{P}_X) = \emptyset$, $D(c, \vec{P}_X) = \{b\}$, $D(d, \vec{P}_X) = \{a, c\}$. According to the Copeland rule 3 the alternative $b$ will be chosen, i.e., $C(\vec{P}_X, X) = \{b\}$.

## 26. The super-threshold choice rule

Let the criterion $\varphi(x)$ be defined over the set $A$, and the threshold function $V(X)$ over $2^A$ assigning to each subset $X \in 2^A$ a threshold level $V(X)$, $V(\cdot): 2^A \to \mathbb{R}^1$. Define the super-threshold rule as follows

$$y \in C(\vec{P}_X, X) \Leftrightarrow (y \in X: \varphi(y) \geq V(X)).$$

*For the previous example* let us write the profile $\vec{P}_X$ in the following form

|   | $\varphi_1$ | $\varphi_2$ | $\varphi_3$ | $\varphi_4$ | $\varphi_5$ | $\varphi_6$ | $\varphi_7$ | $\varphi_8$ |
|---|---|---|---|---|---|---|---|---|
| a | 4 | 4 | 4 | 2 | 2 | 1 | 1 | 1 |
| b | 3 | 3 | 1 | 4 | 4 | 2 | 3 | 3 |
| c | 2 | 2 | 3 | 3 | 3 | 3 | 2 | 2 |
| d | 1 | 1 | 2 | 1 | 1 | 4 | 4 | 4 |

Suppose that the choice is made on the second criterion and the threshold value is equal to 3. According to the super-threshold choice rule the alternatives $a$ and $b$ will be chosen, i.e., $C(\vec{P}_X, X) = \{a, b\}$.

## 27. The minmax procedure.

Construct the matrix $S^-(\vec{P}_X, X)$, such that

$$\forall a, b \in X, S^-(\vec{P}_X, X) = (n(a, b, \vec{P}_X)),$$

where $n(a, b, \vec{P}_X) = card\{i \in N \mid aP_i b\}$, $n(a, a, \vec{P}_X) = -\infty$.

The choice is defined as

$$x \in C(\vec{P}_X, X) \Leftrightarrow x \in \arg\min_{a \in X} \max_{b \in X}\{n(b, a, \vec{P}_X)\}.$$

*For the previous example* the matrix $S^-(\vec{P}_X, X)$ is the following

$$S^-(\vec{P}_X, X) = \begin{array}{c|cccc} & a & b & c & d \\ \hline a & -\infty & 3 & 3 & 5 \\ b & 5 & -\infty & 6 & 4 \\ c & 5 & 2 & -\infty & 5 \\ d & 3 & 4 & 3 & -\infty \end{array}$$

$\max_{b \in X}\{n(b, a, \vec{P}_X)\}$   5   **4**   6   5

According to the minmax procedure the alternative $b$ will be chosen.

## 28. The Simpson procedure

Construct the matrix $S^+(\vec{P}_X, X)$, such that

$$\forall a, b \in X, S^+(\vec{P}_X, X) = (n(a, b, \vec{P}_X)),$$

where $n(a, b, \vec{P}_X) = card\{i \in N \mid aP_i b\}$, $n(a, a, \vec{P}_X) = +\infty$.

The choice is defined as

$$x \in C(\vec{P}_X, X) \Leftrightarrow x \in \arg\max_{a \in X} \min_{b \in X}\{n(a, b, \vec{P}_X)\}.$$

*For the previous example* the matrix $S^+(\vec{P}_X, X)$ is the following

$$S^+(\vec{P}_X, X) = \begin{array}{c|cccc} & a & b & c & d \\ \hline a & +\infty & 3 & 3 & 5 \\ b & 5 & +\infty & 6 & 4 \\ c & 5 & 2 & +\infty & 5 \\ d & 3 & 4 & 3 & +\infty \end{array} \quad \begin{array}{c} \min_{b \in X}\{n(a,b,\vec{P}_X)\} \\ 3 \\ \mathbf{4} \\ 2 \\ 3 \end{array}$$

According to the Simpson procedure the alternative *b* will be chosen.

Now let us discuss normative conditions.

## Normative conditions

All normative conditions, which characterize different choice procedures, can be divided into the following groups [1, 9, 12-14]
1. Rationality conditions;
2. Monotonicity conditions;
3. The non-compensatory condition.

**Rationality conditions**

There are four rationality conditions for choice functions.
1. Heredity condition (**H**)
$$\forall X, X' \in 2^A, X' \subseteq X \Rightarrow C(\vec{P}_{X'}, X') \supseteq C(\vec{P}_X, X) \cap X'.$$

If the presented set is narrowed by eliminating some alternatives, those chosen from the initial set and remaining in the narrowed set will be chosen from the narrowed set.

*Example.* Let $X = \{a, b, c\}$ and the profile $\vec{P}_X$ looks as

| $P_1$ | $P_2$ | $P_3$ | $P_4$ | $P_5$ |
|---|---|---|---|---|
| a | a | c | a | b |
| b | b | b | b | a |
| c | c | a | c | c |

Let for the choice the inverse plurality rule is applied. According to this rule the alternative *b* will be chosen, i.e., $C(\vec{P}_X, X) = \{b\}$.

Consider now the subset $X' = X \setminus \{c\}$. A contraction of the profile $\vec{P}_X$ onto a set $X'$, i.e., $\vec{P}_{X'}$, looks as

| $P_1$ | $P_2$ | $P_3$ | $P_4$ | $P_5$ |
|---|---|---|---|---|
| a | a | b | a | b |
| b | b | a | b | a |

According to the inverse plurality rule the alternative *a* will be chosen, i.e., $C(\vec{P}_{X'}, X') = \{a\}$.

Thus, $C(\vec{P}_{X'}, X') \not\supseteq C(\vec{P}_X, X) \cap X'$. Hence, the inverse plurality rule does not satisfy the condition **H**.

2. Concordance condition (**C**)
$$\forall X', X'' \in 2^A \to C(\vec{P}_{X' \cup X''}, X' \cup X'') \supseteq C(\vec{P}_{X'}, X') \cap C(\vec{P}_{X''}, X'').$$

The Concordance condition (**C**) requires that all the alternatives chosen simultaneously from $X'$ and $X''$ be included in the choice when their union $X = X' \cup X''$ is presented.

*Consider the previous example*. According to the inverse plurality rule the alternative $b$ will be chosen, i.e., $C(\vec{P}_X, X) = \{b\}$. Consider the subset $X' = X\setminus\{c\}$. A contraction of the profile $\vec{P}_X$ onto a set $X'$, i.e., $\vec{P}_{X'}$, looks as

| P₁ | P₂ | P₃ | P₄ | P₅ |
|---|---|---|---|---|
| a | a | b | a | b |
| b | b | a | b | a |

According to the rule the alternative $a$ will be chosen, i.e., $C(\vec{P}_{X'}, X') = \{a\}$.

Finally, let us consider the subset $X'' = X\setminus\{b\}$. A contraction of the profile $\vec{P}_X$ onto a set $X''$, i.e., $\vec{P}_{X''}$, looks as

| P₁ | P₂ | P₃ | P₄ | P₅ |
|---|---|---|---|---|
| a | a | c | a | a |
| c | c | a | c | c |

According to the rule the alternative $a$ will be chosen, i.e., $C(\vec{P}_{X''}, X'') = \{a\}$.

Then $C(\vec{P}_{X'}, X') \cap C(\vec{P}_{X''}, X'') = \{a\} \nsubseteq C(\vec{P}_X, X)$. Thus, the inverse plurality rule does not satisfy the condition **C**.

3. Outcast condition (**O**)
$$\forall X, X' \in 2^A, X' \subseteq X\setminus C(\vec{P}_X, X) \Rightarrow C(\vec{P}_{X\setminus X'}, X\setminus X') = C(\vec{P}_X, X).$$

Choice functions for which narrowing of $X$ by rejection of some or all the alternatives not chosen from the initial set $X$ does not change the choice, satisfy the condition **O**.

*Consider the previous example*. According to the inverse plurality rule the alternative $b$ will be chosen, i.e., $C(\vec{P}_X, X) = \{b\}$. Now let us consider the set $X' = \{c\}$. A contraction of the profile $\vec{P}_X$ onto a set $X\setminus X'$, i.e., $\vec{P}_{X\setminus X'}$, looks as

| P₁ | P₂ | P₃ | P₄ | P₅ |
|---|---|---|---|---|
| a | a | b | a | b |
| b | b | a | b | a |

According to the rule the alternative $a$ will be chosen, i.e. $C(\vec{P}_{X\setminus X'}, X\setminus X') = \{a\}$.

Then $C(\vec{P}_{X\setminus X'}, X\setminus X') \neq C(\vec{P}_X, X)$. Thus, the inverse plurality rule does not satisfy the condition **O**.

4. Arrow's choice axiom (**ACA**)
$$\forall X, X' \in 2^A, X' \subseteq X \Rightarrow \begin{cases} \text{if } C(\vec{P}_X, X) = \emptyset, \text{then } C(\vec{P}_{X'}, X') = \emptyset, \\ \text{if } C(\vec{P}_X, X) \cap X' \neq \emptyset, \text{then } C(\vec{P}_{X'}, X') = C(\vec{P}_X, X) \cap X'. \end{cases}$$

For the case where the choice is not empty, the condition **ACA** requires that the alternatives chosen from the initial set $X$ and left in the narrowed one $X'$ and only such alternatives be chosen from $X'$.

*Consider the previous example*. The condition **ACA** is not satisfied for the inverse plurality rule as the condition **H** is not satisfied for this choice procedure.

Classical rationality includes conditions **H**, **C** and **O**, i.e., Heredity condition, Concordance condition and Outcast condition. The condition **ACA** is a stronger version of the conditions **H**, **C**, **O**.

**Monotonicity conditions**
*1. Monotonicity condition 1*

$\forall X \in 2^A, x \in C(\vec{P}_X, X), \forall \vec{P}_X, \vec{P'}_X: (\forall a, b \in X, aP_i b \Leftrightarrow aP'_i b \: \& \: \exists y \in X, yP_i x \Rightarrow xP'_i y) \Rightarrow x \in C(\vec{P'}_X, X)$.

Let the alternative *x* be chosen from initial set *X*. Suppose now that the relative position of the alternative *x* was improved while the relative comparison of any pair of other alternatives remains unchanged. The Monotonicity condition 1 is satisfied if the alternative *x* is still in choice.

*Example*. Let $X = \{a, b, c\}$ and the profile $\vec{P}_X$ looks as

| $P_1$ | $P_2$ | $P_3$ | $P_4$ | $P_5$ |
|---|---|---|---|---|
| a | a | c | a | b |
| b | b | b | b | a |
| c | c | a | c | c |

Let for the choice the inverse plurality rule is applied. According to this rule the alternative *b* will be chosen, i.e., $C(\vec{P}_X, X) = \{b\}$.

Consider now the profile $\vec{P'}_X$, which differs from the profile $\vec{P}_X$ only by improved position of the alternative *b* in $P'_2$:

| $P'_1$ | $P'_2$ | $P'_3$ | $P'_4$ | $P'_5$ |
|---|---|---|---|---|
| a | b | c | a | b |
| b | a | b | b | a |
| c | c | a | c | c |

According to the inverse plurality rule the alternative *b* is chosen, i.e., $C(\vec{P'}_X, X) = \{b\}$.

Then $b \in C(\vec{P}_X, X)$ and $b \in C(\vec{P'}_X, X)$. Thus, the Monotonicity condition 1 is satisfied for this case.

*2. Monotonicity condition 2*

$\forall X \in 2^A, x, y \in C(\vec{P}_X, X), X' = X \setminus \{x\}, X'' = X \setminus \{y\} \rightarrow x \in C(\vec{P}_{X''}, X'') \: \& \: y \in C(\vec{P}_{X'}, X')$.

Let two alternatives *x,y* be chosen from the initial set *X*. The Monotonicity condition 2 is satisfied if one of the chosen alternatives (*x* or *y*) is still in choice when the other chosen alternative (*y* or *x*) is eliminated.

It is necessary to mention that since some choice procedures choose no more than one best alternative, the Monotonicity condition 2 is not applicable to these procedures as it considers the choice of more than two alternatives. In other words, such procedures obey the Monotonicity condition 2 trivially.

*Example.* Let $X = \{a, b, c\}$ and the profile $\vec{P}_X$ looks as

| $P_1$ | $P_2$ | $P_3$ | $P_4$ | $P_5$ |
|---|---|---|---|---|
| b | a | a | a | a |
| c | b | b | b | c |
| a | c | c | c | b |

Let for the choice the inverse plurality rule is applied. According to this rule the alternatives *a* and *b* will be chosen, i.e., $C(\vec{P}_X, X) = \{a, b\}$.

Consider now the set $X' = X \setminus \{b\}$. A contraction of the profile $\vec{P}_X$ onto a set $X'$, i.e., $\vec{P}_{X'}$, looks as

| $P_1$ | $P_2$ | $P_3$ | $P_4$ | $P_5$ |
|---|---|---|---|---|
| c | a | a | a | a |
| a | c | c | c | c |

According to the inverse plurality rule the alternative *a* will be chosen, i.e., $C(\vec{P}_{X'}, X') = \{a\}$.

Finally, consider the set $X'' = X \setminus \{a\}$. A contraction of the profile $\vec{P}_X$ onto a set $X''$, i.e., $\vec{P}_{X''}$, looks as

| $P_1$ | $P_2$ | $P_3$ | $P_4$ | $P_5$ |
|---|---|---|---|---|
| b | b | b | b | c |
| c | c | c | c | b |

According to the rule the alternative *b* will be chosen, i.e., $C(\vec{P}_{X''}, X'') = \{b\}$.

Then $\{a, b\} \in C(\vec{P}_X, X)$, $\{a\} \in C(\vec{P}_{X'}, X')$ and $\{b\} \in C(\vec{P}_{X''}, X'')$. Thus, the Monotonicity condition 2 is satisfied for this example.

### 3. Strict monotonicity condition

$$\forall X \in 2^A, \forall y \in X, y \notin C(\vec{P}_X, X), \forall \vec{P}_X, \vec{P}'_X : (\forall a, b \in X, aP_i b \Leftrightarrow aP'_i b \ \& aP_i y \Rightarrow yP'_i a)$$

$$\to C(\vec{P}'_X, X) = \begin{bmatrix} C(\vec{P}_X, X) \text{ or} \\ \{y\} \text{ or} \\ C(\vec{P}_X, X) \cup \{y\}. \end{bmatrix}$$

The change of the relative position of unchosen alternative $y \in X \setminus C(\vec{P}_X, X)$ so that its position will be improved while the relative comparison of any pair of other alternatives remains unchanged leads to the choice of the alternative *y* or/and all alternatives that were in the initial choice $C(\vec{P}_X, X)$.

*Example.* Let $X = \{a, b, c\}$ and the profile $\vec{P}_X$ looks as

| $P_1$ | $P_2$ | $P_3$ | $P_4$ |
|---|---|---|---|
| a | a | b | b |
| c | b | a | a |
| b | c | c | c |

Let for the choice the Borda rule is applied. The Borda counts are $r(a, \vec{P}_X) = 6$, $r(b, \vec{P}_X) = 5$, $r(c, \vec{P}_X) = 1$. According to this rule the alternative *a* will be chosen, i.e., $C(\vec{P}_X, X) = \{a\}$.

Consider now the profile $\vec{P_X'}$, which differs from the profile $\vec{P_X}$ only by improved position of the alternative $c$ in $P_1'$:

|  $P_1'$ | $P_2'$ | $P_3'$ | $P_4'$ |
|---|---|---|---|
| c | a | b | b |
| a | b | a | a |
| b | c | c | c |

The Borda counts are now $r\left(a, \vec{P_X'}\right) = r\left(b, \vec{P_X'}\right) = 5, r\left(c, \vec{P_X'}\right) = 2$. According to the Borda rule the alternatives $a$ and $b$ will be chosen, i.e., $C(\vec{P_X'}, X) = \{a, b\}$.

$$C\left(\vec{P_X'}, X\right) \neq \begin{bmatrix} C(\vec{P_X}, X) \text{ or} \\ \{c\} \text{ or} \\ C(\vec{P_X}, X) \cup \{c\}. \end{bmatrix}$$

Thus, the Borda rule does not satisfy the strict monotonicity condition.

**The non-compensatory condition**

Consider two alternatives $x, y \in X$ that characterized by a set of criterial estimates $\varphi_1(x)$, $\varphi_2(x), \ldots, \varphi_n(x)$ and $\varphi_1(y), \varphi_2(y), \ldots, \varphi_n(y)$, which may take $m$ different values ($n$ is the number of criteria). Denote by $v_j(x)$ the number of criteria for which the alternative $x$ take value $j$, where $j \in [1, m]$.

$$v_j(x) = |\{i \in [1, n]: \varphi_i(x) = j\}|.$$

The non-compensatory condition can be formulated as

$$\varphi(x) > \varphi(y) \Leftrightarrow \begin{cases} v_1(x) < v_1(y) \text{ if } m = 2, \\ \exists k \in [2, m-1]: \forall j \in [1, k-1] \ v_j(x) = v_j(y) \text{ and } v_k(x) < v_k(y) \text{ if } m \geq 3. \end{cases}$$

In other words, low estimates on one criterion cannot be compensated by high estimates on other criteria.

*Example.* Let $X = \{a, b, c\}$ and the profile $\vec{P_X}$ looks as

| $P_1$ | $P_2$ | $P_3$ | $P_4$ |
|---|---|---|---|
| c | a | b | b |
| a | b | a | a |
| b | c | c | c |

Let for the choice the Borda rule is applied. The Borda counts are $r(a, \vec{P_X}) = r(b, \vec{P_X}) = 5, r(c, \vec{P_X}) = 2$. According to this rule the alternatives $a$ and $b$ will be chosen, i.e., $C(\vec{P_X}, X) = \{a, b\}$.

Now we write the profile $\vec{P_X}$ in the following form

| X | $\varphi_1$ | $\varphi_2$ | $\varphi_3$ | $\varphi_4$ |
|---|---|---|---|---|
| a | 2 | 3 | 2 | 2 |
| b | 1 | 2 | 3 | 3 |
| c | 3 | 1 | 1 | 1 |

According to the non-compensatory condition the alternative $a$ is better than the alternative $b$ while the alternative $b$ is better than the alternative $c$. Thus, the non-compensatory condition is not satisfied as $b \in C(\vec{P_X}, X)$.

# A study of the properties of choice procedures

A study of the properties is conducted as follows. If a choice procedure does not satisfy given normative condition, a counter-example is provided. On the country, if a choice procedure satisfies given normative condition a necessary proof is followed.

The results of the study of the properties are given in Theorem 1.

**Theorem 1.** Information on which choice procedures satisfy given normative conditions is provided in Table 2.

**Table 2. Properties of existing choice procedures («+» - choice procedure satisfies given normative condition, «-» - choice procedure does not satisfy given normative condition)**

| № | Choice procedure | Normative conditions | | | | | | | |
|---|---|---|---|---|---|---|---|---|---|
| | | Rationality conditions | | | | Monotonicity conditions | | | Non-compensatory condition |
| | | Heredity condition (H) | Concordance condition (C) | Outcast condition (O) | Arrow's choice axiom (ACA) | Monotonicity condition 1 | Monotonicity condition 2 | Strict monotonicity condition | |
| 1 | The simple majority rule | + | - | + | - | + | + | - | - |
| 2 | The plurality rule | - | - | - | - | + | - | - | - |
| 3 | The inverse plurality rule | - | - | - | - | + | - | - | - |
| 4 | The q-Approval rule | - | - | - | - | + | - | - | - |
| 5 | The run-off procedure | - | - | - | - | - | + | - | - |
| 6 | The Hare rule (the Ware procedure) | - | - | - | - | - | + | - | - |
| 7 | The Borda rule | - | - | - | - | + | - | - | - |
| 8 | The Black procedure | - | - | - | - | + | - | - | - |
| 9 | The inverse Borda rule | - | - | - | - | - | - | - | - |
| 10 | The Nanson rule | - | - | - | - | - | - | - | - |
| 11 | The Coombs procedure | - | - | - | - | - | + | - | - |
| 12 | Minimal dominant set | - | + | + | - | + | - | - | - |
| 13 | Minimal undominant set | - | - | - | - | + | - | - | - |
| 14 | Minimal weakly stable set | - | - | - | - | + | - | - | - |
| 15 | The Fishburn rule | - | - | - | - | + | - | - | - |
| 16 | Uncovered set I | - | - | - | - | + | - | - | - |
| 17 | Uncovered set II | - | + | - | - | + | - | - | - |
| 18 | The Richelson rule | - | - | - | - | + | - | - | - |
| 19 | The Condorcet winner | + | + | - | - | + | + | - | - |
| 20 | The core | + | + | - | - | + | + | - | - |
| 21 | k-stable set (k>1) | - | - | - | - | + | - | - | - |
| 22 | The threshold rule | - | - | - | - | + | - | - | + |
| 23 | The Copeland rule 1 | - | - | - | - | + | - | - | - |
| 24 | The Copeland rule 2 | - | - | - | - | + | - | - | - |
| 25 | The Copeland rule 3 | - | - | - | - | + | - | - | - |
| 26a | The super-threshold rule (fixed threshold) | + | + | + | + | + | + | + | - |
| 26b | The super-threshold rule (threshold depends on X) | - | - | - | - | + | + | - | - |
| 27 | The minimax procedure | - | - | - | - | + | - | - | - |
| 28 | The Simpson procedure | - | - | - | - | + | - | - | - |

The proof of the theorem is provided in Appendix 1.

## Conclusion

We have studied the properties of 28 existing choice procedures, which can be used in various social and multi-criteria problems. It was defined which choice procedures satisfy given normative conditions, showing how a final choice is changed due to the changes of preferences or a set of feasible alternatives. Such information leads to a better understanding of different choice procedures and how stable and sensible is a set of alternatives obtained after applying some choice procedure.

The results show that only the simple majority rule, the Condorcet winner, the core and the threshold rule with fixed threshold level satisfy the condition **H**, only minimal dominant set, uncovered set II, the Condorcet winner, the core and the threshold rule with fixed threshold level satisfy the condition **C**, and only the simple majority rule, minimal dominant set and the threshold rule with fixed threshold level satisfy the condition **O**. More detailed information is provided in Table 2.

# Appendix 1. Properties of existing choice procedures

## 1. The simple majority rule

### 1.1 Heredity condition (H)

Let $C(\vec{P}_X, X) = \{a\}$, i.e., $n^+(a, \vec{P}_X) > \frac{1}{2} \cdot |\vec{P}_X|$. For a contraction of the profile $\vec{P}_X$ onto a set $X' \subseteq X$ ($a \in X'$) it is true that $n^+(a, \vec{P}_{X'}) \geq n^+(a, \vec{P}_X)$ and $|\vec{P}_X| = |\vec{P}_{X'}|$. Then $n^+(a, \vec{P}_{X'}) > \frac{1}{2} \cdot |\vec{P}_{X'}|$ and $C(\vec{P}_{X'}, X') = \{a\}$.

Thus, $C(\vec{P}_X, X) \cap X' = \{a\} = C(\vec{P}_{X'}, X')$. The simple majority rule satisfies the condition **H**.

### 1.2 Concordance condition (C)

Let $X = \{a, b, c\}$ and the profile $\vec{P}_X$ is the following

| $P_1$ | $P_2$ | $P_3$ |
|---|---|---|
| a | b | c |
| b | a | b |
| c | c | a |

According to the simple majority rule the choice is empty, i.e., $C(\vec{P}_X, X) = \emptyset$.

Now let us consider the subset $X' = X \setminus \{c\}$. A contraction of the profile $\vec{P}_X$ onto a set $X'$, i.e., $\vec{P}_{X'}$, looks as

| $P_1$ | $P_2$ | $P_3$ |
|---|---|---|
| a | b | b |
| b | a | a |

According to the simple majority rule the alternative $b$ will be chosen, i.e., $C(\vec{P}_{X'}, X') = \{b\}$.

Finally, let us consider the subset $X'' = X \setminus \{a\}$. A contraction of the profile $\vec{P}_X$ onto a set $X''$, i.e., $\vec{P}_{X''}$, looks as

| $P_1$ | $P_2$ | $P_3$ |
|---|---|---|
| b | b | c |
| c | c | b |

According to the rule the alternative $b$ will be chosen, i.e., $C(\vec{P}_{X''}, X'') = \{b\}$.

Then $C(\vec{P}_{X'}, X') \cap C(\vec{P}_{X''}, X'') = \{b\} \nsubseteq C(\vec{P}_X, X)$. The simple majority rule does not satisfy the condition **C**.

### 1.3 Outcast condition (O)

The simple majority rule satisfies the condition **O** (see paragraph 1.1 of this section).

### 1.4 Arrow's choice axiom (ACA)

The simple majority rule does not satisfy the condition **ACA** as the condition **C** is not satisfied.

### 1.5 Monotonicity condition 1

Let $C(\vec{P}_X, X) = \{a\}$, i.e., $n^+(a, \vec{P}_X) > \frac{1}{2} \cdot |\vec{P}_X|$.

Consider now a profile $\overrightarrow{P'_X}$, which differs from the profile $\vec{P}_X$ only by improved position of the alternative $a$. Then, $n^+\left(a, \overrightarrow{P'_X}\right) \geq n^+(a, \vec{P}_X)$ and $|\vec{P}_X| = |\overrightarrow{P'_X}|$. Thus, $n^+\left(a, \overrightarrow{P'_X}\right) > \frac{1}{2} \cdot |\overrightarrow{P'_X}|$ and the alternative $a$ will be chosen, i.e., $C\left(\overrightarrow{P'_X}, X\right) = \{a\}$.

Since $a \in C(\vec{P}_X, X)$ and $a \in C(\overrightarrow{P'_X}, X)$, the simple majority rule satisfies the Monotonicity condition 1.

### 1.6 The Monotonicity condition 2
Since the simple majority rule always chooses no more than one best alternative, the Monotonicity condition 2 is not applicable to it as it considers the choice of more than two alternatives. In other words, the simple majority rule obeys the Monotonicity condition 2 trivially.

### 1.7 The strict monotonicity condition
Let $X = \{a, b, c\}$ and the profile $\vec{P}_X$ is the following

| $P_1$ | $P_2$ | $P_3$ |
|---|---|---|
| a | a | b |
| c | b | c |
| b | c | a |

According to the simple majority rule the alternative $a$ will be chosen, i.e., $C(\vec{P}_X, X) = \{a\}$.
Consider now a profile $\overrightarrow{P'_X}$, which differs from the profile $\vec{P}_X$ only by improved position of the alternative $c$ in $P'_1$:

| $P'_1$ | $P'_2$ | $P'_3$ |
|---|---|---|
| c | a | b |
| a | b | c |
| b | c | a |

According to the rule the choice is empty, i.e., $C(\overrightarrow{P'_X}, X) = \emptyset$.

$$C\left(\overrightarrow{P'_X}, X\right) \neq \begin{bmatrix} C(\vec{P}_X, X) \text{ or} \\ \{c\} \text{ or} \\ C(\vec{P}_X, X) \cup \{c\}. \end{bmatrix}$$

Thus, the simple majority rule does not satisfy the strict monotonicity condition.

### 1.8 The non-compensatory condition
Let $X = \{a, b, c\}$ and the profile $\vec{P}_X$ is the following

| $P_1$ | $P_2$ | $P_3$ |
|---|---|---|
| a | a | b |
| b | b | c |
| c | c | a |

According to the rule the alternative $a$ will be chosen, i.e., $C(\vec{P}_X, X) = \{a\}$.

Let us write the profile $\vec{P}_X$ in the following form

|   | $\varphi_1$ | $\varphi_2$ | $\varphi_3$ |
|---|---|---|---|
| a | 3 | 3 | 1 |
| b | 2 | 2 | 3 |
| c | 1 | 1 | 2 |

According to the non-compensatory condition the alternative $b$ is better than the alternative $a$ while the alternative $a$ is better than the alternative $c$. Thus, the non-compensatory condition is not satisfied as $\{b\} \neq C(\vec{P}_X, X)$.

## 2. The plurality rule

### 2.1 Heredity condition (H)

Let $X = \{a, b, c\}$ and the profile $\vec{P}_X$ is the following

| $P_1$ | $P_2$ | $P_3$ | $P_4$ | $P_5$ |
|---|---|---|---|---|
| a | a | c | b | b |
| c | c | b | a | a |
| b | b | a | c | c |

Let us calculate the value $n^+(x, \vec{P}_X)$ for each alternative $x \in X$: $n^+(a, \vec{P}_X) = 2$, $n^+(b, \vec{P}_X) = 2$, $n^+(c, \vec{P}_X) = 1$. According to the plurality rule the alternatives $a$ and $b$ will be chosen, i.e., $C(\vec{P}_X, X) = \{a, b\}$.

Consider now the subset $X' = X \setminus \{c\}$. A contraction of the profile $\vec{P}_X$ onto a set $X'$, i.e., $\vec{P}_{X'}$, looks as

| $P_1$ | $P_2$ | $P_3$ | $P_4$ | $P_5$ |
|---|---|---|---|---|
| a | a | b | b | b |
| b | b | a | a | a |

According to the rule the alternative $b$ will be chosen, i.e., $C(\vec{P}_{X'}, X') = \{b\}$.

Then $C(\vec{P}_{X'}, X') \not\supseteq C(\vec{P}_X, X) \cap X'$. Thus, the plurality rule does not satisfy the condition **H**.

### 2.2 Concordance condition (C)

Let $X = \{a, b, c, d\}$ and the profile $\vec{P}_X$ is the following

| $P_1$ | $P_2$ | $P_3$ | $P_4$ | $P_5$ | $P_6$ | $P_7$ | $P_8$ | $P_9$ | $P_{10}$ | $P_{11}$ | $P_{12}$ | $P_{13}$ | $P_{14}$ |
|---|---|---|---|---|---|---|---|---|---|---|---|---|---|
| c | c | c | b | b | b | d | d | d | D | a | a | a | a |
| a | a | a | d | d | d | b | b | c | C | b | b | b | c |
| d | d | d | a | a | a | c | c | b | b | c | c | c | b |
| b | b | b | c | c | c | a | a | a | a | d | d | d | d |

Let us calculate the value $n^+(x, \vec{P}_X)$ for each alternative $x \in X$: $n^+(a, \vec{P}_X) = 4$, $n^+(b, \vec{P}_X) = 3$, $n^+(c, \vec{P}_X) = 3$, $n^+(d, \vec{P}_X) = 4$. According to the plurality rule the alternatives $a$ and $d$ will be chosen, i.e., $C(\vec{P}_X, X) = \{a, d\}$.

Consider now the subset $X' = X\setminus\{d\}$. A contraction of the profile $\vec{P}_X$ onto a set $X'$, i.e., $\vec{P}_{X'}$, looks as

| $P_1$ | $P_2$ | $P_3$ | $P_4$ | $P_5$ | $P_6$ | $P_7$ | $P_8$ | $P_9$ | $P_{10}$ | $P_{11}$ | $P_{12}$ | $P_{13}$ | $P_{14}$ |
|---|---|---|---|---|---|---|---|---|---|---|---|---|---|
| c | c | c | b | b | b | b | b | c | c | a | a | a | a |
| a | a | a | a | a | a | c | c | b | b | b | b | b | c |
| b | b | b | c | c | c | a | a | a | a | c | c | c | b |

Let us calculate the value $n^+(x, \vec{P}_{X'})$ for each alternative $x \in X'$: $n^+(a, \vec{P}_{X'}) = 4$, $n^+(b, \vec{P}_{X'}) = 5$, $n^+(c, \vec{P}_{X'}) = 5$. According to the plurality rule the alternatives $b$ and $c$ will be chosen, i.e., $C(\vec{P}_{X'}, X') = \{b, c\}$.

Finally, let us consider the subset $X'' = X\setminus\{a\}$. A contraction of the profile $\vec{P}_X$ onto a set $X''$, i.e., $\vec{P}_{X''}$, looks as

| $P_1$ | $P_2$ | $P_3$ | $P_4$ | $P_5$ | $P_6$ | $P_7$ | $P_8$ | $P_9$ | $P_{10}$ | $P_{11}$ | $P_{12}$ | $P_{13}$ | $P_{14}$ |
|---|---|---|---|---|---|---|---|---|---|---|---|---|---|
| c | c | c | b | b | b | d | d | d | d | b | b | b | c |
| d | d | d | d | d | d | b | b | c | c | c | c | c | b |
| b | b | b | c | c | c | c | c | b | b | d | d | d | d |

Let us calculate the value $n^+(x, \vec{P}_{X''})$ for each alternative $x \in X''$: $n^+(b, \vec{P}_{X''}) = 6$, $n^+(c, \vec{P}_{X''}) = 5$, $n^+(d, \vec{P}_{X''}) = 4$. According to the plurality rule the alternative $b$ will be chosen, i.e., $C(\vec{P}_{X''}, X'') = \{b\}$.

Then $C(\vec{P}_{X'}, X') \cap C(\vec{P}_{X''}, X'') = \{b\} \nsubseteq C(\vec{P}_X, X)$. Thus, the plurality rule does not satisfy the condition **C**.

### 2.3 Outcast condition (O)

Let us consider the example from paragraph 2.1 of this section. Then $C(\vec{P}_{X'}, X') \neq C(\vec{P}_X, X)$. Thus, the plurality rule does not satisfy the condition **O**.

### 2.4 Arrow's choice axiom (ACA)

The plurality rule does not satisfy the condition **ACA** as the condition **H** is not satisfied.

### 2.5 Monotonicity condition 1

Let $C(\vec{P}_X, X) = \{a\}$, i.e., $\forall x \in X\; n^+(a, \vec{P}_X) \geq n^+(x, \vec{P}_X)$.

Consider now a profile $\vec{P'_X}$, which differs from the profile $\vec{P}_X$ only by improved position of the alternative $a$. Then $n^+(a, \vec{P'_X}) \geq n^+(a, \vec{P}_X)$ and $\forall x \in X\setminus\{a\}\; n^+(x, \vec{P'_X}) \leq n^+(x, \vec{P}_X)$. Thus, the alternative $a$ will be chosen, i.e., $a \in C(\vec{P'_X}, X)$, as $n^+(a, \vec{P'_X}) \geq n^+(x, \vec{P'_X})$.

Since $a \in C(\vec{P}_X, X)$ and $a \in C(\vec{P'_X}, X)$, the plurality rule satisfies the Monotonicity condition 1.

### 2.6 The Monotonicity condition 2

Let $X = \{a, b, c\}$ and the profile $\vec{P}_X$ is the following

| $P_1$ | $P_2$ | $P_3$ | $P_4$ | $P_5$ |
|---|---|---|---|---|
| a | a | c | b | b |
| c | c | a | c | c |
| b | b | b | a | a |

Let us calculate the value $n^+(x, \vec{P}_X)$ for each alternative $x \in X$: $n^+(a, \vec{P}_X) = 2$, $n^+(b, \vec{P}_X) = 2$, $n^+(c, \vec{P}_X) = 1$. According to the plurality rule the alternatives $a$ and $b$ will be chosen, i.e., $C(\vec{P}_X, X) = \{a, b\}$.

Consider now the subset $X' = X \setminus \{b\}$. A contraction of the profile $\vec{P}_X$ onto a set $X'$, i.e., $\vec{P}_{X'}$, looks as

| $P_1$ | $P_2$ | $P_3$ | $P_4$ | $P_5$ |
|---|---|---|---|---|
| a | a | c | c | c |
| c | c | a | a | a |

According to the plurality rule the alternative $c$ will be chosen, i.e., $C(\vec{P}_{X'}, X') = \{c\}$.

Finally, let us consider the subset $X'' = X \setminus \{a\}$. A contraction of the profile $\vec{P}_X$ onto a set $X''$, i.e., $\vec{P}_{X''}$, looks as

| $P_1$ | $P_2$ | $P_3$ | $P_4$ | $P_5$ |
|---|---|---|---|---|
| c | c | c | b | b |
| b | b | b | c | c |

According to the plurality rule the alternative $c$ will be chosen, i.e., $C(\vec{P}_{X''}, X'') = \{c\}$.

Then $\{a, b\} \in C(\vec{P}_X, X)$, $\{a\} \notin C(\vec{P}_{X'}, X')$ and $\{b\} \notin C(\vec{P}_{X''}, X'')$. Thus, the plurality rule does not satisfy the Monotonicity condition 2.

*2.7 The strict monotonicity condition*

Let $X = \{a, b, c\}$ and the profile $\vec{P}_X$ is the following

| $P_1$ | $P_2$ | $P_3$ |
|---|---|---|
| a | a | b |
| c | c | c |
| b | b | a |

Let us calculate the value $n^+(x, \vec{P}_X)$ for each alternative $x \in X$: $n^+(a, \vec{P}_X) = 2$, $n^+(b, \vec{P}_X) = 1$, $n^+(c, \vec{P}_X) = 0$. According to the plurality rule the alternative $a$ will be chosen, i.e., $C(\vec{P}_X, X) = \{a\}$.

Consider now a profile $\vec{P'_X}$, which differs from the profile $\vec{P}_X$ only by improved position of the alternative $c$ in $P'_2$:

| $P'_1$ | $P'_2$ | $P'_3$ |
|---|---|---|
| a | c | b |
| c | a | c |
| b | b | a |

According to the rule the alternatives $a, b, c$ will be chosen, i.e., $C(\vec{P'_X}, X) = \{a, b, c\}$.

$$C(\vec{P'_X}, X) \neq \begin{bmatrix} C(\vec{P}_X, X) \text{ or} \\ \{c\} \text{ or} \\ C(\vec{P}_X, X) \cup \{c\}. \end{bmatrix}$$

Thus, the plurality rule does not satisfy the strict monotonicity condition.

*2.8 The non-compensatory condition*

The plurality rule does not satisfy the non-compensatory condition (see paragraph 1.8 of this section).

## 3. The inverse plurality rule
### 3.1 Heredity condition (H)
Let $X = \{a, b, c, d\}$ and the profile $\vec{P}_X$ is the following

| $P_1$ | $P_2$ | $P_3$ |
|---|---|---|
| b | d | a |
| a | b | d |
| c | a | c |
| d | c | b |

According to the inverse plurality rule the alternative $a$ will be chosen, i.e., $C(\vec{P}_X, X) = \{a\}$.
Consider now the subset $X' = X \setminus \{c, d\}$. A contraction of the profile $\vec{P}_X$ onto a set $X'$, i.e., $\vec{P}_{X'}$, looks as

| $P_1$ | $P_2$ | $P_3$ |
|---|---|---|
| b | b | a |
| a | a | b |

According to the inverse plurality rule the alternative $b$ will be chosen, i.e., $C(\vec{P}_{X'}, X') = \{b\}$.
Then $C(\vec{P}_{X'}, X') \not\supseteq C(\vec{P}_X, X) \cap X'$. Thus, the inverse plurality rule does not satisfy the condition **H**.

### 3.2 Concordance condition (C)
Consider the example from paragraph 3.1 of this section. According to the inverse plurality rule the alternative $a$ will be chosen, i.e., $C(\vec{P}_X, X) = \{a\}$.
Consider now the subset $X' = X \setminus \{c\}$. A contraction of the profile $\vec{P}_X$ onto a set $X'$, i.e., $\vec{P}_{X'}$, looks as

| $P_1$ | $P_2$ | $P_3$ |
|---|---|---|
| b | d | a |
| a | b | d |
| d | a | b |

According to the inverse plurality rule the alternatives $a,b,d$ will be chosen, i.e., $C(\vec{P}_{X'}, X') = \{a, b, d\}$.
Finally, let us consider the subset $X'' = X \setminus \{a, d\}$. A contraction of the profile $\vec{P}_X$ onto a set $X''$, i.e., $\vec{P}_{X''}$, looks as

| $P_1$ | $P_2$ | $P_3$ |
|---|---|---|
| b | b | c |
| c | c | b |

According to the inverse plurality rule the alternative $b$ will be chosen, i.e., $C(\vec{P}_{X''}, X'') = \{b\}$.
Then $C(\vec{P}_{X'}, X') \cap C(\vec{P}_{X''}, X'') = \{b\} \not\subseteq C(\vec{P}_X, X)$. Thus, the inverse plurality rule does not satisfy the condition **C**.

### 3.3 Outcast condition (O)
Consider the example from paragraph 3.1 of this section. $C(\vec{P}_{X'}, X') \neq C(\vec{P}_X, X)$, consequently, the inverse plurality rule does not satisfy the condition **O**.

### 3.4 Arrow's choice axiom (ACA)
The inverse plurality rule does not satisfy the condition **ACA** as the condition **H** is not satisfied.

### 3.5 Monotonicity condition 1

Let $C(\vec{P}_X, X) = \{a\}$, i.e., $\forall x \in X\ n^-(a, \vec{P}_X) \leq n^-(x, \vec{P}_X)$.

Consider now a profile $\vec{P}_X'$, which differs from the profile $\vec{P}_X$ only by improved position of the alternative $a$. Then $n^-(a, \vec{P}_X') \leq n^-(a, \vec{P}_X)$ and $\forall x \in X \setminus \{a\}\ n^-(x, \vec{P}_X') \geq n^-(x, \vec{P}_X)$. Thus, the alternative $a$ will be chosen, i.e., $a \in C(\vec{P}_X', X)$, as $n^-(a, \vec{P}_X') \leq n^-(x, \vec{P}_X')$.

Then $a \in C(\vec{P}_X, X)$ and $a \in C(\vec{P}_X', X)$. Thus, the inverse plurality rule satisfies the Monotonicity condition 1.

### 3.6 The Monotonicity condition 2

Let $X = \{a, b, c\}$ and the profile $\vec{P}_X$ is the following

| $P_1$ | $P_2$ | $P_3$ | $P_4$ | $P_5$ | $P_6$ | $P_7$ |
|---|---|---|---|---|---|---|
| c | c | a | a | a | c | c |
| a | a | b | b | b | b | b |
| b | b | c | c | c | a | a |

According to the inverse plurality rule the alternatives $a,b$ will be chosen, i.e., $C(\vec{P}_X, X) = \{a, b\}$. Consider now the subset $X' = X \setminus \{b\}$. A contraction of the profile $\vec{P}_X$ onto a set $X'$, i.e., $\vec{P}_{X'}$, looks as

| $P_1$ | $P_2$ | $P_3$ | $P_4$ | $P_5$ | $P_6$ | $P_7$ |
|---|---|---|---|---|---|---|
| c | c | a | a | a | c | c |
| a | a | c | c | c | a | a |

According to the inverse plurality rule the alternative $c$ will be chosen, i.e., $C(\vec{P}_{X'}, X') = \{c\}$.

Finally, let us consider the subset $X'' = X \setminus \{a\}$. A contraction of the profile $\vec{P}_X$ onto a set $X''$, i.e., $\vec{P}_{X''}$, looks as

| $P_1$ | $P_2$ | $P_3$ | $P_4$ | $P_5$ | $P_6$ | $P_7$ |
|---|---|---|---|---|---|---|
| c | c | b | b | b | c | c |
| b | b | c | c | c | b | b |

According to the inverse plurality rule the alternative $c$ will be chosen, i.e., $C(\vec{P}_{X''}, X'') = \{c\}$.

Then $\{a, b\} \in C(\vec{P}_X, X)$, $\{a\} \notin C(\vec{P}_{X'}, X')$ and $\{b\} \notin C(\vec{P}_{X''}, X'')$. Thus, the inverse plurality rule does not satisfy the Monotonicity condition 2.

### 3.7 The strict monotonicity condition

Let $X = \{a, b, c\}$ and the profile $\vec{P}_X$ is the following

| $P_1$ | $P_2$ | $P_3$ |
|---|---|---|
| c | a | b |
| a | c | c |
| b | b | a |

According to the inverse plurality rule the alternative $c$ will be chosen, i.e., $C(\vec{P}_X, X) = \{c\}$.

Consider now a profile $\vec{P'_X}$, which differs from the profile $\vec{P_X}$ only by improved position of the alternative $b$ in $P'_2$:

|  $P'_1$ | $P'_2$ | $P'_3$ |
| --- | --- | --- |
| c | a | b |
| a | b | c |
| b | c | a |

According to the inverse plurality rule the alternatives $a,b,c$ will be chosen, i.e., $C(\vec{P'_X}, X) = \{a, b, c\}$.

$$C(\vec{P'_X}, X) \neq \begin{bmatrix} C(\vec{P_X}, X) \ or \\ \{b\} \ or \\ C(\vec{P_X}, X) \cup \{b\}. \end{bmatrix}$$

Thus, the inverse plurality rule does not satisfy the strict monotonicity condition.

### 3.8 The non-compensatory condition

Let $X = \{a, b, c\}$ and the profile $\vec{P_X}$ is the following

| $P_1$ | $P_2$ | $P_3$ |
| --- | --- | --- |
| a | a | b |
| b | b | a |
| c | c | c |

According to the inverse plurality rule the alternatives $a,b$ will be chosen, i.e., $C(\vec{P_X}, X) = \{a, b\}$.
Let us write the profile $\vec{P_X}$ in the following form

|  | $\varphi_1$ | $\varphi_2$ | $\varphi_3$ |
| --- | --- | --- | --- |
| a | 3 | 3 | 2 |
| b | 2 | 2 | 3 |
| c | 1 | 1 | 1 |

According to the non-compensatory condition the alternative $a$ is better than the alternative $b$ while the alternative $b$ is better than the alternative $c$. Thus, the non-compensatory condition is not satisfied as $b \in C(\vec{P_X}, X)$.

## 4. The q-Approval rule

The q-Approval rule satisfies the same conditions as the plurality rule regardless the value of parameter $q$. For the case $q = 1$ the proof is provided in paragraph 2 of this section. For the case $q > 1$ we can use the same examples but with larger number of alternatives. Let us prove it for the condition **H**.

### 4.1 Heredity condition (H)

Suppose the set $X$ is finite and the profile $\vec{P_X}$ is the following

| $P_1$ | $P_2$ | $P_3$ | $P_4$ | $P_5$ |
| --- | --- | --- | --- | --- |
| a | a | c | b | b |
| b | b | b | a | a |
| c | c | a | c | c |

$\forall y \in X\setminus\{a, b, c\}: n^+(y, \vec{P_X}, q) = 1$

$\forall y \in X\setminus\{a, b, c\}$

Let us calculate the value $n^+(x, \vec{P}_X, q)$ for each alternative $x \in X$: $n^+(a, \vec{P}_X, q) = 2$, $n^+(b, \vec{P}_X, q) = 2$, $n^+(c, \vec{P}_X, q) = 1$, $n^+(y, \vec{P}_X, q) = 1$. According to the q-Approval rule the alternatives $a$ and $b$ will be chosen, i.e., $C(\vec{P}_X, X) = \{a, b\}$.

Consider now the subset $X' = X \setminus \{c\}$. A contraction of the profile $\vec{P}_X$ onto a set $X'$, i.e., $\vec{P}_{X'}$, looks as

|   | $P_1$ | $P_2$ | $P_3$ | $P_4$ | $P_5$ |
|---|---|---|---|---|---|
|   | $a$ | $a$ | ... | $b$ | $b$ |
| $\forall y \in X' \setminus \{a,b\}: n^+(y, \vec{P}_{X'}, q) = 1$ | ... | ... | ...$b$... | ... | ... |
|   | $b$ | $b$ | $a$ | $a$ | $a$ |
| $\forall y \in X' \setminus \{a,b\}$ | ... | ... | ... | ... | ... |

$\} q$

Let us calculate the value $n^+(x, \vec{P}_{X'}, q)$ for each alternative $x \in X'$: $n^+(a, \vec{P}_{X'}, q) = 2$, $n^+(b, \vec{P}_{X'}, q) = 3$, $n^+(y, \vec{P}_{X'}, q) = 1$. According to the rule the alternative $b$ will be chosen, i.e., $C(\vec{P}_{X'}, X', q) = \{b\}$.

Then $C(\vec{P}_{X'}, X', q) \not\supseteq C(\vec{P}_X, X, q) \cap X'$. Thus, the q-Approval rule does not satisfy the condition **H**.

## 5. The run-off procedure

### 5.1 Heredity condition (H)

Let $X = \{a, b, c, d, e, f\}$ and the profile $\vec{P}_X$ is the following

| $P_1$ | $P_2$ | $P_3$ | $P_4$ | $P_5$ | $P_6$ | $P_7$ | $P_8$ |
|---|---|---|---|---|---|---|---|
| $a$ | $a$ | $b$ | $b$ | $c$ | $d$ | $e$ | $f$ |
| $b$ | $b$ | $a$ | $a$ | $b$ | $c$ | $a$ | $c$ |
| $c$ | $c$ | $c$ | $c$ | $a$ | $b$ | $b$ | $b$ |
| $d$ | $d$ | $d$ | $d$ | $d$ | $a$ | $c$ | $a$ |
| $e$ | $e$ | $e$ | $e$ | $e$ | $e$ | $d$ | $d$ |
| $f$ | $f$ | $f$ | $f$ | $f$ | $f$ | $f$ | $e$ |

Let us calculate the value $n^+(x, \vec{P}_X)$ for each alternative $x \in X$: $n^+(a, \vec{P}_X) = 2$, $n^+(b, \vec{P}_X) = 2$, $n^+(c, \vec{P}_X) = 1$, $n^+(d, \vec{P}_X) = 1$, $n^+(e, \vec{P}_X) = 1$, $n^+(f, \vec{P}_X) = 1$.

According to the simple majority rule the choice is empty. Consider now the subset $X' = X \setminus \{c, d, e, f\}$. A contraction of the profile $\vec{P}_X$ onto a set $X'$, i.e., $\vec{P}_{X'}$, looks as

| $P_1$ | $P_2$ | $P_3$ | $P_4$ | $P_5$ | $P_6$ | $P_7$ | $P_8$ |
|---|---|---|---|---|---|---|---|
| $a$ | $a$ | $b$ | $b$ | $b$ | $b$ | $a$ | $b$ |
| $b$ | $b$ | $a$ | $a$ | $a$ | $a$ | $b$ | $a$ |

Thus, according to the run-off procedure the alternative $b$ will be chosen, i.e., $C(\vec{P}_X, X) = \{b\}$.

Consider now the subset $X'' = X \setminus \{d, e, f\}$. A contraction of the profile $\vec{P}_X$ onto a set $X''$, i.e., $\vec{P}_{X''}$, looks as

| $P_1$ | $P_2$ | $P_3$ | $P_4$ | $P_5$ | $P_6$ | $P_7$ | $P_8$ |
|---|---|---|---|---|---|---|---|
| $a$ | $a$ | $b$ | $b$ | $c$ | $c$ | $a$ | $c$ |
| $b$ | $b$ | $a$ | $a$ | $b$ | $b$ | $b$ | $b$ |
| $c$ | $c$ | $c$ | $c$ | $a$ | $a$ | $c$ | $a$ |

According to the simple majority rule the choice is empty. Consider the subset $X''' = X''\backslash\{b\}$. A contraction of the profile $\vec{P}_{X''}$ onto a set $X'''$, i.e., $\vec{P}_{X'''}$, looks as

| $P_1$ | $P_2$ | $P_3$ | $P_4$ | $P_5$ | $P_6$ | $P_7$ | $P_8$ |
|---|---|---|---|---|---|---|---|
| a | a | a | a | c | c | a | c |
| c | c | c | c | a | a | c | a |

According to the rule the alternative $a$ will be chosen, consequently, $C(\vec{P}_{X''}, X'') = \{a\}$.
Then $C(\vec{P}_{X''}, X'') \not\supseteq C(\vec{P}_X, X) \cap X''$. Thus, the run-off procedure does not satisfy the condition **H**.

*5.2 Concordance condition (C)*
The run-off procedure does not satisfy the condition **C** (see paragraph 1.2 of this section).

*5.3 Outcast condition (O)*
The run-off procedure does not satisfy the condition **O** (see paragraph 5.1 of this section).

*5.4 Arrow's choice axiom (ACA)*
The run-off procedure does not satisfy the condition **ACA** as the condition **H** is not satisfied.

*5.5 Monotonicity condition 1*
Let $X = \{a, b, c, d\}$ and the profile $\vec{P}_X$ is the following

| $P_1$ | $P_2$ | $P_3$ | $P_4$ | $P_5$ | $P_6$ | $P_7$ |
|---|---|---|---|---|---|---|
| a | b | c | a | b | c | d |
| c | a | b | c | a | b | c |
| b | c | a | b | c | a | b |
| d | d | d | d | d | d | a |

According to the simple majority rule the choice is empty. Thus, the alternatives $c$ and $d$ are omitted. Consider now the subset $X' = X\backslash\{c, d\}$. A contraction of the profile $\vec{P}_X$ onto a set $X'$, i.e., $\vec{P}_{X'}$, looks as

| $P_1$ | $P_2$ | $P_3$ | $P_4$ | $P_5$ | $P_6$ | $P_7$ |
|---|---|---|---|---|---|---|
| a | b | b | a | b | b | b |
| b | a | a | b | a | a | a |

According to the rule the alternative $b$ will be chosen, consequently, $C(\vec{P}_X, X) = \{b\}$.
Suppose now that the position of the alternative $b$ was improved in $P_1'$ while the relative comparison of any pair of other alternatives remained unchanged.

| $P_1'$ | $P_2'$ | $P_3'$ | $P_4'$ | $P_5'$ | $P_6'$ | $P_7'$ |
|---|---|---|---|---|---|---|
| b | b | c | a | b | c | d |
| a | a | b | c | a | b | c |
| c | c | a | b | c | a | b |
| d | d | d | d | d | d | a |

According to the simple majority rule the choice is empty. Thus, the alternatives $a$ and $d$ are omitted. Next, consider the subset $X'' = X\backslash\{a, d\}$. A contraction of the profile $\vec{P}_X'$ onto a set $X''$, i.e., $\overrightarrow{P'_{X''}}$, looks as

| $P_1'$ | $P_2'$ | $P_3'$ | $P_4'$ | $P_5'$ | $P_6'$ | $P_7'$ |
|---|---|---|---|---|---|---|
| b | b | c | c | b | c | c |
| c | c | b | b | c | b | b |

According to the rule the alternative $c$ will be chosen, i.e., $C\left(\vec{P'_X}, X\right) = C\left(\vec{P''_{X''}}, X''\right) = \{c\}$.

Then $\{b\} \in C(\vec{P}_X, X)$, $\{b\} \notin C\left(\vec{P'_X}, X\right)$. Thus, the run-off procedure does not satisfy the Monotonicity condition 1.

*5.6 The Monotonicity condition 2*
Since the run-off procedure always chooses no more than one best alternative, the Monotonicity condition 2 is not applicable to it as it considers the choice of more than two alternatives. In other words, the run-off procedure obeys the Monotonicity condition 2 trivially.

*5.7 The strict monotonicity condition*
Let $X = \{a, b, c\}$ and the profile $\vec{P}_X$ is the following

| $P_1$ | $P_2$ | $P_3$ |
|---|---|---|
| c | c | a |
| a | b | b |
| b | a | c |

According to the run-off procedure the alternative $c$ will be chosen, i.e., $C(\vec{P}_X, X) = \{c\}$.

Consider now a profile $\vec{P'_X}$, which differs from the profile $\vec{P}_X$ only by improved position of the alternative $b$ in $P'_2$:

| $P'_1$ | $P'_2$ | $P'_3$ |
|---|---|---|
| c | b | a |
| a | c | b |
| b | a | c |

According to the rule the alternative $a$ will be chosen, i.e., $C(\vec{P'_X}, X) = \{a\}$.

$$C\left(\vec{P'_X}, X\right) \neq \begin{bmatrix} C(\vec{P}_X, X) \text{ or} \\ \{b\} \text{ or} \\ C(\vec{P}_X, X) \cup \{b\}. \end{bmatrix}$$

Thus, the run-off procedure does not satisfy the strict monotonicity condition.

*5.8 The non-compensatory condition*
The run-off procedure does not satisfy the non-compensatory condition (see paragraph 1.8 of this section).

# 6. The Hare rule (the Ware procedure)
*6.1 Heredity condition (H)*
Let $X = \{a, b, c, d\}$ and the profile $\vec{P}_X$ is the following

| $P_1$ | $P_2$ | $P_3$ | $P_4$ | $P_5$ | $P_6$ | $P_7$ |
|---|---|---|---|---|---|---|
| a | a | a | c | c | d | d |
| b | b | b | b | b | b | b |
| c | c | c | d | d | a | a |
| d | d | d | a | a | c | c |

According to the Hare rule the alternative $a$ will be chosen, i.e., $C(\vec{P}_X, X) = \{a\}$.

Consider now the subset $X' = X\setminus\{c, d\}$. A contraction of the profile $\vec{P}_X$ onto a set $X'$, i.e., $\vec{P}_{X'}$, looks as

| $P_1$ | $P_2$ | $P_3$ | $P_4$ | $P_5$ | $P_6$ | $P_7$ |
|---|---|---|---|---|---|---|
| a | a | a | b | b | b | b |
| b | b | b | a | a | a | a |

According to the rule the alternative $b$ will be chosen, i.e., $C(\vec{P}_{X'}, X') = \{b\}$.

Then $C(\vec{P}_{X'}, X') \not\supseteq C(\vec{P}_X, X) \cap X'$. Thus, the Hare rule (the Ware procedure) does not satisfy the condition **H**.

### 6.2 Concordance condition (C)

Consider the previous example. According to the Hare rule (the Ware procedure) the alternative $a$ will be chosen, i.e., $C(\vec{P}_X, X) = \{a\}$.

Consider now the subset $X' = X \setminus \{c, d\}$. A contraction of the profile $\vec{P}_X$ onto a set $X'$, i.e., $\vec{P}_{X'}$, looks as

| $P_1$ | $P_2$ | $P_3$ | $P_4$ | $P_5$ | $P_6$ | $P_7$ |
|---|---|---|---|---|---|---|
| a | a | a | b | b | b | b |
| b | b | b | a | a | a | a |

According to the Hare rule the alternative $b$ will be chosen, i.e., $C(\vec{P}_{X'}, X') = \{b\}$.

Finally, let us consider the subset $X'' = X \setminus \{a\}$. A contraction of the profile $\vec{P}_X$ onto a set $X''$, i.e., $\vec{P}_{X''}$, looks as

| $P_1$ | $P_2$ | $P_3$ | $P_4$ | $P_5$ | $P_6$ | $P_7$ |
|---|---|---|---|---|---|---|
| b | b | b | c | c | d | d |
| c | c | c | b | b | b | b |
| d | d | d | d | d | c | c |

According to the rule the alternative $b$ will be chosen, i.e., $C(\vec{P}_{X''}, X'') = \{b\}$.

Then $C(\vec{P}_{X'}, X') \cap C(\vec{P}_{X''}, X'') = \{b\} \nsubseteq C(\vec{P}_X, X)$. Thus, the Hare rule (the Ware procedure) does not satisfy the condition **C**.

### 6.3 Outcast condition (O)

Consider the example from paragraph 6.1 of this section. Then $C(\vec{P}_{X'}, X') \neq C(\vec{P}_X, X)$. Thus, the Hare rule (the Ware procedure) does not satisfy the condition **O**.

### 6.4 Arrow's choice axiom (ACA)

The Hare rule (the Ware procedure) does not satisfy the condition **ACA** as the condition **C** is not satisfied.

### 6.5 Monotonicity condition 1

Let $X = \{a, b, c, d\}$ and the profile $\vec{P}_X$ is the following

| $P_1$ | $P_2$ | $P_3$ | $P_4$ | $P_5$ | $P_6$ | $P_7$ | $P_8$ | $P_9$ | $P_{10}$ | $P_{11}$ | $P_{12}$ | $P_{13}$ |
|---|---|---|---|---|---|---|---|---|---|---|---|---|
| a | a | a | a | a | c | c | c | c | d | d | d | d |
| b | b | b | b | b | b | b | b | b | b | b | b | b |
| c | c | c | c | c | d | d | d | d | a | a | a | a |
| d | d | d | d | d | a | a | a | a | c | c | c | c |

According to the rule the alternative $a$ will be chosen, i.e., $C(\vec{P}_X, X) = \{a\}$.

Suppose now that the position of the alternative $a$ was improved in $P'_6$ while the relative comparison of any pair of other alternatives remained unchanged.

| $P'_1$ | $P'_2$ | $P'_3$ | $P'_4$ | $P'_5$ | $P'_6$ | $P'_7$ | $P'_8$ | $P'_9$ | $P'_{10}$ | $P'_{11}$ | $P'_{12}$ | $P'_{13}$ |
|---|---|---|---|---|---|---|---|---|---|---|---|---|
| a | a | a | a | a | a | c | c | c | d | d | d | d |
| b | b | b | b | b | c | b | b | b | b | b | b | b |
| c | c | c | c | c | b | d | d | d | a | a | a | a |
| d | d | d | d | d | d | a | a | a | c | c | c | c |

According to the Hare rule (the Ware procedure) the alternative $c$ is omitted on the first stage. Next, consider the subset $X' = X\setminus\{c\}$. A contraction of the profile $\vec{P'_X}$ onto a set $X'$, i.e., $\vec{P'_{X'}}$, looks as

| $P'_1$ | $P'_2$ | $P'_3$ | $P'_4$ | $P'_5$ | $P'_6$ | $P'_7$ | $P'_8$ | $P'_9$ | $P'_{10}$ | $P'_{11}$ | $P'_{12}$ | $P'_{13}$ |
|---|---|---|---|---|---|---|---|---|---|---|---|---|
| a | a | a | a | a | a | b | b | b | b | b | b | b |
| b | b | b | b | b | b | a | a | a | a | a | a | a |

According to the rule the alternative $b$ will be chosen, i.e., $C(\vec{P'_X}, X) = C(\vec{P'_{X'}}, X') = \{b\}$.

Then $\{a\} \in C(\vec{P_X}, X)$, $\{a\} \notin C(\vec{P'_X}, X)$. Thus, the Hare rule (the Ware procedure) does not satisfy the Monotonicity condition 1.

### 6.6 The Monotonicity condition 2
Since the Hare rule (the Ware procedure) always chooses no more than one best alternative, the Monotonicity condition 2 is not applicable to it as it considers the choice of more than two alternatives. In other words, the Hare rule (the Ware procedure) obeys the Monotonicity condition 2 trivially.

### 6.7 The strict monotonicity condition
The Hare rule (the Ware procedure) does not satisfy the strict monotonicity condition (see the proof in paragraph 5.7 of this section).

### 6.8 The non-compensatory condition
The Hare rule (the Ware procedure) does not satisfy the non-compensatory condition (see paragraph 1.8 of this section).

## 7. The Borda rule

### 7.1 Heredity condition (H)
Let $X = \{a, b, c, d, e\}$ and the profile $\vec{P_X}$ is the following

| $P_1$ | $P_2$ | $P_3$ | $P_4$ | $P_5$ |
|---|---|---|---|---|
| e | e | a | b | b |
| a | a | d | c | c |
| b | b | c | e | a |
| c | c | b | a | d |
| d | d | e | d | e |

Let us calculate the Borda count for each alternative: $r(a, \vec{P_X}) = r(b, \vec{P_X}) = 13, r(c, \vec{P_X}) = 10, r(d, \vec{P_X}) = 4, r(e, \vec{P_X}) = 10$. According to the Borda rule the alternatives $a,b$ will be chosen, i.e., $C(\vec{P_X}, X) = \{a, b\}$.

Consider now the subset $X' = X\setminus\{d\}$. A contraction of the profile $\vec{P_X}$ onto a set $X'$, i.e., $\vec{P_{X'}}$, looks as

|     | P₁ | P₂ | P₃ | P₄ | P₅ |
| --- | --- | --- | --- | --- | --- |
|     | e | e | a | b | b |
|     | a | a | c | c | c |
|     | b | b | b | e | a |
|     | c | c | e | a | e |

Let us calculate the Borda count for each alternative: $r(a, \vec{P}_{X'}) = 8, r(b, \vec{P}_{X'}) = 9, r(c, \vec{P}_{X'}) = 6, r(e, \vec{P}_{X'}) = 7$. According to the rule the alternative $b$ will be chosen, i.e., $C(\vec{P}_{X'}, X') = \{b\}$.

Then $C(\vec{P}_{X'}, X') \not\supseteq C(\vec{P}_X, X) \cap X'$. Thus, the Borda rule does not satisfy the condition **H**.

### 7.2 Concordance condition (C)

Let $X = \{a, b, c, d\}$ and the profile $\vec{P}_X$ is the following

|     | P₁ | P₂ | P₃ | P₄ |
| --- | --- | --- | --- | --- |
|     | a | c | b | c |
|     | d | d | a | b |
|     | b | a | d | a |
|     | c | b | c | d |

Let us calculate the Borda count for each alternative: $r(a, \vec{P}_X) = 7, r(b, \vec{P}_X) = r(c, \vec{P}_X) = 6, r(d, \vec{P}_X) = 5$. According to the Borda rule the alternative $a$ will be chosen, i.e., $C(\vec{P}_X, X) = \{a\}$.

Consider now the subset $X' = X \setminus \{d\}$. A contraction of the profile $\vec{P}_X$ onto a set $X'$, i.e., $\vec{P}_{X'}$, looks as

|     | P₁ | P₂ | P₃ | P₄ |
| --- | --- | --- | --- | --- |
|     | a | c | b | c |
|     | b | a | a | b |
|     | c | b | c | a |

Let us calculate the Borda count for each alternative: $r(a, \vec{P}_{X'}) = r(b, \vec{P}_{X'}) = r(c, \vec{P}_{X'}) = 4$. According to the Borda rule the alternatives $a,b,c$ will be chosen, i.e., $C(\vec{P}_{X'}, X') = \{a, b, c\}$.

Finally, let us consider the subset $X'' = X \setminus \{a\}$. A contraction of the profile $\vec{P}_X$ onto a set $X''$, i.e., $\vec{P}_{X''}$, looks as

|     | P₁ | P₂ | P₃ | P₄ |
| --- | --- | --- | --- | --- |
|     | d | c | b | c |
|     | b | d | d | b |
|     | c | b | c | d |

Let us calculate the Borda count for each alternative: $r(b, \vec{P}_{X''}) = r(c, \vec{P}_{X''}) = r(d, \vec{P}_{X''}) = 4$. According to the rule the alternatives $b,c,d$ will be chosen, i.e., $C(\vec{P}_{X''}, X'') = \{b, c, d\}$.

Then $C(\vec{P}_{X'}, X') \cap C(\vec{P}_{X''}, X'') = \{b, c\} \not\subseteq C(\vec{P}_X, X)$. Thus, the Borda rule does not satisfy the condition **C**.

### 7.3 Outcast condition (O)

Consider the example from paragraph 7.2 of this section. Then $C(\vec{P}_{X'}, X') \neq C(\vec{P}_X, X)$, consequently, the Borda rule does not satisfy the condition **O**.

## 7.4 Arrow's choice axiom (ACA)

The Borda rule does not satisfy the condition **ACA** as the condition **H** is not satisfied.

## 7.5 Monotonicity condition 1

Let $a \in C(\vec{P}_X, X)$, i.e., $\forall x \in X \setminus \{a\}\ r(a, \vec{P}_X) \geq r(x, \vec{P}_X)$.

Consider now a profile $\vec{P'_X}$, which differs from the profile $\vec{P}_X$ only by improved position of the alternative $a$. Then $r\left(a, \vec{P'_X}\right) > r(a, \vec{P}_X)$ and $r\left(x, \vec{P'_X}\right) \leq r(x, \vec{P}_X)$. Thus, $r\left(a, \vec{P'_X}\right) > r(x, \vec{P'_X})$ and $a \in C(\vec{P'_X}, X)$.

$a \in C(\vec{P}_X, X)$ and $a \in C(\vec{P'_X}, X)$. Thus, the Borda rule satisfies the Monotonicity condition 1.

## 7.6 The Monotonicity condition 2

Let $X = \{a, b, c, d, e, f\}$ and the profile $\vec{P}_X$ is the following

| $P_1$ | $P_2$ | $P_3$ | $P_4$ | $P_5$ | $P_6$ | $P_7$ | $P_8$ |
|---|---|---|---|---|---|---|---|
| a | a | a | f | f | f | f | c |
| b | b | b | b | b | b | c | f |
| c | c | c | a | a | a | e | e |
| d | d | d | c | c | c | d | d |
| e | e | e | d | d | d | b | a |
| f | f | f | e | e | e | a | b |

Let us calculate the Borda count for each alternative: $r(a, \vec{P}_X) = r(b, \vec{P}_X) = 25$, $r(c, \vec{P}_X) = 24$, $r(d, \vec{P}_X) = 13$, $r(e, \vec{P}_X) = 9$, $r(f, \vec{P}_X) = 24$. According to the Borda rule the alternatives $a,b$ will be chosen, i.e., $C(\vec{P}_X, X) = \{a, b\}$.

Consider now the subset $X' = X \setminus \{b\}$. A contraction of the profile $\vec{P}_X$ onto a set $X'$, i.e., $\vec{P}_{X'}$, looks as

| $P_1$ | $P_2$ | $P_3$ | $P_4$ | $P_5$ | $P_6$ | $P_7$ | $P_8$ |
|---|---|---|---|---|---|---|---|
| a | a | a | f | f | f | f | c |
| c | c | c | a | a | a | c | f |
| d | d | d | c | c | c | e | e |
| e | e | e | d | d | d | d | d |
| f | f | f | e | e | e | a | a |

Let us calculate the Borda count for each alternative: $r(a, \vec{P}_{X'}) = 21$, $r(c, \vec{P}_{X'}) = 22$, $r(d, \vec{P}_{X'}) = 13$, $r(e, \vec{P}_{X'}) = 7$, $r(f, \vec{P}_{X'}) = 19$. According to the Borda rule the alternative $c$ will be chosen, i.e., $C(\vec{P}_{X'}, X') = \{c\}$.

Finally, let us consider the subset $X'' = X \setminus \{a\}$. A contraction of the profile $\vec{P}_X$ onto a set $X''$, i.e., $\vec{P}_{X''}$, looks as

| $P_1$ | $P_2$ | $P_3$ | $P_4$ | $P_5$ | $P_6$ | $P_7$ | $P_8$ |
|---|---|---|---|---|---|---|---|
| b | b | b | f | f | f | f | c |
| c | c | c | b | b | b | c | f |
| d | d | d | c | c | c | e | e |
| e | e | e | d | d | d | d | d |
| f | f | f | e | e | e | b | b |

Let us calculate the Borda count for each alternative: $r(b, \vec{P}_{X''}) = 21$, $r(c, \vec{P}_{X''}) = 22$, $r(d, \vec{P}_{X''}) = 13$, $r(e, \vec{P}_{X''}) = 7$, $r(f, \vec{P}_{X''}) = 19$. According to the rule the alternative $c$ will be chosen, i.e., $C(\vec{P}_{X''}, X'') = \{c\}$.

Then $\{a, b\} \in C(\vec{P}_X, X)$, $\{a\} \notin C(\vec{P}_{X'}, X')$ and $\{b\} \notin C(\vec{P}_{X''}, X'')$. Thus, the Borda rule does not satisfy the Monotonicity condition 2.

### 7.7 The strict monotonicity condition

Let $X = \{a, b, c\}$ and the profile $\vec{P}_X$ is the following

| $P_1$ | $P_2$ | $P_3$ | $P_4$ |
|---|---|---|---|
| a | a | b | b |
| c | b | a | a |
| b | c | c | c |

Let us calculate the Borda count for each alternative: $r(a, \vec{P}_X) = 6, r(b, \vec{P}_X) = 5, r(c, \vec{P}_X) = 1$. According to the Borda rule the alternative $a$ will be chosen, i.e., $C(\vec{P}_X, X) = \{a\}$.

Consider now a profile $\vec{P'_X}$, which differs from the profile $\vec{P}_X$ only by improved position of the alternative $c$ in $P'_1$:

| $P'_1$ | $P'_2$ | $P'_3$ | $P'_4$ |
|---|---|---|---|
| c | a | b | b |
| a | b | a | a |
| b | c | c | c |

Let us calculate the Borda count for each alternative: $r(a, \vec{P'_X}) = 5, r(b, \vec{P'_X}) = 5, r(c, \vec{P'_X}) = 2$.

According to the rule the alternatives $a,b$ will be chosen, i.e., $C(\vec{P'_X}, X) = \{a, b\}$.

$$C(\vec{P'_X}, X) \neq \begin{bmatrix} C(\vec{P}_X, X) \text{ or} \\ \{c\} \text{ or} \\ C(\vec{P}_X, X) \cup \{c\}. \end{bmatrix}$$

Thus, the Borda rule does not satisfy the strict monotonicity condition.

### 7.8 The non-compensatory condition

Let $X = \{a, b, c\}$ and the profile $\vec{P}_X$ is the following

| $P_1$ | $P_2$ | $P_3$ |
|---|---|---|
| a | a | b |
| b | b | c |
| c | c | a |

According to the rule the alternatives $a,b$ will be chosen, i.e., $C(\vec{P}_X, X) = \{a, b\}$.

Let us write the profile $\vec{P}_X$ in the following form

|   | $\varphi_1$ | $\varphi_2$ | $\varphi_3$ |
|---|---|---|---|
| a | 3 | 3 | 1 |
| b | 2 | 2 | 3 |
| c | 1 | 1 | 2 |

According to the non-compensatory condition the alternative $b$ is better than the alternative $a$ and the alternative $a$ is better than the alternative $c$. Thus, the non-compensatory condition is not satisfied as $a \in C(\vec{P}_X, X)$.

## 8. The Black procedure
*8.1 Heredity condition (H)*

Let $X = \{a, b, c\}$ and the profile $\vec{P}_X$ is the following

| $P_1$ | $P_2$ | $P_3$ | $P_4$ | $P_5$ |
|---|---|---|---|---|
| a | a | d | b | b |
| b | b | a | d | d |
| d | d | b | a | a |

In this example none of alternatives is a Condorcet winner ($a\mu b, d\mu a, b\mu d$). Thus, the Borda rule is applied.

Let us calculate the Borda count for each alternative: $r(a, \vec{P}_X) = 5, r(b, \vec{P}_X) = 6, r(c, \vec{P}_X) = 4$.

According to the Black procedure the alternative $b$ will be chosen, i.e., $C(\vec{P}_X, X) = \{b\}$.

Consider now the subset $X' = X \backslash \{d\}$. A contraction of the profile $\vec{P}_X$ onto a set $X'$, i.e., $\vec{P}_{X'}$, looks as

| $P_1$ | $P_2$ | $P_3$ | $P_4$ | $P_5$ |
|---|---|---|---|---|
| a | a | a | b | b |
| b | b | b | a | a |

According to the rule the alternative $a$ will be chosen, i.e., $C(\vec{P}_{X'}, X') = \{a\}$.

Then $C(\vec{P}_{X'}, X') \not\supseteq C(\vec{P}_X, X) \cap X'$. Thus, the Black procedure does not satisfy the condition **H**.

*8.2 Concordance condition (C)*

Let $X = \{a, b, c, d\}$ and the profile $\vec{P}_X$ is the following

| $P_1$ | $P_2$ | $P_3$ | $P_4$ | $P_5$ |
|---|---|---|---|---|
| a | a | c | d | d |
| d | d | b | b | b |
| c | c | a | a | a |
| b | b | d | c | c |

In this example none of alternatives is a Condorcet winner ($a\mu c, a\mu d, b\mu a, c\mu b, d\mu b, d\mu c$). Thus, the Borda rule is used for this case. Let us calculate the Borda count for each alternative: $r(a, \vec{P}_X) = 9, r(b, \vec{P}_X) = 6, r(c, \vec{P}_X) = 5, r(d, \vec{P}_X) = 10$. According to the Black procedure the alternative $d$ will be chosen, i.e., $C(\vec{P}_X, X) = \{d\}$.

Consider now the subset $X' = X \backslash \{d\}$. A contraction of the profile $\vec{P}_X$ onto a set $X'$, i.e., $\vec{P}_{X'}$, looks as

| $P_1$ | $P_2$ | $P_3$ | $P_4$ | $P_5$ |
|---|---|---|---|---|
| a | a | c | b | b |
| c | c | b | a | a |
| b | b | a | c | c |

According to the Black procedure the alternative $a$ will be chosen, i.e., $C(\vec{P}_{X'}, X') = \{a\}$.

Finally, let us consider the subset $X'' = X \backslash \{b, c\}$. A contraction of the profile $\vec{P}_X$ onto a set $X''$, i.e., $\vec{P}_{X''}$, looks as

| $P_1$ | $P_2$ | $P_3$ | $P_4$ | $P_5$ |
|---|---|---|---|---|
| a | a | a | d | d |
| d | d | d | a | a |

According to the rule the alternative $a$ will be chosen, i.e., $C(\vec{P}_{X''}, X'') = \{a\}$.

Then $C(\vec{P}_{X'}, X') \cap C(\vec{P}_{X''}, X'') = \{a\} \nsubseteq C(\vec{P}_X, X)$. Thus, the Black procedure does not satisfy the condition **C**.

### 8.3 Outcast condition (O)
Consider the example from paragraph 8.1 of this section. Then $C(\vec{P}_{X'}, X') \neq C(\vec{P}_X, X)$. Thus, the Black procedure does not satisfy the condition **O**.

### 8.4 Arrow's choice axiom (ACA)
The Black procedure does not satisfy the condition **ACA** as the condition **H** is not satisfied.

### 8.5 Monotonicity condition 1
Let $a \in C(\vec{P}_X, X)$. Consider two possible situations:
1. None of alternatives is a Condorcet winner. Thus, the Borda rule is used for this case. Then $\forall x \in X \setminus \{a\}\, r(a, \vec{P}_X) \geq r(x, \vec{P}_X)$. Consider now a profile $\overrightarrow{P'_X}$, which differs from the profile $\vec{P}_X$ only by improved position of the alternative $a$. Then $r(a, \overrightarrow{P'_X}) > r(a, \vec{P}_X)$ and $r(x, \overrightarrow{P'_X}) \leq r(x, \vec{P}_X)$. Thus, $r(a, \overrightarrow{P'_X}) > r(x, \overrightarrow{P'_X})$ and $a \in C(\overrightarrow{P'_X}, X)$.
2. The alternative $a$ is a Condorcet winner, i.e., $a \in CW(\vec{P}_X, X)$. Then $\nexists x \in X, x \mu a$. Since for any profile $\overrightarrow{P'_X}$ which differs from the profile $\vec{P}_X$ only by improved position of the alternative $a$ the relative comparison of the alternative $a$ and any other alternative remains unchanged, the alternative $a$ will be in choice, i.e., $a \in C(\overrightarrow{P'_X}, X)$.

Then $a \in C(\vec{P}_X, X)$ and $a \in C(\overrightarrow{P'_X}, X)$. Thus, the Black procedure satisfies the Monotonicity condition 1.

### 8.6 The Monotonicity condition 2
Let $X = \{a, b, c, d, e\}$ and the profile $\vec{P}_X$ is the following

| $P_1$ | $P_2$ | $P_3$ | $P_4$ | $P_5$ | $P_6$ | $P_7$ | $P_8$ | $P_9$ | $P_{10}$ | $P_{11}$ | $P_{12}$ | $P_{13}$ |
|---|---|---|---|---|---|---|---|---|---|---|---|---|
| e | b | b | a | c | c | c | a | e | a | e | d | d |
| a | a | a | b | d | d | d | b | b | b | b | e | b |
| b | c | c | c | e | e | e | c | a | c | a | a | c |
| c | e | e | e | a | b | a | e | c | e | c | c | e |
| d | d | d | d | b | a | b | d | d | d | d | b | a |

In this example none of alternatives is a Condorcet winner ($a\mu b, a\mu c, a\mu d, e\mu a, c\mu e$). Thus, the Borda rule is used for this case. Let us calculate the Borda count for each alternative: $r(a, \vec{P}_X) = r(b, \vec{P}_X) = 29, r(c, \vec{P}_X) = 28, r(d, \vec{P}_X) = 17, r(e, \vec{P}_X) = 27$. According to the Black procedure the alternatives $a,b$ will be chosen, i.e., $C(\vec{P}_X, X) = \{a, b\}$.

Consider now the subset $X' = X \setminus \{a\}$. A contraction of the profile $\vec{P}_X$ onto a set $X'$, i.e., $\vec{P}_{X'}$, looks as

| $P_1$ | $P_2$ | $P_3$ | $P_4$ | $P_5$ | $P_6$ | $P_7$ | $P_8$ | $P_9$ | $P_{10}$ | $P_{11}$ | $P_{12}$ | $P_{13}$ |
|---|---|---|---|---|---|---|---|---|---|---|---|---|
| e | b | b | b | c | c | c | b | e | b | e | d | d |
| b | c | c | c | d | d | d | c | b | c | b | e | b |
| c | e | e | e | e | e | e | e | c | e | c | c | c |
| d | d | d | d | b | b | b | d | d | d | d | b | e |

In this example none of alternatives is a Condorcet winner ($b\mu c, b\mu d, e\mu b, c\mu e$). Thus, the Borda rule is used for this case. Let us calculate the Borda count for each alternative: $r(b, \vec{P}_{X'}) = 23, r(c, \vec{P}_{X'}) = 24, r(d, \vec{P}_{X'}) = 12, r(e, \vec{P}_{X'}) = 19$. According to the Black procedure the alternative $c$ will be chosen, i.e., $C(\vec{P}_{X'}, X') = \{c\}$.

Finally, let us consider the subset $X'' = X \setminus \{b\}$. A contraction of the profile $\vec{P}_X$ onto a set $X''$, i.e., $\vec{P}_{X''}$, looks as

| $P_1$ | $P_2$ | $P_3$ | $P_4$ | $P_5$ | $P_6$ | $P_7$ | $P_8$ | $P_9$ | $P_{10}$ | $P_{11}$ | $P_{12}$ | $P_{13}$ |
|---|---|---|---|---|---|---|---|---|---|---|---|---|
| e | a | a | a | c | c | c | a | e | a | e | d | d |
| a | c | c | c | d | d | d | c | a | c | a | e | c |
| c | e | e | e | e | e | e | e | c | e | c | a | e |
| d | d | d | d | a | a | a | d | d | d | d | c | a |

In this example none of alternatives is a Condorcet winner ($a\mu c, a\mu d, e\mu a, c\mu e$). Thus, the Borda rule is used for this case. Let us calculate the Borda count for each alternative: $r(a, \vec{P}_{X''}) = 22, r(c, \vec{P}_{X''}) = 24, r(d, \vec{P}_{X''}) = 12, r(e, \vec{P}_{X''}) = 20$. According to the rule the alternative $c$ will be chosen, i.e., $C(\vec{P}_{X''}, X'') = \{c\}$.

Then $\{a, b\} \in C(\vec{P}_X, X)$, $\{b\} \notin C(\vec{P}_{X'}, X')$ and $\{a\} \notin C(\vec{P}_{X''}, X'')$. Thus, the Black procedure does not satisfy the Monotonicity condition 2.

### 8.7 The strict monotonicity condition

Let $X = \{a, b, c, d\}$ and the profile $\vec{P}_X$ is the following

| $P_1$ | $P_2$ | $P_3$ | $P_4$ | $P_5$ | $P_6$ | $P_7$ | $P_8$ | $P_9$ |
|---|---|---|---|---|---|---|---|---|
| d | a | a | a | a | c | b | b | b |
| a | d | d | c | c | d | c | c | c |
| b | b | b | d | d | b | d | d | d |
| c | c | c | b | b | a | a | a | a |

In this example none of alternatives is a Condorcet winner ($a\mu b, a\mu c, d\mu a, b\mu c, d\mu b, c\mu d$). Thus, the Borda rule is used for this case. Let us calculate the Borda count for each alternative: $r(a, \vec{P}_X) = 14, r(b, \vec{P}_X) = r(c, \vec{P}_X) = 13, r(d, \vec{P}_X) = 14$. According to the Black procedure the alternatives $a, d$ will be chosen, i.e., $C(\vec{P}_X, X) = \{a, d\}$.

Consider now a profile $\vec{P'_X}$, which differs from the profile $\vec{P}_X$ only by improved position of the alternative $c$ in $P'_4$:

| $P'_1$ | $P'_2$ | $P'_3$ | $P'_4$ | $P'_5$ | $P'_6$ | $P'_7$ | $P'_8$ | $P'_9$ |
|---|---|---|---|---|---|---|---|---|
| d | a | a | c | a | c | b | b | b |
| a | d | d | a | c | d | c | c | c |
| b | b | b | d | d | b | d | d | d |
| c | c | c | b | b | a | a | a | a |

In this example none of alternatives is a Condorcet winner ($a\mu b, c\mu a, d\mu a, b\mu c, d\mu b, c\mu d$). Thus, the Borda rule is used for this case. Let us calculate the Borda count for each alternative: $r(a, \vec{P'_X}) = r(b, \vec{P'_X}) = 13, r(c, \vec{P'_X}) = 14, r(d, \vec{P'_X}) = 14$. According to the rule the alternatives $c, d$ will be chosen, i.e., $C(\vec{P'_X}, X) = \{c, d\}$.

$$C\left(\vec{P}'_X, X\right) \neq \begin{bmatrix} C(\vec{P}_X, X) \text{ or} \\ \{c\} \text{ or} \\ C(\vec{P}_X, X) \cup \{c\}. \end{bmatrix}$$

Thus, the Black procedure does not satisfy the strict monotonicity condition.

### *8.8 The non-compensatory condition*

The Black procedure does not satisfy the non-compensatory condition (see paragraph 1.8 of this section).

## 9. The inverse Borda rule
### *9.1 Heredity condition (H)*

Let $X = \{a, b, c, d, e\}$ and the profile $\vec{P}_X$ is the following

| P₁ | P₂ | P₃ | P₄ | P₅ | P₆ |
|---|---|---|---|---|---|
| a | b | c | a | d | c |
| b | d | a | b | c | a |
| d | c | b | d | b | b |
| c | a | d | c | a | d |

Step 1. Let us calculate the Borda count for each alternative: $r(a, \vec{P}_X) = 10$, $r(b, \vec{P}_X) = 10$, $r(c, \vec{P}_X) = 9$, $r(d, \vec{P}_X) = 7$. According to the Inverse Borda rule the alternative $d$ is omitted.

Step 2. Consider the subset $X' = X \setminus \{d\}$. A contraction of the profile $\vec{P}_X$ onto a set $X'$, i.e., $\vec{P}_{X'}$, looks as

| P₁ | P₂ | P₃ | P₄ | P₅ | P₆ | Borda count |
|---|---|---|---|---|---|---|
| a | b | c | a | c | c | $r(a, \vec{P}_{X'}) = 6$ |
| b | c | a | b | b | a | $r(b, \vec{P}_{X'}) = 5$ |
| c | a | b | c | a | b | $r(c, \vec{P}_{X'}) = 7$ |

According to the Inverse Borda rule the alternative $b$ is omitted.

Step 3. Consider the subset $X'' = X' \setminus \{b\}$. A contraction of the profile $\vec{P}_{X'}$ onto a set $X''$, i.e., $\vec{P}_{X''}$, looks as

| P₁ | P₂ | P₃ | P₄ | P₅ | P₆ | Borda count |
|---|---|---|---|---|---|---|
| a | c | c | a | c | c | $r(a, \vec{P}_{X''}) = 2$ |
| c | a | a | c | a | a | $r(c, \vec{P}_{X''}) = 4$ |

According to the Inverse Borda rule the alternative $c$ will be chosen, i.e., $C(\vec{P}, X) = \{c\}$.

Consider now the subset $X''' = X \setminus \{a, b\} = \{c, d\}$. A contraction of the profile $\vec{P}_X$ onto a set $X'''$, i.e., $\vec{P}_{X'''}$, looks as

| P₁ | P₂ | P₃ | P₄ | P₅ | P₆ | Borda count |
|---|---|---|---|---|---|---|
| d | d | c | d | d | c | $r(c, \vec{P}_{X'''}) = 2$ |
| c | c | d | c | c | d | $r(d, \vec{P}_{X'''}) = 4$ |

According to the rule the alternative $d$ will be chosen, i.e., $C(\vec{P}_{X'''}, X''') = \{d\}$.

Then $C(\vec{P}_{X'''}, X''') \not\supseteq C(\vec{P}_X, X) \cap X'''$. Thus, the Inverse Borda rule does not satisfy the condition **H**.

### 9.2 Concordance condition (C)

Let $X = \{a, b, c, d\}$ and the profile $\vec{P}_X$ is the following

| $P_1$ | $P_2$ | $P_3$ | $P_4$ | $P_5$ |
|---|---|---|---|---|
| a | a | b | b | c |
| c | b | c | d | d |
| d | c | d | a | a |
| b | d | a | c | b |

Step 1. Let us calculate the Borda count for each alternative: $r(a, \vec{P}_{\{a,b,c,d\}}) = r(b, \vec{P}_{\{a,b,c,d\}}) = r(c, \vec{P}_{\{a,b,c,d\}}) = 8$, $r(d, \vec{P}_{\{a,b,c,d\}}) = 6$. The alternative $d$ is omitted.

Step 2. Let us calculate the Borda count for each alternative: $r(a, \vec{P}_{\{a,b,c\}}) = 6, r(b, \vec{P}_{\{a,b,c\}}) = 5, r(c, \vec{P}_{\{a,b,c,d\}}) = 4$. The alternative $c$ is omitted.

Step 3. Let us calculate the Borda count for each alternative: $r(a, \vec{P}_{\{a,b\}}) = 3, r(b, \vec{P}_{\{a,b\}}) = 2$.

According to the Inverse Borda rule the alternative $a$ will be chosen, i.e., $C(\vec{P}_X, X) = \{a\}$.

Consider now the subset $X' = X \setminus \{a\} = \{b, c, d\}$. A contraction of the profile $\vec{P}_X$ onto a set $X'$, i.e., $\vec{P}_{X'}$, looks as

| $P_1$ | $P_2$ | $P_3$ | $P_4$ | $P_5$ |
|---|---|---|---|---|
| c | b | b | b | c |
| d | c | c | d | d |
| b | d | d | c | b |

Step 1. Let us calculate the Borda count for each alternative: $r(b, \vec{P}_{\{b,c,d\}}) = r(c, \vec{P}_{\{b,c,d\}}) = 6$, $r(d, \vec{P}_{\{b,c,d\}}) = 3$. According to the rule the alternative $d$ is omitted.

Step 2. Let us calculate the Borda count for each alternative: $r(b, \vec{P}_{\{b,c\}}) = 3, r(c, \vec{P}_{\{b,c\}}) = 2$.

According to the Inverse Borda rule the alternative $b$ will be chosen, i.e., $C(\vec{P}_{X'}, X') = \{b\}$.

Finally, let us consider the subset $X'' = X \setminus \{c\}$. A contraction of the profile $\vec{P}_X$ onto a set $X''$, i.e., $\vec{P}_{X''}$, looks as

| $P_1$ | $P_2$ | $P_3$ | $P_4$ | $P_5$ |
|---|---|---|---|---|
| a | a | b | b | d |
| d | b | d | d | a |
| b | d | a | a | b |

Step 1. Let us calculate the Borda count for each alternative: $r(a, \vec{P}_{X''}) = r(b, \vec{P}_{X''}) = r(d, \vec{P}_{X''}) = 5$. According to the rule the alternatives $a,b,d$ will be chosen, i.e., $C(\vec{P}_{X''}, X'') = \{a, b, d\}$.

Then $C(\vec{P}_{X'}, X') \cap C(\vec{P}_{X''}, X'') = \{b\} \not\subseteq C(\vec{P}_X, X)$. Thus, the Inverse Borda rule does not satisfy the condition **C**.

### 9.3 Outcast condition (O)

Consider the example from paragraph 9.1 of this section. Then $C(\vec{P}_{X'}, X') \neq C(\vec{P}_X, X)$, consequently, the Inverse Borda rule does not satisfy the condition **O**.

### 9.4 Arrow's choice axiom (ACA)

The Inverse Borda rule does not satisfy the condition **ACA** as the condition **H** is not satisfied.

## 9.5 Monotonicity condition 1

Let $X = \{a, b, c, d\}$ and the profile $\vec{P}_X$ is the following

| $P_1$ | $P_2$ | $P_3$ | $P_4$ | $P_5$ | $P_6$ | $P_7$ |
|---|---|---|---|---|---|---|
| b | c | b | d | b | a | c |
| a | a | a | a | c | d | d |
| c | d | d | c | d | c | b |
| d | b | c | b | a | b | a |

Step 1. Let us calculate the Borda count for each alternative: $r(a, \vec{P}_X) = 11$, $r(b, \vec{P}_X) = 10$, $r(c, \vec{P}_X) = 11$, $r(d, \vec{P}_X) = 10$. According to the Inverse Borda rule the alternatives $b$ and $d$ are omitted.

Step 2. Consider the subset $X' = X \setminus \{b, d\}$. A contraction of the profile $\vec{P}_X$ onto a set $X'$, i.e., $\vec{P}_{X'}$, looks as

| $P_1$ | $P_2$ | $P_3$ | $P_4$ | $P_5$ | $P_6$ | $P_7$ | | Borda count |
|---|---|---|---|---|---|---|---|---|
| a | c | a | a | c | a | c | | $r(a, \vec{P}_{X'}) = 4$ |
| c | a | c | c | a | c | a | | $r(c, \vec{P}_{X'}) = 3$ |

The alternative $c$ is omitted and the alternative $a$ will be chosen, i.e., $C(\vec{P}, X) = \{a\}$.

Consider now a profile $\vec{P'_X}$, which differs from the profile $\vec{P}_X$ only by improved position of the alternative $a$ in $P'_5$.

| $P'_1$ | $P'_2$ | $P'_3$ | $P'_4$ | $P'_5$ | $P'_6$ | $P'_7$ |
|---|---|---|---|---|---|---|
| b | c | b | d | b | a | c |
| a | a | a | a | a | d | d |
| c | d | d | c | c | c | b |
| d | b | c | b | d | b | a |

Step 1. Let us calculate the Borda count for each alternative: $r(a, \vec{P'_X}) = 13$, $r(b, \vec{P'_X}) = 10$, $r(c, \vec{P'_X}) = 10$, $r(d, \vec{P'_X}) = 9$. The alternative $d$ is omitted.

Step 2. Consider the subset $X' = X \setminus \{d\}$. A contraction of the profile $\vec{P'_X}$ onto a set $X'$, i.e., $\vec{P'}_{X'}$, looks as

| $P'_1$ | $P'_2$ | $P'_3$ | $P'_4$ | $P'_5$ | $P'_6$ | $P'_7$ | | Borda count |
|---|---|---|---|---|---|---|---|---|
| b | c | b | a | b | a | c | | $r(a, \vec{P'}_{X'}) = 8$ |
| a | a | a | c | a | c | b | | $r(b, \vec{P'}_{X'}) = 7$ |
| c | b | c | b | c | b | a | | $r(c, \vec{P'}_{X'}) = 6$ |

The alternative $c$ is omitted.

Step 3. Consider the subset $X'' = X' \setminus \{c\}$. A contraction of the profile $\vec{P'_X}$ onto a set $X''$, i.e., $\vec{P'}_{X''}$, looks as

| $P'_1$ | $P'_2$ | $P'_3$ | $P'_4$ | $P'_5$ | $P'_6$ | $P'_7$ | | Borda count |
|---|---|---|---|---|---|---|---|---|
| b | a | b | a | b | a | b | | $r(a, \vec{P'}_{X''}) = 3$ |
| a | b | a | b | a | b | a | | $r(b, \vec{P'}_{X''}) = 4$ |

The alternative $a$ is omitted. Thus, according to the rule the alternative $b$ will be chosen, i.e., $C(\vec{P'_X}, X) = \{b\}$.

Then $\{a\} \in C(\vec{P}_X, X)$, $\{a\} \notin C(\vec{P'_X}, X)$. Thus, the Inverse Borda rule does not satisfy the Monotonicity condition 1.

### 9.6 The Monotonicity condition 2

Let $X = \{a, b, c, d\}$ and the profile $\vec{P}_X$ is the following

| $P_1$ | $P_2$ | $P_3$ | $P_4$ |
|---|---|---|---|
| a | a | b | c |
| d | d | c | b |
| b | c | d | d |
| c | b | a | a |

Let us calculate the Borda count for each alternative: $r(a, \vec{P}_X) = r(b, \vec{P}_X) = r(c, \vec{P}_X) = r(d, \vec{P}_X) = 6$. According to the Inverse Borda rule the alternatives $a,b,c,d$ will be chosen, i.e., $C(\vec{P}_X, X) = \{a, b, c, d\}$.

Consider now the subset $X' = X \setminus \{a\}$. A contraction of the profile $\vec{P}_X$ onto a set $X'$, i.e., $\vec{P}_{X'}$, looks as

| $P_1$ | $P_2$ | $P_3$ | $P_4$ |
|---|---|---|---|
| d | d | b | c |
| b | c | c | b |
| c | b | d | d |

Let us calculate the Borda count for each alternative: $r(b, \vec{P}_{X'}) = r(c, \vec{P}_{X'}) = 4, r(d, \vec{P}_{X'}) = 2$. According to the rule the alternatives $b,c$ will be chosen, i.e., $C(\vec{P}_{X'}, X') = \{b, c\}$.

Then $\{d\} \in C(\vec{P}_X, X)$, $\{d\} \notin C(\vec{P}_{X'}, X')$. Thus, the Inverse Borda rule does not satisfy the Monotonicity condition 2.

### 9.7 The strict monotonicity condition

The Inverse Borda rule does not satisfy the strict monotonicity condition (see the proof in paragraph 2.7 of this section).

### 9.8 The non-compensatory condition

The Inverse Borda rule does not satisfy the non-compensatory condition (see paragraph 1.8 of this section).

## 10. The Nanson rule
### 10.1 Heredity condition (H)

Let $X = \{a, b, c, d, e\}$ and the profile $\vec{P}_X$ is the following

| $P_1$ | $P_2$ | $P_3$ | $P_4$ | $P_5$ | $P_6$ |
|---|---|---|---|---|---|
| a | b | c | a | d | c |
| b | d | a | b | c | a |
| d | c | b | d | b | b |
| c | a | e | c | a | d |
| e | e | d | e | e | e |

Step 1. Let us calculate the Borda count for each alternative: $r(a, \vec{P}_X) = 16$, $r(b, \vec{P}_X) = 16$, $r(c, \vec{P}_X) = 15$, $r(d, \vec{P}_X) = 12$, $r(e, \vec{P}_X) = 1$, $\bar{r}(\vec{P}_X, X) = 12$. According to the Nanson rule the alternatives $d,e$ are omitted.

Step 2. Consider the subset $X' = X \setminus \{d, e\}$. A contraction of the profile $\vec{P}_X$ onto a set $X'$, i.e., $\vec{P}_{X'}$, looks as

| $P_1$ | $P_2$ | $P_3$ | $P_4$ | $P_5$ | $P_6$ | Borda count |
|---|---|---|---|---|---|---|
| a | b | c | a | c | c | $r(a, \vec{P}_{X'}) = 6$ |
| b | c | a | b | b | a | $r(b, \vec{P}_{X'}) = 5$ |
| c | a | b | c | a | b | $r(c, \vec{P}_{X'}) = 7$ |

The average Borda count is equal to 6, i.e., $\bar{r}(\vec{P}_{X'}, X') = 6$. According to the Nanson rule the alternative $b$ is omitted.

Step 3. Consider the subset $X'' = X' \setminus \{b\}$. A contraction of the profile $\vec{P}_X$ onto a set $X''$, i.e., $\vec{P}_{X''}$, looks as

| $P_1$ | $P_2$ | $P_3$ | $P_4$ | $P_5$ | $P_6$ | Borda count |
|---|---|---|---|---|---|---|
| a | c | c | a | c | c | $r(a, \vec{P}_{X''}) = 2$ |
| c | a | a | c | a | a | $r(c, \vec{P}_{X''}) = 4$ |

The average Borda count is equal to 6, i.e., $\bar{r}(\vec{P}_{X''}, X'') = 6$. Thus, according to the Nanson rule the alternative $c$ will be chosen, i.e., $C(\vec{P}_X, X) = \{c\}$.

Consider now the subset $X''' = X \setminus \{a, b\} = \{a, d, e\}$. A contraction of the profile $\vec{P}$ onto a set $X'''$, i.e., $\vec{P}_{X'''}$, looks as

| $P_1$ | $P_2$ | $P_3$ | $P_4$ | $P_5$ | $P_6$ | Borda count |
|---|---|---|---|---|---|---|
| d | d | c | d | d | c | $r(c, \vec{P}_{X'''}) = 2$ |
| c | c | d | c | c | d | $r(d, \vec{P}_{X'''}) = 4$ |

According to the rule the alternative $d$ will be chosen, i.e., $C(\vec{P}_{X'''}, X''') = \{d\}$.

Then $C(\vec{P}_{X'''}, X''') \not\supseteq C(\vec{P}_X, X) \cap X'''$. Thus, the Nanson rule does not satisfy the condition **H**.

*10.2 Concordance condition (C)*

Let $X = \{a, b, c, d\}$ and the profile $\vec{P}_X$ is the following

| $P_1$ | $P_2$ | $P_3$ | $P_4$ | $P_5$ |
|---|---|---|---|---|
| a | a | b | b | c |
| c | b | c | d | d |
| b | c | d | a | a |
| d | d | a | c | b |

Step 1. Let us calculate the Borda count for each alternative: $r(a, \vec{P}_{\{a,b,c,d\}}) = 8$, $r(b, \vec{P}_{\{a,b,c,d\}}) = 9$, $r(c, \vec{P}_{\{a,b,c,d\}}) = 8$, $r(d, \vec{P}_{\{a,b,c,d\}}) = 5$, $\bar{r}(\vec{P}_X, X) = 7.5$. According to the Nanson rule the alternative $d$ is omitted.

Step 2. Let us calculate the Borda count for each alternative: $r(a, \vec{P}_{\{a,b,c\}}) = 6, r(b, \vec{P}_{\{a,b,c\}}) = 5$, $r(c, \vec{P}_{\{a,b,c\}}) = 4, \bar{r}(\vec{P}_{\{a,b,c\}}, \{a,b,c\}) = 5.$. According to the Nanson rule $C(\vec{P}_X, X) = \{a\}$.

Consider now the subset $X' = X \setminus \{a\} = \{b, c, d\}$. A contraction of the profile $\vec{P}_X$ onto a set $X'$, i.e., $\vec{P}_{X'}$, looks as

| $P_1$ | $P_2$ | $P_3$ | $P_4$ | $P_5$ |
|---|---|---|---|---|
| c | b | b | b | c |
| b | c | c | d | d |
| d | d | d | c | b |

Step 1. Let us calculate the Borda count for each alternative: $r(b, \vec{P}_{\{b,c,d\}}) = 7, r(c, \vec{P}_{\{b,c,d\}}) = 5, r(d, \vec{P}_{\{b,c,d\}}) = 2, \bar{r}(\vec{P}_{\{b,c,d\}}, \{b,c,d\}) \approx 4{,}7$. According to the Nanson rule the alternative $d$ is omitted.

Step 2. Let us calculate the Borda count for each alternative: $r(b, \vec{P}_{\{b,c\}}) = 3, r(c, \vec{P}_{\{b,c\}}) = 2, \bar{r}(\vec{P}_{\{b,c\}}, \{b,c\}) = 2{,}5$. According to the Nanson rule $C(\vec{P}_{X'}, X') = \{b\}$.

Finally, let us consider the subset $X'' = X\setminus\{c\}$. A contraction of the profile $\vec{P}_X$ onto a set $X''$, i.e., $\vec{P}_{X''}$, looks as

| $P_1$ | $P_2$ | $P_3$ | $P_4$ | $P_5$ |
|---|---|---|---|---|
| a | a | b | b | d |
| b | b | d | d | a |
| d | d | a | a | b |

Step 1. Let us calculate the Borda count for each alternative: $r(a, \vec{P}_{X''}) = 5, r(b, \vec{P}_{X''}) = 6, r(d, \vec{P}_{X''}) = 4$. According to the Nanson rule the alternative $b$ will be chosen, i.e., $C(\vec{P}_{X''}, X'') = \{b\}$.

Then $C(\vec{P}_{X'}, X') \cap C(\vec{P}_{X''}, X'') = \{b\} \nsubseteq C(\vec{P}_X, X)$. Thus, the Nanson rule does not satisfy the condition **C**.

### *10.3 Outcast condition (O)*
Consider the example from paragraph 10.1 of this section. Then $C(\vec{P}_{X'''}, X''') \neq C(\vec{P}_X, X)$, consequently, the Nanson rule does not satisfy the condition **O**.

### *10.4 Arrow's choice axiom (ACA)*
The Nanson rule does not satisfy the condition **ACA** as the condition **H** is not satisfied.

### *10.5 Monotonicity condition 1*
Let $X = \{a, b, c, d\}$ and the profile $\vec{P}_X$ is the following

| $P_1$ | $P_2$ | $P_3$ | $P_4$ | $P_5$ |
|---|---|---|---|---|
| a | c | b | c | b |
| b | a | c | a | a |
| c | b | d | d | d |
| d | d | a | b | c |

Step 1. Let us calculate the Borda count for each alternative: $r(a, \vec{P}_X) = 7$, $r(b, \vec{P}_X) = 7$, $r(c, \vec{P}_X) = 7$, $r(d, \vec{P}_X) = 3$. According to the rule the alternatives $a,b,c$ will be chosen, i.e., $C(\vec{P}, X) = \{a, b, c\}$.

Consider now a profile $\vec{P}'_X$, which differs from the profile $\vec{P}_X$ only by improved position of the alternative $a$ in $P'_5$.

| $P'_1$ | $P'_2$ | $P'_3$ | $P'_4$ | $P'_5$ |
|---|---|---|---|---|
| a | c | b | c | a |
| b | a | c | a | b |
| c | b | d | d | d |
| d | d | a | b | c |

Step 1. Let us calculate the Borda count for each alternative: $r(a, \overrightarrow{P'_X}) = 8$, $r(b, \overrightarrow{P'_X}) = 6$, $r(c, \overrightarrow{P'_X}) = 7$, $r(d, \overrightarrow{P'_X}) = 3$, $\bar{r}(\overrightarrow{P'_X}, X) = 6$. The alternatives $b, d$ are omitted.

Step 2. Consider the subset $X' = X \setminus \{b\}$. A contraction of the profile $\overrightarrow{P'_X}$ onto a set $X'$, i.e., $\overrightarrow{P'}_{X'}$, looks as

| $P'_1$ | $P'_2$ | $P'_3$ | $P'_4$ | $P'_5$ | Borda count |
|---|---|---|---|---|---|
| a | c | c | c | a | $r(a, \overrightarrow{P'}_{X'}) = 2$ |
| c | a | a | a | c | $r(c, \overrightarrow{P'}_{X'}) = 3$ |

The alternative $a$ is omitted. Thus, according to the rule the alternative $c$ will be chosen, i.e., $C(\overrightarrow{P'_X}, X) = \{c\}$.

Then $\{a\} \in C(\vec{P}_X, X)$, $\{a\} \notin C(\overrightarrow{P'_X}, X)$. Thus, the Nanson rule satisfies the Monotonicity condition 1.

### 10.6 The Monotonicity condition 2

The Nanson rule does not satisfy the Monotonicity condition 2 (see paragraph 9.6 of this section).

### 10.7 The strict monotonicity condition

Consider the example from paragraph 9.7. According to the rule the alternative $c$ will be chosen, i.e., $C(\vec{P}_X, X) = \{c\}$. Consider now a profile $\overrightarrow{P'_X}$, which differs from the profile $\vec{P}_X$ only by improved position of the alternative $b$ in $P'_5$:

| $P'_1$ | $P'_2$ | $P'_3$ | $P'_4$ | $P'_5$ | $P'_6$ |
|---|---|---|---|---|---|
| a | b | c | a | b | c |
| b | d | a | b | d | a |
| d | c | b | d | c | b |
| c | a | d | c | a | d |

Step 1. Let us calculate the Borda count for each alternative: $r(a, \overrightarrow{P'_X}) = 10$, $r(b, \overrightarrow{P'_X}) = 12$, $r(c, \overrightarrow{P'_X}) = 8$, $r(d, \overrightarrow{P'_X}) = 6$, $\bar{r}(\overrightarrow{P'_X}, X) = 9$. According to the Nanson rule the alternatives $c$ and $d$ are omitted.

Step 2. Consider the subset $X' = X \setminus \{c, d\}$. A contraction of the profile $\overrightarrow{P'_X}$ onto a set $X'$, i.e., $\overrightarrow{P'}_{X'}$, looks as

| $P'_1$ | $P'_2$ | $P'_3$ | $P'_4$ | $P'_5$ | $P'_6$ | Borda count |
|---|---|---|---|---|---|---|
| a | b | a | a | b | a | $r(a, \overrightarrow{P'}_{X'}) = 4$ |
| b | a | b | b | a | b | $r(b, \overrightarrow{P'}_{X'}) = 2$ |

According to the rule the alternative $a$ will be chosen, i.e., $C(\overrightarrow{P'_X}, X) = \{a\}$.

$$C(\overrightarrow{P'_X}, X) \neq \begin{bmatrix} C(\vec{P}_X, X) \text{ or} \\ \{b\} \text{ or} \\ C(\vec{P}_X, X) \cup \{b\}. \end{bmatrix}$$

Thus, the Nanson rule does not satisfy the strict monotonicity condition.

### 10.8 The non-compensatory condition

The Nanson rule does not satisfy the non-compensatory condition (see paragraph 1.8 of this section).

# 11. The Coombs procedure
## 11.1 Heredity condition (H)
Let $X = \{a, b, c, d\}$ and the profile $\vec{P}_X$ is the following

| $P_1$ | $P_2$ | $P_3$ | $P_4$ | $P_5$ | $P_6$ | $P_7$ |
|---|---|---|---|---|---|---|
| c | c | c | d | d | a | a |
| a | d | d | b | a | c | c |
| d | b | b | a | b | b | b |
| b | a | a | c | c | d | d |

According to the Coombs procedure the alternative $b$ will be chosen, i.e., $C(\vec{P}_X, X) = \{b\}$.
Consider now the subset $X' = X\setminus\{d\}$. A contraction of the profile $\vec{P}_X$ onto a set $X'$, i.e., $\vec{P}_{X'}$, looks as

| $P_1$ | $P_2$ | $P_3$ | $P_4$ | $P_5$ | $P_6$ | $P_7$ |
|---|---|---|---|---|---|---|
| c | c | c | b | a | a | a |
| d | b | b | a | b | c | c |
| b | a | a | c | c | b | b |

$n^-(a, \vec{P}_{X'}) = 2$, $n^-(b, \vec{P}_{X'}) = 3$, $n^-(c, \vec{P}_{X'}) = 2$. The alternative $b$ is omitted.
Thus, according to the rule the alternative $a$ will be chosen, i.e., $C(\vec{P}_{X'}, X') = \{a\}$.
Since $C(\vec{P}_{X'}, X') \not\supseteq C(\vec{P}_X, X) \cap X'$, the Coombs procedure does not satisfy the condition **H**.

## 11.2 Concordance condition (C)
Consider the previous example. According to the Coombs procedure $C(\vec{P}_X, X) = \{b\}$, $C(\vec{P}_{X'}, X') = \{a\}$.
Consider now the subset $X'' = X\setminus\{b\}$. A contraction of the profile $\vec{P}_X$ onto a set $X''$, i.e., $\vec{P}_{X''}$, looks as

| $P_1$ | $P_2$ | $P_3$ | $P_4$ | $P_5$ | $P_6$ | $P_7$ |
|---|---|---|---|---|---|---|
| c | c | c | d | d | a | a |
| a | d | d | a | a | c | c |
| d | a | a | c | c | d | d |

$n^-(a, \vec{P}_{X''}) = 2$, $n^-(c, \vec{P}_{X''}) = 2$, $n^-(d, \vec{P}_{X''}) = 3$. The alternative $d$ is omitted.
Thus, according to the rule the alternatives $a$ and $c$ will be chosen, i.e., $C(\vec{P}_{X''}, X'') = \{a, c\}$. Then $C(\vec{P}_{X'}, X') \cap C(\vec{P}_{X''}, X'') = \{a\} \nsubseteq C(\vec{P}_X, X)$. Thus, the Coombs procedure does not satisfy the condition **C**.

## 11.3 Outcast condition (O)
Consider the example from paragraph 11.1 of this section. Then $C(\vec{P}_{X'}, X') \neq C(\vec{P}_X, X)$, consequently, the Coombs procedure does not satisfy the condition **O**.

## 11.4 Arrow's choice axiom (ACA)
The Coombs procedure does not satisfy the condition **ACA** as the condition **H** is not satisfied.

*11.5 Monotonicity condition 1*

Let $X = \{a, b, c\}$ and the profile $\vec{P}_X$ is the following

| $P_1$ | $P_2$ | $P_3$ | $P_4$ | $P_5$ | $P_6$ | $P_7$ | $P_8$ |
|---|---|---|---|---|---|---|---|
| a | c | c | c | c | a | a | a |
| c | a | b | b | b | b | b | b |
| b | b | a | a | a | c | c | c |

According to the rule the alternative $b$ will be chosen, i.e., $C(\vec{P}, X) = \{b\}$.

Suppose now that the position of the alternative $b$ was improved in $P_1'$ while the relative comparison of any pair of other alternatives remained unchanged.

| $P_1$ | $P_2$ | $P_3$ | $P_4$ | $P_5$ | $P_6$ | $P_7$ | $P_8$ |
|---|---|---|---|---|---|---|---|
| a | c | c | c | c | a | a | a |
| b | a | b | b | b | b | b | b |
| c | b | a | a | a | c | c | c |

$n^-(a, \vec{P}_X) = 3$, $n^-(b, \vec{P}_X) = 1$, $n^-(c, \vec{P}_X) = 4$. According to the Coombs procedure the alternative $c$ is omitted. Next, consider the subset $X' = X\setminus\{c\}$. A contraction of the profile $\vec{P}_X'$ onto a set $X'$, i.e., $\vec{P'}_{X'}$, looks as

| $P_1'$ | $P_2'$ | $P_3'$ | $P_4'$ | $P_5'$ | $P_6'$ | $P_7'$ | $P_8'$ |
|---|---|---|---|---|---|---|---|
| a | a | b | b | b | a | a | a |
| b | b | a | a | a | b | b | b |

According to the rule the alternative $a$ will be chosen, i.e., $C\left(\vec{P_X'}, X\right) = C\left(\vec{P_{X'}'}, X'\right) = \{a\}$.

Then $\{b\} \in C(\vec{P}_X, X)$, $\{b\} \notin C\left(\vec{P_X'}, X\right)$. Thus, the Coombs procedure does not satisfy the Monotonicity condition 1.

*11.6 The Monotonicity condition 2*

Let $X = \{a, b, c\}$ and the profile $\vec{P}_X$ is the following

| $P_1$ | $P_2$ | $P_3$ | $P_4$ | $P_5$ | $P_6$ |
|---|---|---|---|---|---|
| c | c | c | a | a | a |
| d | b | d | d | b | b |
| a | d | b | b | d | c |
| b | a | a | c | c | d |

According to the Coombs procedure the alternatives $b$ and $d$ will be chosen, i.e., $C(\vec{P}_X, X) = \{b, d\}$.

Consider the subset $X' = X\setminus\{b\}$. A contraction of the profile $\vec{P}_X$ onto a set $X'$, i.e., $\vec{P}_{X'}$, looks as

| $P_1$ | $P_2$ | $P_3$ | $P_4$ | $P_5$ | $P_6$ |
|---|---|---|---|---|---|
| c | c | c | a | a | a |
| d | d | d | d | d | c |
| a | a | a | c | c | d |

$n^-(a, \vec{P}_{X'}) = 3$, $n^-(c, \vec{P}_{X'}) = 2$, $n^-(d, \vec{P}_{X'}) = 1$. According to the Coombs procedure the alternative $a$ is omitted. Next, consider the subset $X'' = X'\setminus\{a\}$. A contraction of the profile $\vec{P}_{X'}$ onto a set $X''$, i.e., $\vec{P}_{X''}$, looks as

| $P_1$ | $P_2$ | $P_3$ | $P_4$ | $P_5$ | $P_6$ |
|---|---|---|---|---|---|
| c | c | c | d | d | c |
| d | d | d | c | c | d |

According to the Coombs procedure the alternative $c$ will be chosen, i.e., $C(\vec{P}_{X'}, X') = C(\vec{P}_{X''}, X'') = \{c\}$.

Finally, consider the subset $X''' = X \setminus \{d\}$. A contraction of the profile $\vec{P}_X$ onto a set $X'$, i.e., $\vec{P}_{X'}$, looks as

| $P_1$ | $P_2$ | $P_3$ | $P_4$ | $P_5$ | $P_6$ |
|---|---|---|---|---|---|
| c | c | c | a | a | a |
| a | b | b | b | b | b |
| b | a | a | c | c | c |

$n^-(a, \vec{P}_{X'''}) = 2$, $n^-(b, \vec{P}_{X'''}) = 1$, $n^-(c, \vec{P}_{X'''}) = 3$. According to the Coombs procedure the alternative $c$ is omitted.

Next, consider the subset $X'''' = X''' \setminus \{c\}$. A contraction of the profile $\vec{P}_{X'''}$ onto a set $X''''$, i.e., $\vec{P}_{X''''}$, looks as

| $P_1$ | $P_2$ | $P_3$ | $P_4$ | $P_5$ | $P_6$ |
|---|---|---|---|---|---|
| a | b | b | a | a | a |
| b | a | a | b | b | b |

According to the rule the alternative $a$ will be chosen, i.e., $C(\vec{P}_{X'''}, X''') = C(\vec{P}_{X''''}, X'''') = \{a\}$. Then $\{b,d\} \in C(\vec{P}_X, X)$, $\{d\} \notin C(\vec{P}_{X \setminus \{b\}}, X \setminus \{b\})$ and $\{b\} \notin C(\vec{P}_{X \setminus \{d\}}, X \setminus \{d\})$. Thus, the Coombs procedure does not satisfy the Monotonicity condition 2.

*11.7 The strict monotonicity condition*
The Coombs procedure does not satisfy the strict monotonicity condition (see paragraph 3.7 of this section).

*11.8 The non-compensatory condition*
The Coombs procedure does not satisfy the non-compensatory condition (see paragraph 2.8 of this section).

## 12. Minimal dominant set
*12.1 Heredity condition (H)*
Let $X = \{a, b, c\}$ and the profile $\vec{P}_X$ is the following

| $P_1$ | $P_2$ | $P_3$ |
|---|---|---|
| a | c | b |
| b | a | c |
| c | b | a |

For this case a matrix of majority relation μ is the following

|   | a | b | c |
|---|---|---|---|
| a | - | 1 | 0 |
| b | 0 | - | 1 |
| c | 1 | 0 | - |

According to the rule the alternatives $a,b,c$ are included in minimal dominant set $Q$. Thus, $C(\vec{P}_X, X) = Q = \{a, b, c\}$.

Consider now the subset $X' = X\setminus\{b\}$. A contraction of the profile $\vec{P}_X$ onto a set $X'$, i.e., $\vec{P}_{X'}$, looks as

| $P_1$ | $P_2$ | $P_3$ |
|---|---|---|
| a | c | c |
| c | a | a |

According to the rule the alternative $c$ will be chosen, i.e., $C(\vec{P}_{X'}, X') = \{c\}$.
Then $C(\vec{P}_{X'}, X') \not\supseteq C(\vec{P}_X, X) \cap X'$. Thus, the condition **H** is not satisfied.

## 12.2 Concordance condition (C)

The condition **C** is satisfied if $\forall X', X'' \in 2^A \rightarrow C(\vec{P}_{X' \cup X''}, X' \cup X'') \supseteq C(\vec{P}_{X'}, X') \cap C(\vec{P}_{X''}, X'')$.
Let us proof it by contradiction.

1. Suppose that $\exists x \in C(\vec{P}_{X'}, X') \cap C(\vec{P}_{X''}, X''): x \notin C(\vec{P}_{X' \cup X''}, X' \cup X'')$. Then $\forall y \in C(\vec{P}_{X' \cup X''}, X' \cup X'') \rightarrow y\mu x$. It is also true that $y \in X'$ and/or $y \in X''$.
2. Let $y \in X'$. Since $x \in C(\vec{P}_{X'}, X')$, $\forall z \in X' \setminus C(\vec{P}_{X'}, X') \rightarrow x\mu z$. Thus, $y \in C(\vec{P}_{X'}, X')$ and $y\mu z$.
3. The alternatives $x$ and $y$ can be chosen from the set $X'$ iff $\exists u \in C(\vec{P}_{X'}, X'): u\mu y$ and $x\mu u$.
4. Since $u\mu y$ and and $y \in C(\vec{P}_{X' \cup X''}, X' \cup X'')$, $u \in C(\vec{P}_{X' \cup X''}, X' \cup X'')$. Thus, $u\mu x$ which is a contradiction. Thus, the assumption 1 is incorrect.

Thus, the condition **C** is satisfied.

## 12.3 Outcast condition (O)
The condition **O** is satisfied (the proof follows from definition of the rule).

## 12.4 Arrow's choice axiom (ACA)
Minimal dominant set does not satisfy the condition **ACA** as the condition **H** is not satisfied.

## 12.5 Monotonicity condition 1
Let $a \in C(\vec{P}_X, X)$. Then $\forall x \in X \setminus C(\vec{P}_X, X) \rightarrow a\mu x$ and $\nexists X' \subseteq C(\vec{P}_X, X) \setminus \{a\}: \forall y \in X'$ and $\forall z \in X \setminus X'\ y\mu z$. Consider now a profile $\vec{P}'_X$, which differs from the profile $\vec{P}_X$ only by improved position of the alternative $a$. Then $a\mu x$ and $\nexists X' \subseteq C(\vec{P}_X, X) \setminus \{a\}: \forall y \in X'$ and $\forall z \in X \setminus X'\ y\mu z$. Consequently, $a \in C(\vec{P}'_X, X)$. Thus, the Monotonicity condition 1 is satisfied.

## 12.6 The Monotonicity condition 2
The Monotonicity condition 2 is not satisfied (see paragraph 12.1 of this section).

## 12.7 The strict monotonicity condition
The strict monotonicity condition is not satisfied (see paragraph 1.7 of this section).

## 12.8 The non-compensatory condition
The non-compensatory condition is not satisfied (see paragraph 1.8 of this section).

## 13. Minimal undominated set

### 13.1 Heredity condition (H)
The condition **H** is not satisfied (see paragraph 12.1 of this section).

### 13.2 Concordance condition (C)
Let $X = \{a, b, c, d\}$ and the profile $\vec{P}_X$ looks as

| $P_1$ | $P_2$ | $P_3$ | $P_4$ | $P_5$ | $P_6$ |
|---|---|---|---|---|---|
| b | a | a | b | a | d |
| a | c | d | a | c | b |
| c | d | b | c | d | c |
| d | b | c | d | b | a |

For this case a matrix of majority relation μ is the following

|   | a | b | c | d |
|---|---|---|---|---|
| **a** | - | 0 | 1 | 1 |
| **b** | 0 | - | 1 | 0 |
| **c** | 0 | 0 | - | 1 |
| **d** | 0 | 1 | 0 | - |

According to the rule the alternative $a$ is included in minimal undominated set $Q$. Thus, $C(\vec{P}_X, X) = Q = \{a\}$. Consider now the subset $X' = X \setminus \{a\}$. A contraction of the profile $\vec{P}_X$ onto a set $X'$, i.e., $\vec{P}_{X'}$, looks as

| $P_1$ | $P_2$ | $P_3$ | $P_4$ | $P_5$ | $P_6$ |
|---|---|---|---|---|---|
| b | c | d | b | c | d |
| c | d | b | c | d | b |
| d | b | c | d | b | c |

According to the rule the alternatives $b,c,d$ make the minimal undominated set. Thus, $C(\vec{P}_{X'}, X') = \{b, c, d\}$.

Finally, consider the subset $X'' = X \setminus \{c, d\}$. A contraction of the profile $\vec{P}_X$ onto a set $X''$, i.e., $\vec{P}_{X''}$, looks as

| $P_1$ | $P_2$ | $P_3$ | $P_4$ | $P_5$ | $P_6$ |
|---|---|---|---|---|---|
| b | a | a | b | a | b |
| a | b | b | a | b | a |

According to the rule the alternatives $a$ and $b$ will be chosen, i.e., $C(\vec{P}_{X''}, X'') = \{a, b\}$.

Then $C(\vec{P}_{X'}, X') \cap C(\vec{P}_{X''}, X'') = \{b\} \nsubseteq C(\vec{P}_X, X)$. Thus, the condition **C** is not satisfied.

### 13.3 Outcast condition (O)
Consider the previous example. Then $C(\vec{P}_X, X) \neq C(\vec{P}_{X''}, X'')$. Thus, the condition **O** is not satisfied.

### 13.4 Arrow's choice axiom (ACA)
Minimal undominated set does not satisfy the condition **ACA** as the condition **H** is not satisfied.

## 13.5 Monotonicity condition 1

Let $a \in C(\vec{P}_X, X)$. Then $\nexists x \in X \setminus C(\vec{P}_X, X) \to x\mu a$ and $\nexists X' \subseteq C(\vec{P}_X, X) \setminus \{a\}: \forall y \in X'$ and $\forall z \in X \setminus X'$ $z\mu y$. Consider now a profile $\vec{P'_X}$, which differs from the profile $\vec{P}_X$ only by improved position of the alternative $a$. Then $\nexists x \in X \setminus C(\vec{P}_X, X) \to x\mu a$ and $\nexists X' \subseteq C(\vec{P}_X, X) \setminus \{a\}: \forall y \in X'$ and $\forall z \in X \setminus X'$ $z\mu y$. Thus, $a \in C\left(\vec{P'_X}, X\right)$. The Monotonicity condition 1 is satisfied.

## 13.6 The Monotonicity condition 2
The Monotonicity condition 2 is not satisfied (see paragraph 12.1 of this section).

## 13.7 The strict monotonicity condition
The strict monotonicity condition is not satisfied (see paragraph 1.7 of this section).

## 13.8 The non-compensatory condition
The non-compensatory condition is not satisfied (see paragraph 1.8 of this section).

# 14. Minimal weakly stable set
## 14.1 Heredity condition (H)
The condition **H** is not satisfied (see paragraph 12.1 of this section).

## 14.2 Concordance condition (C)
Let $X = \{a, b, c, d, e\}$ and a matrix of majority relation μ is the following

|   | a | b | c | d | e |
|---|---|---|---|---|---|
| a | - | 1 | 0 | 0 | 0 |
| b | 0 | - | 1 | 0 | 0 |
| c | 0 | 0 | - | 1 | 1 |
| d | 1 | 0 | 0 | - | 0 |
| e | 0 | 0 | 0 | 1 | - |

According to the rule the alternatives *a* and *c* are included in minimal weakly stable set $Q$. Thus, $C(\vec{P}_X, X) = Q = \{a, c\}$.

Consider now the subset $X' = X \setminus \{e\}$. Then a matrix of majority relation μ is the following

|   | a | b | c | d |
|---|---|---|---|---|
| a | - | 1 | 0 | 0 |
| b | 0 | - | 1 | 0 |
| c | 0 | 0 | - | 1 |
| d | 1 | 0 | 0 | - |

According to the rule the alternatives *a,b,c,d* are included in minimal weakly stable set, i.e., $C(\vec{P}_{X'}, X') = \{a, b, c, d\}$.

Finally, consider the subset $X'' = X \setminus \{a\}$. Then a matrix of majority relation μ is the following

|   | b | c | d | e |
|---|---|---|---|---|
| b | - | 1 | 0 | 0 |
| c | 0 | - | 1 | 1 |
| d | 0 | 0 | - | 0 |
| e | 0 | 0 | 1 | - |

According to the rule the alternative $b$ is included in minimal weakly stable set, i.e., $C(\vec{P}_{X''}, X'') = \{b\}$.

Then $C(\vec{P}_{X'}, X') \cap C(\vec{P}_{X''}, X'') = \{b\} \nsubseteq C(\vec{P}_X, X)$. Thus, the condition **C** is not satisfied.

### 14.3 Outcast condition (O)

Consider the previous example. $C(\vec{P}_{X''}, X'') \neq C(\vec{P}_X, X)$, consequently, the condition **O** is not satisfied.

### 14.4 Arrow's choice axiom (ACA)

Minimal weakly stable set does not satisfy the condition **ACA** as the condition **H** is not satisfied.

### 14.5 Monotonicity condition 1

Let $a \in C(\vec{P}_X, X)$. Then

$$\begin{cases} \forall b \in X'' = (\{a\} \cup X') \to [\exists c \in X \setminus X''\ c\mu b \Rightarrow \exists d \in X''\ d\mu c]\ and \\ \nexists X''' \subseteq X \setminus \{a\}: |X'''| < |X''|\ and\ \forall e \in X''' \to [\exists f \in X \setminus X'''\ f\mu e \Rightarrow \exists g \in X'''\ g\mu e], \end{cases}$$

where $X' \subseteq C(\vec{P}_X, X) \setminus \{a\}: \nexists X'''' \subset X'\ \forall x \in X'''' \cup \{a\} \to [\exists y \in X \setminus \{X'''' \cup \{a\}\}\ y\mu x \Rightarrow \exists z \in X'''' \cup \{a\}\ z\mu y]$.

Consider now a profile $\vec{P'_X}$, which differs from the profile $\vec{P}_X$ only by improved position of the alternative $a$. Then the cardinality of a set $X''$ which includes the alternative $a$ is reduced or remained the same while the cardinality of a set $X'''$ is higher than the cardinality of a set $X''$. Thus, $a \in C(\vec{P'_X}, X)$ and the Monotonicity condition 1 is satisfied.

### 14.6 The Monotonicity condition 2

The Monotonicity condition 2 is not satisfied (see paragraph 12.1 of this section).

### 14.7 The strict monotonicity condition

The strict monotonicity condition is not satisfied (see paragraph 1.7 of this section).

### 14.8 The non-compensatory condition

The non-compensatory condition is not satisfied (see paragraph 1.8 of this section).

## 15. The Fishburn rule

### 15.1 Heredity condition (H)

Condition **H** is not satisfied (see paragraph 12.1 of this section).

### 15.2 Concordance condition (C)

Let $X = \{a, b, c, d, e\}$ and a matrix of majority relation μ is the following

|   | a | b | c | d | e | f |
|---|---|---|---|---|---|---|
| a | - | 1 | 0 | 0 | 1 | 1 |
| b | 0 | - | 1 | 1 | 0 | 1 |
| c | 1 | 0 | - | 1 | 1 | 0 |
| d | 1 | 0 | 0 | - | 0 | 0 |
| e | 0 | 0 | 0 | 1 | - | 1 |
| f | 0 | 0 | 1 | 1 | 0 | - |

Let us define the upper contour sets for each alternative. $D(a, \vec{P}_X) = \{c, d\}, D(b, \vec{P}_X) = \{a\}$, $D(c, \vec{P}_X) = \{b, f\}$, $D(d, \vec{P}_X) = \{b, c, e, f\}$, $D(e, \vec{P}_X) = \{a, c\}$, $D(f, \vec{P}_X) = \{a, b, e\}$. Then $c\gamma d$, $b\gamma e$, $b\gamma f$. Thus, $C(\vec{P}_X, X) = \{a, b, c\}$.

Consider now the subset $X' = X \setminus \{c, f\}$. Then a matrix of majority relation µ is the following

|   | a | b | d | e |
|---|---|---|---|---|
| a | - | 1 | 0 | 1 |
| b | 0 | - | 1 | 0 |
| d | 1 | 0 | - | 0 |
| e | 0 | 0 | 1 | - |

Let us define the upper contour sets for each alternative. $D(a, \vec{P}_{X'}) = \{d\}, D(b, \vec{P}_{X'}) = \{a\}$, $D(d, \vec{P}_{X'}) = \{b, e\}, D(e, \vec{P}_{X'}) = \{a\}$. Thus, $C(\vec{P}_{X'}, X') = \{a, b, d, e\}$.

Finally, consider the subset $X'' = X \setminus \{a, b\}$. Then a matrix of majority relation µ is the following

|   | c | d | e | f |
|---|---|---|---|---|
| c | - | 1 | 1 | 0 |
| d | 0 | - | 0 | 0 |
| e | 0 | 1 | - | 1 |
| f | 1 | 1 | 0 | - |

Let us define the upper contour sets for each alternative. $D(c, \vec{P}_{X''}) = \{f\}$, $D(d, \vec{P}_{X''}) = \{c, e, f\}, D(e, \vec{P}_{X''}) = \{c\}, D(f, \vec{P}_{X''}) = \{e\}$. Then $c\gamma d$, $e\gamma d$, $f\gamma d$. Thus, $C(\vec{P}_{X''}, X'') = \{c, e, f\}$.
Then $C(\vec{P}_{X'}, X') \cap C(\vec{P}_{X''}, X'') = \{e\} \nsubseteq C(\vec{P}_X, X)$. Thus, the condition **C** is not satisfied.

*15.3 Outcast condition (O)*

Let $X = \{a, b, c, d, e\}$ and a matrix of majority relation µ is the following

|   | a | b | c | d | e |
|---|---|---|---|---|---|
| a | - | 1 | 0 | 1 | 1 |
| b | 0 | - | 1 | 1 | 0 |
| c | 1 | 0 | - | 1 | 0 |
| d | 0 | 0 | 0 | - | 1 |
| e | 0 | 0 | 0 | 0 | - |

Let us define the upper contour sets for each alternative. $D(a, \vec{P}_X) = \{c\}, D(b, \vec{P}_X) = \{a\}$, $D(c, \vec{P}_X) = \{b\}$, $D(d, \vec{P}_X) = \{a, b, c\}$, $D(e, \vec{P}_X) = \{a, d\}$. Then $c\gamma d$, $b\gamma e$. Thus, $C(\vec{P}_X, X) = \{a, b, c\}$.

Consider now the subset $X' = X \setminus \{d\}$. Then a matrix of majority relation µ is the following

|   | a | b | c | e |
|---|---|---|---|---|
| a | - | 1 | 0 | 1 |
| b | 0 | - | 1 | 0 |
| c | 1 | 0 | - | 0 |
| e | 0 | 0 | 0 | - |

Let us define the upper contour sets for each alternative. $D(a, \vec{P}_{X'}) = \{c\}, D(b, \vec{P}_{X'}) = \{a\}$, $D(c, \vec{P}_{X'}) = \{b\}, D(e, \vec{P}_{X'}) = \{a\}$. Thus, $C(\vec{P}_{X'}, X') = \{a, b, c, e\}$.
Then $C(\vec{P}_{X'}, X') \neq C(\vec{P}_X, X)$, Thus, the condition **O** is not satisfied.

### 15.4 Arrow's choice axiom (ACA)
The Fishburn rule does not satisfy the condition **ACA** as the condition **H** is not satisfied.

### 15.5 Monotonicity condition 1
Let $a \in C(\vec{P}_X, X)$. Then $\nexists y \in X\ D(y, \vec{P}_X) \subset D(a, \vec{P}_X)$.
Consider now a profile $\vec{P'_X}$, which differs from the profile $\vec{P}_X$ only by improved position of the alternative $a$. Then $(a, \vec{P'_X}) \subseteq D(a, \vec{P}_X)$ and $\forall y \in X\setminus\{a\}\ D(y, \vec{P'_X}) \supseteq D(y, \vec{P}_X)$. Thus, $\nexists y \in X$ $D(y, \vec{P'_X}) \subset D(a, \vec{P'_X})$. Thus, $a \in C(\vec{P'_X}, X)$ and monotonicity condition 1 is satisfied.

### 15.6 The Monotonicity condition 2
The Monotonicity condition 2 is not satisfied (see paragraph 12.1 of this section).

### 15.7 The strict monotonicity condition
The strict monotonicity condition is not satisfied (see paragraph 1.7 of this section).

### 15.8 The non-compensatory condition
The non-compensatory condition is not satisfied (see paragraph 1.8 of this section).

## 16. Uncovered set I
### 16.1 Heredity condition (H)
The condition **H** is not satisfied (see paragraph 12.1 of this section).

### 16.2 Concordance condition (C)
Let $X = \{a, b, c, d, e, f\}$ and a matrix of majority relation μ is the following

|   | a | b | c | d | e | f |
|---|---|---|---|---|---|---|
| a | - | 1 | 0 | 0 | 1 | 1 |
| b | 0 | - | 1 | 1 | 0 | 1 |
| c | 1 | 0 | - | 1 | 1 | 0 |
| d | 1 | 0 | 0 | - | 0 | 0 |
| e | 0 | 0 | 0 | 1 | - | 1 |
| f | 0 | 0 | 1 | 1 | 0 | - |

Let us define the lower contour sets for each alternative. $L(a, \vec{P}_X) = \{b, e, f\}, L(b, \vec{P}_X) = \{c, d, f\}$, $L(c, \vec{P}_X) = \{a, d, e\}$, $L(d, \vec{P}_X) = \{a\}$, $L(e, \vec{P}_X) = \{d, f\}$, $L(f, \vec{P}_X) = \{c, d\}$. Then $c\delta d$, $b\delta e$, $b\delta f$. Thus, $C(\vec{P}_X, X) = \{a, b, c\}$.
Consider now the subset $X' = X\setminus\{c, f\}$. Then a matrix of majority relation μ is the following

|   | a | b | d | e |
|---|---|---|---|---|
| a | - | 1 | 0 | 1 |
| b | 0 | - | 1 | 0 |
| d | 1 | 0 | - | 0 |
| e | 0 | 0 | 1 | - |

Let us define the lower contour sets for each alternative. $L(a, \vec{P}_{X'}) = \{b, e\}, L(b, \vec{P}_{X'}) = \{d\}$, $L(d, \vec{P}_{X'}) = \{a\}, L(e, \vec{P}_{X'}) = \{d\}$. Thus, $C(\vec{P}_{X'}, X') = \{a, b, d, e\}$.

Finally, consider the subset $X'' = X \backslash \{a,b\}$. Then a matrix of majority relation µ is the following

|   | c | d | e | f |
|---|---|---|---|---|
| c | - | 1 | 1 | 0 |
| d | 0 | - | 0 | 0 |
| e | 0 | 1 | - | 1 |
| f | 1 | 1 | 0 | - |

Let us define the lower contour sets for each alternative. $L(c,\vec{P}_{X''}) = \{d,e\}$, $L(d,\vec{P}_{X''}) = \emptyset$, $L(e,\vec{P}_{X''}) = \{d,f\}$, $L(f,\vec{P}_{X''}) = \{c,d\}$. Then $c\delta d$, $e\delta d$, $f\gamma d$. Thus, $C(\vec{P}_{X''}, X'') = \{c,e,f\}$.
Then $C(\vec{P}_{X'}, X') \cap C(\vec{P}_{X''}, X'') = \{e\} \nsubseteq C(\vec{P}_X, X)$. Thus, the condition **C** is not satisfied.

### 16.3 Outcast condition (O)
Let $X = \{a,b,c,d,e\}$ and a matrix of majority relation µ is the following

|   | a | b | c | d | e |
|---|---|---|---|---|---|
| a | - | 0 | 1 | 1 | 0 |
| b | 1 | - | 0 | 1 | 0 |
| c | 0 | 1 | - | 1 | 1 |
| d | 0 | 0 | 0 | - | 1 |
| e | 1 | 0 | 0 | 0 | - |

Let us define the lower contour sets for each alternative. $L(a,\vec{P}_X) = \{c,d\}$, $D(b,\vec{P}_X) = \{a,d\}$, $D(c,\vec{P}_X) = \{b,d,e\}$, $D(d,\vec{P}_X) = \{e\}$, $D(e,\vec{P}_X) = \{a\}$. Then $b\gamma e$, $c\gamma d$. Thus, $C(\vec{P}_X, X) = \{a,b,c\}$.

Consider now the subset $X' = X \backslash \{d\}$. Then a matrix of majority relation µ is the following

|   | a | b | c | e |
|---|---|---|---|---|
| a | - | 0 | 1 | 0 |
| b | 1 | - | 0 | 0 |
| c | 0 | 1 | - | 1 |
| e | 1 | 0 | 0 | - |

Let us define the lower contour sets for each alternative. $D(a,\vec{P}_{X'}) = \{c\}$, $D(b,\vec{P}_{X'}) = \{a\}$, $D(c,\vec{P}_{X'}) = \{b,e\}$, $D(e,\vec{P}_{X'}) = \{a\}$. Thus, $C(\vec{P}_{X'}, X') = \{a,b,c,e\}$.
Then $C(\vec{P}_{X'}, X') \neq C(\vec{P}_X, X)$, Thus, the condition **O** is not satisfied.

### 16.4 Arrow's choice axiom (ACA)
Uncovered set I does not satisfy the condition **ACA** as the condition **H** is not satisfied.

### 16.5 Monotonicity condition 1
Let $a \in C(\vec{P}_X, X)$. Then $\nexists y \in X\; L(y,\vec{P}_X) \supset L(a,\vec{P}_X)$.

Consider now a profile $\vec{P'_X}$, which differs from the profile $\vec{P}_X$ only by improved position of the alternative $a$. Then $(a,\vec{P'_X}) \supseteq L(a,\vec{P}_X)$ and $\forall y \in X\backslash\{a\}\; L(y,\vec{P'_X}) \subseteq L(y,\vec{P}_X)$. Thus, $\nexists y \in X\; L(y,\vec{P'_X}) \supset L(a,\vec{P'_X})$. Consequently, $a \in C(\vec{P'_X}, X)$ and monotonicity condition 1 is satisfied.

### 16.6 The Monotonicity condition 2
The Monotonicity condition 2 is not satisfied (see paragraph 12.1 of this section).

### 16.7 The strict monotonicity condition
The strict monotonicity condition is not satisfied (see paragraph 1.7 of this section).

### 16.8 The non-compensatory condition
The non-compensatory condition is not satisfied (see paragraph 1.8 of this section).

## 17. Uncovered set II
### 17.1 Heredity condition (H)
The condition **H** is not satisfied (see paragraph 12.1 of this section).

### 17.2 Concordance condition (C)
The condition **C** is satisfied iff $\forall X', X'' \epsilon 2^A \rightarrow C(\vec{P}_{X' \cup X''}, X' \cup X'') \supseteq C(\vec{P}_{X'}, X') \cap C(\vec{P}_{X''}, X'')$.
Let us proof it by contradiction.
1. Suppose that $\exists x \in C(\vec{P}_{X'}, X') \cap C(\vec{P}_{X''}, X''): x \notin C(\vec{P}_{X' \cup X''}, X' \cup X'')$. Then $\exists y \in C(\vec{P}_{X' \cup X''}, X' \cup X'') \rightarrow y \mu x$ and $D(y, \vec{P}_X) \subset D(x, \vec{P}_X)$. Moreover, $y \in X'$ and/or $y \in X''$.
2. Let $y \in X'$. Since $y\mu x$, $D(y, \vec{P}_{X'}) \subset D(x, \vec{P}_{X'})$. Thus, $y \in C(\vec{P}_{X'}, X')$ and $x \notin C(\vec{P}_{X'}, X')$ which is a contradiction. Consequently, the assumption 1 is incorrect.

Thus, Uncovered set II satisfies condition **C**.

### 17.3 Outcast condition (O)
Let $X = \{a, b, c, d, e\}$ and a matrix of majority relation μ is the following

|   | a | b | c | d |
|---|---|---|---|---|
| a | - | 0 | 1 | 0 |
| b | 1 | - | 0 | 0 |
| c | 0 | 0 | - | 1 |
| d | 1 | 1 | 0 | - |

Let us define the upper contour sets for each alternative. $D(a, \vec{P}_X) = \{b, d\}, D(b, \vec{P}_X) = \{d\}$, $D(c, \vec{P}_X) = \{a\}, D(d, \vec{P}_X) = \{c\}$. Then $bBa$. Thus, $C(\vec{P}_X, X) = \{b, c, d\}$.
Consider now the subset $X' = X \setminus \{a\}$. Then a matrix of majority relation μ is the following

|   | b | c | d |
|---|---|---|---|
| b | - | 0 | 0 |
| c | 0 | - | 1 |
| d | 1 | 0 | - |

Let us define the upper contour sets for each alternative. $D(b, \vec{P}_{X'}) = \{d\}, D(c, \vec{P}_{X'}) = \emptyset$, $D(d, \vec{P}_{X'}) = \{c\}$. Thus, $C(\vec{P}_{X'}, X') = \{b, c\}$.
Then $C(\vec{P}_{X'}, X') \neq C(\vec{P}_X, X)$, Thus, the condition **O** is not satisfied.

### 17.4 Arrow's choice axiom (ACA)
Uncovered set II does not satisfy the condition **ACA** as the condition **H** is not satisfied.

### 17.5 Monotonicity condition 1
Let $a \in C(\vec{P}_X, X)$. Then $\nexists y \in X \ y\mu x \ \& \ D(y, \vec{P}_X) \subset D(a, \vec{P}_X)$.

Consider now a profile $\overrightarrow{P'_X}$, which differs from the profile $\vec{P}_X$ only by improved position of the alternative $a$. Then $(a, \overrightarrow{P'_X}) \subseteq D(a, \vec{P}_X)$ and $\forall y \in X \setminus \{a\}$ $D(y, \overrightarrow{P'_X}) \supseteq D(y, \vec{P}_X)$. Thus, $\not\exists y \in X$ $y\mu x$ & $D(y, \overrightarrow{P'_X}) \subset D(a, \overrightarrow{P'_X})$. Thus, $a \in C(\overrightarrow{P'_X}, X)$ and monotonicity condition 1 is satisfied.

*17.6 The Monotonicity condition 2*
The Monotonicity condition 2 is not satisfied (see paragraph 12.1 of this section).

*17.7 The strict monotonicity condition*
The strict monotonicity condition is not satisfied (see paragraph 1.7 of this section).

*17.8 The non-compensatory condition*
The non-compensatory condition is not satisfied (see paragraph 1.8 of this section).

## 18. The Richelson rule
*18.1 Heredity condition (H)*
The condition **H** is not satisfied (see paragraph 12.1 of this section).

*18.2 Concordance condition (C)*
The condition **C** is not satisfied (see paragraph 15.2 of this section).

*18.3 Outcast condition (O)*
The condition **O** is not satisfied (see paragraph 15.3 of this section).

*18.4 Arrow's choice axiom (ACA)*
The Richelson rule does not satisfy the condition **ACA** as the condition **H** is not satisfied.

*18.5 Monotonicity condition 1*
The Monotonicity condition 1 is satisfied (see paragraph 15.5 and 16.5 of this section).

*18.6 The Monotonicity condition 2*
The Monotonicity condition 2 is not satisfied (see paragraph 12.1 of this section).

*18.7 The strict monotonicity condition*
The strict monotonicity condition is not satisfied (see paragraph 1.7 of this section).

*18.8 The non-compensatory condition*
The non-compensatory condition is not satisfied (see paragraph 1.8 of this section).

## 19. The Condorcet winner
*19.1 Heredity condition (H)*
Since a contraction of the profile $\vec{P}_X$ onto any set $X' \subseteq X$ does not affect the pairwise comparison of any alternatives, $C(\vec{P}_{X'}, X') = C(\vec{P}_X, X)$ if $C(\vec{P}_X, X) \neq \emptyset$. Consequently, $C(\vec{P}_{X'}, X') \supseteq C(\vec{P}_X, X) \cap X'$ and condition **H** is satisfied.

*19.2 Concordance condition (C)*
The condition **C** is satisfied iff $\forall X', X'' \in 2^A \to C(\vec{P}_{X' \cup X''}, X' \cup X'') \supseteq C(\vec{P}_{X'}, X') \cap C(\vec{P}_{X''}, X'')$.
Let us proof it by contradiction.

1. Suppose that $\exists x \in C(\vec{P}_{X'}, X') \cap C(\vec{P}_{X''}, X''): x \notin C(\vec{P}_{X' \cup X''}, X' \cup X'')$. Then $\exists y \in C(\vec{P}_{X' \cup X''}, X' \cup X''): \forall z \notin C(\vec{P}_{X' \cup X''}, X' \cup X'') \to y\mu z$. Moreover, $y \in X'$ and/or $y \in X''$.
2. Let $y \in X'$. Since $y\mu x$, $x \notin C(\vec{P}_{X'}, X')$ which is a contradiction. Thus, the assumption 1 is incorrect.

Thus, the condition **C** is satisfied.

*19.3 Outcast condition (O)*

The condition **O** is not satisfied as the choice $C(\vec{P}_{X'}, X')$ of the subset $X' \subseteq X$ can be non-empty when $C(\vec{P}_X, X) = \emptyset$.

*19.4 Arrow's choice axiom (ACA)*

The condition **ACA** is not satisfied since the first condition of the axiom ($if\ C(\vec{P}_X, X) = \emptyset$ then $C(\vec{P}_{X'}, X') = \emptyset$) is not satisfied.

*19.5 Monotonicity condition 1*

The Monotonicity condition 1 is satisfied since the advancement of the chosen alternative does not worsen its pairwise comparison with other alternatives.

*19.6 The Monotonicity condition 2*

Since the Condorcet winner always chooses no more than one best alternative, the Monotonicity condition 2 is not applicable to it as it considers the choice of more than two alternatives. In other words, the Condorcet winner obeys the Monotonicity condition 2 trivially.

*19.7 The strict monotonicity condition*

The strict monotonicity condition is not satisfied (see paragraph 1.7 of this section).

*19.8 The non-compensatory condition*

The non-compensatory condition is not satisfied (see paragraph 1.8 of this section).

## 20. The core

*20.1 Heredity condition (H)*

Let us proof it by contradiction. Suppose that $\exists x \in C(\vec{P}_X, X) \cap X': x \notin C(\vec{P}_{X'}, X')$. Then $\exists y \in X': y\mu x$. Since $y \in X\ and\ y\mu x$, $x \notin C(\vec{P}_X, X)$ which is a contradition. Thus, the previous assumption is incorrect and, consequently, the condition **H** is satisfied.

*20.2 Concordance condition (C)*

The condition **C** is satisfied (see paragraph 20.1 of this section).

*20.3 Outcast condition (O)*

Let $X = \{a, b, c\}$ and a matrix of majority relation μ is the following

|   | a | b | c |
|---|---|---|---|
| a | - | 0 | 1 |
| b | 0 | - | 0 |
| c | 0 | 1 | - |

According to the rule the alternative *a* will be chosen, i.e., $C(\vec{P}_X, X) = \{a\}$.

Consider now the subset $X' = X \setminus \{c\}$. Then a matrix of majority relation µ is the following

|   | a | b |
|---|---|---|
| a | - | 0 |
| b | 0 | - |

According to the rule the alternatives *a* and *b* will be chosen, i.e., $C(\vec{P}_{X'}, X') = \{a, b\}$.
Then $C(\vec{P}_{X'}, X') \neq C(\vec{P}_X, X)$. Thus, the condition **O** is not satisfied.

### *20.4 Arrow's choice axiom (ACA)*
The core does not satisfy the condition **ACA** as the condition **O** is not satisfied.

### *20.5 Monotonicity condition 1*
The Monotonicity condition 1 is satisfied since the advancement of the chosen alternative does not worsen its pairwise comparison with other alternatives.

### *20.6 The Monotonicity condition 2*
The Monotonicity condition 2 is not satisfied (see paragraph 20.1 of this section).

### *20.7 The strict monotonicity condition*
The strict monotonicity condition is not satisfied (see paragraph 1.7 of this section).

### *20.8 The non-compensatory condition*
The non-compensatory condition is not satisfied (see paragraph 1.8 of this section).

## 21. k-stable set
To check the properties of k-stable set similar examples with the paragraph 14 of this section can be used. For the case $k = 1$ properties of k-stable set are equal to the properties of minimal weakly stable set (see paragraph 14 of this section). For the case $k > 1$ we can use examples that differs from the paragraph 14 of this section by the addition of k-1 extra vertices (i.e., alternatives) to each edge between alternatives of minimal weakly stable set.

Let us show it for an example from the paragraph 14.2 of this section. Let $X = \{a, b, c, d, e\}$ and a majority graph is the following

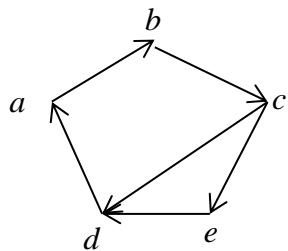

For the case $k = 1$, the alternatives *a* and *c* will be chosen, i.e., $C(\vec{P}_X, X, 1) = \{a, c\}$.
Transform this majority graph for the case $k > 1$

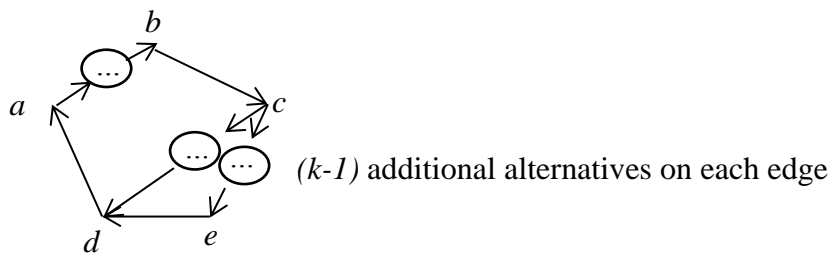

*(k-1)* additional alternatives on each edge

For the case $k > 1$, the alternatives $a$ and $c$ will be chosen, i.e., $C(\vec{P}_X, X, k) = \{a, c\}$.

## 22. The threshold rule

### 22.1 Heredity condition (H)

The condition **H** is not satisfied (see paragraph 12.1 of this section).

### 22.2 Concordance condition (C)

Let $X = \{a, b, c\}$ and the profile $\vec{P}_X$ is the following

| $P_1$ | $P_2$ | $P_3$ |
|---|---|---|
| a | b | b |
| c | a | a |
| b | c | c |

According to the threshold rule the alternative $a$ will be chosen, i.e., $C(\vec{P}_X, X) = \{a\}$.

Consider now the subset $X' = X \setminus \{a\}$. A contraction of the profile $\vec{P}_X$ onto a set $X'$, i.e., $\vec{P}_{X'}$, looks as

| $P_1$ | $P_2$ | $P_3$ |
|---|---|---|
| c | b | b |
| b | c | c |

According to the threshold rule the alternative $b$ will be chosen, i.e., $C(\vec{P}_{X'}, X') = \{b\}$.

Finally, let us consider the subset $X'' = X \setminus \{c\}$. A contraction of the profile $\vec{P}_X$ onto a set $X''$, i.e., $\vec{P}_{X''}$, looks as

| $P_1$ | $P_2$ | $P_3$ |
|---|---|---|
| a | b | b |
| b | a | a |

According to the rule the alternative $b$ will be chosen, i.e., $C(\vec{P}_{X''}, X'') = \{b\}$.

Then $C(\vec{P}_{X'}, X') \cap C(\vec{P}_{X''}, X'') = \{b\} \nsubseteq C(\vec{P}_X, X)$. Thus, the condition **C** is not satisfied.

### 22.3 Outcast condition (O)

Consider the example from paragraph 22.2 of this section. Then $C(\vec{P}_{X''}, X'') \neq C(\vec{P}_X, X)$. Thus, the condition **O** is not satisfied.

### 22.4 Arrow's choice axiom (ACA)

The threshold rule does not satisfy the condition **ACA** as the condition **H** is not satisfied.

### 22.5 Monotonicity condition 1

Let $a \in C(\vec{P}_X, X)$. Then $\nexists x \in X \setminus \{a\}\ xVa$.

Consider now a profile $\vec{P'_X}$, which differs from the profile $\vec{P}_X$ only by improved position of the alternative $a$. Then $\nexists x \in X \setminus \{a\}\ xVa$. Thus, $a \in C(\vec{P'_X}, X)$.

$a \in C(\vec{P}_X, X)$ and $a \in C(\vec{P'_X}, X)$. Thus, threshold rule satisfies the Monotonicity condition 1.

## 22.6 The Monotonicity condition 2

Let $X = \{a, b, c, d\}$ and the profile $\vec{P}_X$ is the following

| $P_1$ | $P_2$ | $P_3$ | $P_4$ | $P_5$ | $P_6$ |
|---|---|---|---|---|---|
| a | b | b | b | b | d |
| d | d | c | c | a | a |
| c | a | d | a | d | c |
| b | c | a | d | c | b |

According to the threshold rule the alternatives *a,d* will be chosen, i.e., $C(\vec{P}_X, X) = \{a, d\}$.

Consider now the subset $X' = X \backslash \{d\}$. A contraction of the profile $\vec{P}_X$ onto a set $X'$, i.e., $\vec{P}_{X'}$, looks as

| $P_1$ | $P_2$ | $P_3$ | $P_4$ | $P_5$ | $P_6$ |
|---|---|---|---|---|---|
| a | b | b | b | b | a |
| c | a | c | c | a | c |
| b | c | a | a | c | b |

According to the threshold rule the alternative *b* will be chosen, i.e., $C(\vec{P}_{X'}, X') = \{b\}$.

Finally, consider the subset $X'' = X \backslash \{a\}$. A contraction of the profile $\vec{P}_X$ onto a set $X''$, i.e., $\vec{P}_{X''}$, looks as

| $P_1$ | $P_2$ | $P_3$ | $P_4$ | $P_5$ | $P_6$ |
|---|---|---|---|---|---|
| d | b | b | b | b | d |
| c | d | c | c | d | c |
| b | c | d | d | c | b |

According to the threshold rule the alternative *b* will be chosen, i.e., $C(\vec{P}_{X''}, X'') = \{b\}$.

Then $\{a, d\} \in C(\vec{P}_X, X)$, $\{a\} \notin C(\vec{P}_{X'}, X')$ and $\{d\} \notin C(\vec{P}_{X''}, X'')$. Thus, threshold rule does not satisfy the Monotonicity condition 2.

## 22.7 The strict monotonicity condition

The threshold rule does not satisfy the strict monotonicity condition (see paragraph 1.7 of this section).

## 22.8 The non-compensatory condition

The non-compensatory condition is satisfied by the definition of the rule.

## 23. The Copeland rule 1

### 23.1 Heredity condition (H)

The condition **H** is not satisfied (see paragraph 12.1 of this section).

### 23.2 Concordance condition (C)

Let $X = \{a, b, c, d, e\}$ and a matrix of majority relation μ is the following

|   | a | b | c | d | e |
|---|---|---|---|---|---|
| a | - | 1 | 0 | 1 | 1 |
| b | 0 | - | 1 | 0 | 1 |
| c | 1 | 0 | - | 0 | 0 |
| d | 0 | 0 | 1 | - | 1 |
| e | 0 | 0 | 1 | 0 | - |

According to the Copeland rule 1 the alternative *a* will be chosen, i.e., $C(\vec{P}_X, X) = \{a\}$.

Consider now the subset $X' = X\setminus\{b\}$. Then a matrix of majority relation µ is the following

|   | a | c | d | e |
|---|---|---|---|---|
| a | - | 0 | 1 | 1 |
| c | 1 | - | 0 | 0 |
| d | 0 | 1 | - | 1 |
| e | 0 | 1 | 0 | - |

According to the Copeland rule 1 the alternatives $a$ and $d$ will be chosen, i.e., $C(\vec{P}_{X'}, X') = \{a, d\}$.

Finally, consider the subset $X'' = X\setminus\{a\}$. Then a matrix of majority relation µ is the following

|   | b | c | d | e |
|---|---|---|---|---|
| b | - | 1 | 0 | 1 |
| c | 0 | - | 0 | 0 |
| d | 0 | 1 | - | 1 |
| e | 0 | 1 | 0 | - |

According to the rule the alternatives $b$ and $d$ will be chosen, i.e., $C(\vec{P}_{X''}, X'') = \{b, d\}$.

Then $C(\vec{P}_{X'}, X') \cap C(\vec{P}_{X''}, X'') = \{d\} \nsubseteq C(\vec{P}_X, X)$. Thus, the condition **C** is not satisfied.

### 23.3 Outcast condition (O)
Consider the previous example. $C(\vec{P}_X, X) = \{a\} \neq C(\vec{P}_{X'}, X')$. Thus, the condition **O** is not satisfied.

### 23.4 Arrow's choice axiom (ACA)
The Copeland rule 1 does not satisfy the condition **ACA** as the condition **H** is not satisfied.

### 23.5 Monotonicity condition 1
Let $a \in C(\vec{P}_X, X)$. Then $\nexists y \in X \; u(y, \vec{P}_X) > u(a, \vec{P}_X)$.

Consider now a profile $\vec{P'_X}$, which differs from the profile $\vec{P}_X$ only by improved position of the alternative $a$. Then $(a, \vec{P'_X}) \geq u(a, \vec{P}_X)$ and $\forall y \in X\setminus\{a\} \; (y, \vec{P'_X}) \leq u(y, \vec{P}_X)$. Thus, $\nexists y \in X \; u(y, \vec{P'_X}) > u(a, \vec{P'_X})$. Consequently, $a \in C(\vec{P'_X}, X)$ and monotonicity condition 1 is satisfied.

### 23.6 The Monotonicity condition 2
The Monotonicity condition 2 is not satisfied (see paragraph 12.1 of this section).

### 23.7 The strict monotonicity condition
The strict monotonicity condition is not satisfied (see paragraph 1.7 of this section).

### 23.8 The non-compensatory condition
The non-compensatory condition is not satisfied (see paragraph 1.8 of this section).

## 24. The Copeland rule 2
### 24.1 Heredity condition (H)
The condition **H** is not satisfied (see paragraph 12.1 of this section).

### 24.2 Concordance condition (C)
The condition **C** is not satisfied (see paragraph 23.2 of this section).

### 24.3 Outcast condition (O)
The condition **O** is not satisfied (see paragraph 23.3 of this section).

### 24.4 Arrow's choice axiom (ACA)
Copeland rule 2 does not satisfy the condition **ACA** as the condition **H** is not satisfied.

### 24.5 Monotonicity condition 1
Let $a \in C(\vec{P}_X, X)$. Then $\forall y \in X$ $|L(a, \vec{P}_X)| \geq |L(y, \vec{P}_X)|$.

Consider now a profile $\vec{P'_X}$, which differs from the profile $\vec{P}_X$ only by improved position of the alternative $a$. Then $|L(a, \vec{P'_X})| \geq |L(a, \vec{P}_X)|$ and $\forall y \in X \backslash \{a\}$ $|L(y, \vec{P'_X})| \leq |L(y, \vec{P}_X)|$. Thus, $\forall y \in X$ $|L(a, \vec{P'_X})| \geq |L(y, \vec{P'_X})|$. Consequently, $a \in C(\vec{P'_X}, X)$ and monotonicity condition 1 is satisfied.

### 24.6 The Monotonicity condition 2
The Monotonicity condition 2 is not satisfied (see paragraph 12.1 of this section).

### 24.7 The strict monotonicity condition
The strict monotonicity condition is not satisfied (see paragraph 1.7 of this section).

### 24.8 The non-compensatory condition
The non-compensatory condition is not satisfied (see paragraph 1.8 of this section).

## 25. The Copeland rule 3

### 25.1 Heredity condition (H)
The condition **H** is not satisfied (see paragraph 12.1 of this section).

### 25.2 Concordance condition (C)
Let $X = \{a, b, c, d, e, f\}$ and a matrix of majority relation μ is the following

|   | a | b | c | d | e | f |
|---|---|---|---|---|---|---|
| a | - | 1 | 0 | 0 | 1 | 1 |
| b | 0 | - | 1 | 1 | 0 | 1 |
| c | 1 | 0 | - | 1 | 1 | 0 |
| d | 1 | 0 | 0 | - | 0 | 0 |
| e | 0 | 0 | 0 | 1 | - | 1 |
| f | 0 | 0 | 1 | 1 | 0 | - |

According to the Copeland rule 3 the alternative $b$ will be chosen, i.e., $C(\vec{P}_X, X) = \{b\}$.

Consider now the subset $X' = X \backslash \{c, f\}$. Then a matrix of majority relation μ is the following

|   | a | b | d | e |
|---|---|---|---|---|
| a | - | 1 | 0 | 1 |
| b | 0 | - | 1 | 0 |
| d | 1 | 0 | - | 0 |
| e | 0 | 0 | 1 | - |

According to the Copeland rule 3 the alternatives $a,b,e$ will be chosen, i.e., $C(\vec{P}_{X'}, X') = \{a, b, e\}$.

Finally, consider the subset $X'' = X \setminus \{a, b\}$. Then a matrix of majority relation μ is the following

|   | c | d | e | f |
|---|---|---|---|---|
| c | - | 1 | 1 | 0 |
| d | 0 | - | 0 | 0 |
| e | 0 | 1 | - | 1 |
| f | 1 | 1 | 0 | - |

According to the rule the alternatives $c,e,f$ will be chosen, i.e., $C(\vec{P}_{X''}, X'') = \{c, e, f\}$.
Then $C(\vec{P}_{X'}, X') \cap C(\vec{P}_{X''}, X'') = \{e\} \nsubseteq C(\vec{P}_X, X)$. Thus, the condition **C** is not satisfied.

### 25.3 Outcast condition (O)
Consider the previous example. $C(\vec{P}_X, X) = \{b\} \neq C(\vec{P}_{X'}, X')$. Thus, the condition **O** is not satisfied.

### 25.4 Arrow's choice axiom (ACA)
The Copeland rule 3 does not satisfy the condition **ACA** as the condition **H** is not satisfied.

### 25.5 Monotonicity condition 1
Let $a \in C(\vec{P}_X, X)$. Then $\forall y \in X \; |D(a, \vec{P}_X)| \leq |D(y, \vec{P}_X)|$.
Consider now a profile $\vec{P'_X}$, which differs from the profile $\vec{P}_X$ only by improved position of the alternative $a$. Then $\left|D(a, \vec{P'_X})\right| \leq |D(a, \vec{P}_X)|$ and $\forall y \in X \setminus \{a\} \; \left|D(y, \vec{P'_X})\right| \geq |D(y, \vec{P}_X)|$. Thus, $\forall y \in X \; \left|D(a, \vec{P'_X})\right| \leq |D(y, \vec{P'_X})|$. Consequently, $a \in C(\vec{P'_X}, X)$ and monotonicity condition 1 is satisfied.

### 25.6 The Monotonicity condition 2
The Monotonicity condition 2 is not satisfied (see paragraph 12.1 of this section).

### 25.7 The strict monotonicity condition
The strict monotonicity condition is not satisfied (see paragraph 1.7 of this section).

### 25.8 The non-compensatory condition
The non-compensatory condition is not satisfied (see paragraph 1.8 of this section).

## 26a. The super-threshold rule (fixed threshold)
The super-threshold rule with fixed threshold value chooses alternatives for which the criterion value is more or less than some fixed threshold value. Hence, the results of such choice procedure do not depend on the initial set $X$. Thus, the super-threshold choice rule with fixed threshold value satisfies all given condition except the non-compensatory condition (see paragraph 1.8 of this section).

## 26b. The super-threshold rule (threshold depends on X)
### 26.1 Heredity condition (H)
The condition **H** is satisfied iff $\forall X, X' \in 2^A, X' \subseteq X \Rightarrow C(\vec{P}_{X'}, X') \supseteq C(\vec{P}_X, X) \cap X'$. Let us proof it by contradiction.
Suppose that the condition **H** is not satisfied, i.e., $\exists x \in C(\vec{P}_X, X) \cap X': x \notin C(\vec{P}_{X'}, X')$. Then $\varphi(x, \vec{P}_X) \geq V(X, \vec{P}_X)$ and $\varphi(x, \vec{P}_{X'}) < V(X', \vec{P}_{X'})$. In other words, there should be a subset $X' \subseteq X$ such that the threshold value $V(X', \vec{P}_{X'})$ will be more that $V(X, \vec{P}_X)$ and some alternative $x$ will not be chosen. Obviously, this situation may occur, consequently, the condition **H** is not satisfied.

*26.2 Concordance condition (C)*

The condition **C** is not satisfied (see paragraph 26.1 of this section).

*26.3 Outcast condition (O)*

The condition **O** is not satisfied (see paragraph 26.1 of this section).

*26.4 Arrow's choice axiom (ACA)*

The condition **ACA** is not satisfied as the condition **H** is not satisfied.

*26.5 Monotonicity condition 1*

Let $a \in C(\vec{P}_X, X)$. Then $\varphi(a, \vec{P}_X) \geq V(X, \vec{P}_X)$. Consider now a profile $\vec{P'_X}$, which differs from the profile $\vec{P}_X$ only by improved position of the alternative $a$. Then $\varphi(a, \vec{P'_X}) \geq V(X, \vec{P'_X})$ since $\varphi(a, \vec{P'_X}) - \varphi(a, \vec{P}_X) \geq V(X, \vec{P'_X}) - V(X, \vec{P}_X)$. Thus, the Monotonicity condition 1 is satisfied.

*26.6 The Monotonicity condition 2*

Let $\{a, b\} \in C(\vec{P}_X, X)$. Then $\varphi(a, \vec{P}_X) \geq V(X, \vec{P}_X)$ and $\varphi(b, \vec{P}_X) \geq V(X, \vec{P}_X)$. By definition it is obvious that the threshold value does not increase when one of the chosen alternative is eliminated. Consequently, $\varphi(a, \vec{P}_{X \setminus \{b\}}) \geq V(X \setminus \{b\}, \vec{P}_{X \setminus \{b\}})$ and $\varphi(b, \vec{P}_{X \setminus \{a\}}) \geq V(X \setminus \{a\}, \vec{P}_{X \setminus \{a\}})$. Then $\{a, b\} \in C(\vec{P}_X, X)$, $\{a\} \in C(\vec{P}_{X \setminus \{b\}}, X \setminus \{b\})$ and $\{b\} \in C(\vec{P}_{X \setminus \{a\}}, X \setminus \{a\})$. Thus, the Monotonicity condition 2 is satisfied.

*26.7 The strict monotonicity condition*

The super-threshold rule with the threshold value that depends on $X$ does not satisfy the strict monotonicity condition (see paragraph 26.1 of this section).

*26.8 The non-compensatory condition*

The non-compensatory condition is not satisfied (see paragraph 1.8 of this section).

## 27. The minimax procedure

*27.1 Heredity condition (H)*

Let $X = \{a, b, c\}$ and the profile $\vec{P}_X$ looks as

| $P_1$ | $P_2$ | $P_3$ | $P_4$ | $P_5$ | $P_6$ | $P_7$ | $P_8$ | $P_9$ | $P_{10}$ |
|---|---|---|---|---|---|---|---|---|---|
| a | a | a | a | b | b | b | c | c | c |
| b | b | b | b | c | c | c | a | a | a |
| c | c | c | c | a | a | a | b | b | b |

Let us construct a matrix $S^-(\vec{P}_X, X)$ for the profile $\vec{P}_X$.

|   | a | b | c |
|---|---|---|---|
| a | - | 6 | 8 |
| b | 7 | - | 5 |
| c | 5 | 8 | - |

According to the minimax procedure the alternative $a$ will be chosen, i.e., $C(\vec{P}_X, X) = \{a\}$.

Consider now the subset $X' = X \setminus \{b\}$. A contraction of the profile $\vec{P}_X$ onto a set $X'$, i.e., $\vec{P}_{X'}$, looks as

| $P_1$ | $P_2$ | $P_3$ | $P_4$ | $P_5$ | $P_6$ | $P_7$ | $P_8$ | $P_9$ | $P_{10}$ |
|---|---|---|---|---|---|---|---|---|---|
| a | a | a | a | c | c | c | c | c | c |
| c | c | c | c | a | a | a | a | a | a |

According to the rule the alternative $c$ will be chosen, i.e., $C(\vec{P}_{X'}, X') = \{c\}$.

Then $C(\vec{P}_{X'}, X') \not\supseteq C(\vec{P}_X, X) \cap X'$. Thus, the condition **H** is not satisfied.

### 27.2 Concordance condition (C)

Let $X = \{a, b, c, d\}$ and the profile $\vec{P}_X$ looks as

| $P_1$ | $P_2$ | $P_3$ | $P_4$ | $P_5$ | $P_6$ | $P_7$ | $P_8$ | $P_9$ | $P_{10}$ | $P_{11}$ |
|---|---|---|---|---|---|---|---|---|---|---|
| d | a | a | a | a | c | b | b | b | d | a |
| a | d | d | c | c | d | c | c | c | a | d |
| b | b | b | d | d | b | d | d | d | b | b |
| c | c | c | b | b | a | a | a | a | c | c |

Let us construct a matrix $S^-(\vec{P}_X, X)$ for the profile $\vec{P}_X$.

|   | a | b | c | d |
|---|---|---|---|---|
| a | - | 5 | 5 | 4 |
| b | 4 | - | 6 | 3 |
| c | 4 | 3 | - | 6 |
| d | 5 | 6 | 3 | - |

According to the minimax procedure the alternative $a$ will be chosen, i.e., $C(\vec{P}_X, X) = \{a\}$.

Consider now the subset $X' = X \setminus \{c\}$. A contraction of the profile $\vec{P}_X$ onto a set $X'$, i.e., $\vec{P}_{X'}$, looks as

| $P_1$ | $P_2$ | $P_3$ | $P_4$ | $P_5$ | $P_6$ | $P_7$ | $P_8$ | $P_9$ | $P_{10}$ | $P_{11}$ |
|---|---|---|---|---|---|---|---|---|---|---|
| d | a | a | a | a | d | b | b | b | d | a |
| a | d | d | d | d | b | d | d | d | a | d |
| b | b | b | b | b | a | a | a | a | b | b |

Let us construct a matrix $S^-(\vec{P}_{X'}, X')$ for the profile $\vec{P}_{X'}$.

|   | a | b | d |
|---|---|---|---|
| a | - | 5 | 4 |
| b | 4 | - | 3 |
| d | 5 | 6 | - |

According to the minimax procedure the alternative $d$ will be chosen, i.e., $C(\vec{P}_{X'}, X') = \{d\}$.

Finally, consider the subset $X'' = X \setminus \{a\}$. A contraction of the profile $\vec{P}_X$ onto a set $X''$, i.e., $\vec{P}_{X''}$, looks as

| $P_1$ | $P_2$ | $P_3$ | $P_4$ | $P_5$ | $P_6$ | $P_7$ | $P_8$ | $P_9$ | $P_{10}$ | $P_{11}$ |
|---|---|---|---|---|---|---|---|---|---|---|
| d | d | d | c | c | c | b | b | b | d | d |
| b | b | b | d | d | d | c | c | c | b | b |
| c | c | c | b | b | b | d | d | d | c | c |

Let us construct a matrix $S^-(\vec{P}_{X''}, X'')$ for the profile $\vec{P}_{X''}$.

|   | b | c | d |
|---|---|---|---|
| b | - | 6 | 3 |
| c | 3 | - | 6 |
| d | 6 | 3 | - |

According to the rule the alternatives $b,c,d$ will be chosen, i.e., $C(\vec{P}_{X''}, X'') = \{b, c, d\}$.

Then $C(\vec{P}_{X'}, X') \cap C(\vec{P}_{X''}, X'') = \{d\} \nsubseteq C(\vec{P}_X, X)$. Thus, the condition **C** is not satisfied.

### 27.3 Outcast condition (O)

Let $X = \{a, b, d\}$ and the profile $\vec{P}_X$ looks as

| $P_1$ | $P_2$ | $P_3$ | $P_4$ | $P_5$ |
|---|---|---|---|---|
| a | a | d | b | b |
| b | b | a | d | d |
| d | d | b | a | a |

Let us construct a matrix $S^-(\vec{P}_X, X)$ for the profile $\vec{P}_X$

|   | a | b | d |
|---|---|---|---|
| a | - | 3 | 2 |
| b | 2 | - | 4 |
| d | 3 | 1 | - |

According to the minimax procedure the alternatives $a,b$ will be chosen, i.e., $C(\vec{P}_X, X) = \{a, b\}$.

Consider now the subset $X' = X \setminus \{d\}$. A contraction of the profile $\vec{P}_X$ onto a set $X'$, i.e., $\vec{P}_{X'}$, looks as

| $P_1$ | $P_2$ | $P_3$ | $P_4$ | $P_5$ |
|---|---|---|---|---|
| a | a | a | b | b |
| b | b | b | a | a |

Let us construct a matrix $S^-(\vec{P}_{X'}, X')$ for the profile $\vec{P}_{X'}$

|   | a | b |
|---|---|---|
| a | - | 3 |
| b | 2 | - |

According to the rule the alternative $a$ will be chosen, i.e., $C(\vec{P}_{X'}, X') = \{a\}$.

Then $C(\vec{P}_{X'}, X') \neq C(\vec{P}_X, X)$. Thus, the condition **O** is not satisfied.

### 27.4 Arrow's choice axiom (ACA)

The minimax procedure does not satisfy the condition **ACA** as the condition **H** is not satisfied.

### 27.5 Monotonicity condition 1

Let $a \in \arg\min_{b \in X} \max_{c \in X}\{n(b, c, \vec{P}_X)\}$, i.e., $a \in C(\vec{P}_X, X)$.

Consider now a profile $\vec{P'_X}$, which differs from the profile $\vec{P}_X$ only by improved position of the alternative $a$. Then $\forall x, y \in X \setminus \{a\}\ n\left(x, a, \vec{P'_X}\right) \leq n(x, a, \vec{P}_X)$, $n\left(a, x, \vec{P'_X}\right) \geq n(a, x, \vec{P}_X)$, $n\left(x, y, \vec{P'_X}\right) = n(x, y, \vec{P}_X)$ and $n\left(y, x, \vec{P'_X}\right) = n(y, x, \vec{P}_X)$. Thus, $a \in \arg\min_{b \in X} \max_{c \in X}\{n(b, c, \vec{P'_X})\}$, i.e., $a \in C\left(\vec{P'_X}, X\right)$.

Then $a \in C(\vec{P}_X, X)$ and $a \in C(\vec{P'_X}, X)$. Thus, the Monotonicity condition 1 is satisfied.

### 27.6 The Monotonicity condition 2

Let $X = \{a, b, c, d\}$ and the profile $\vec{P}_X$ looks as

| $P_1$ | $P_2$ | $P_3$ | $P_4$ | $P_5$ | $P_6$ | $P_7$ | $P_8$ | $P_9$ | $P_{10}$ | $P_{11}$ |
|---|---|---|---|---|---|---|---|---|---|---|
| d | a | a | a | a | c | b | b | b | d | a |
| a | d | d | c | c | d | c | c | c | a | d |
| b | b | b | d | d | b | d | d | d | b | b |
| c | c | c | b | b | a | a | a | a | c | c |

Let us construct a matrix $S^-(\vec{P}_X, X)$ for the profile $\vec{P}_X$

|   | a | b | c | d |
|---|---|---|---|---|
| a | - | 3 | 2 | 1 |
| b | 1 | - | 3 | 2 |
| c | 2 | 1 | - | 3 |
| d | 3 | 2 | 1 | - |

According to the minimax procedure the alternatives $a,b,c,d$ will be chosen, i.e., $C(\vec{P}_X, X) = \{a, b, c, d\}$.

Consider now the subset $X' = X\setminus\{d\}$. A contraction of the profile $\vec{P}_X$ onto a set $X'$, i.e., $\vec{P}_{X'}$, looks as

| $P_1$ | $P_2$ | $P_3$ | $P_4$ | $P_5$ | $P_6$ | $P_7$ | $P_8$ | $P_9$ | $P_{10}$ | $P_{11}$ |
|---|---|---|---|---|---|---|---|---|---|---|
| a | a | a | a | a | c | b | b | b | a | a |
| b | b | b | c | c | b | c | c | c | b | b |
| c | c | c | b | b | a | a | a | a | c | c |

Let us construct a matrix $S^-(\vec{P}_{X'}, X')$ for the profile $\vec{P}_{X'}$

|   | a | b | c |
|---|---|---|---|
| a | - | 3 | 2 |
| b | 1 | - | 3 |
| c | 2 | 1 | - |

According to the rule the alternative $a$ will be chosen, i.e., $C(\vec{P}_{X'}, X') = \{a\}$.

Then $\{b, c\} \in C(\vec{P}_X, X)$ and $\{b, c\} \notin C(\vec{P}_{X'}, X')$. Thus, the Monotonicity condition 2 is not satisfied.

### 27.7 The strict monotonicity condition

Let $X = \{a, b, c\}$ and the profile $\vec{P}_X$ looks as

| $P_1$ | $P_2$ | $P_3$ |
|---|---|---|
| a | a | b |
| c | b | c |
| b | c | a |

Let us construct a matrix $S^-(\vec{P}_X, X)$ for the profile $\vec{P}_X$

|   | a | b | c |
|---|---|---|---|
| a | - | 2 | 2 |
| b | 1 | - | 2 |
| c | 1 | 1 | - |

According to the minimax procedure the alternative $a$ will be chosen, i.e., $C(\vec{P}_X, X) = \{a\}$.

Consider now a profile $\overrightarrow{P'_X}$, which differs from the profile $\vec{P}_X$ only by improved position of the alternative $c$ in $P'_1$:

| $P'_1$ | $P'_2$ | $P'_3$ |
|---|---|---|
| c | a | b |
| a | b | c |
| b | c | a |

Let us construct a matrix $S^-(\overrightarrow{P'_X}, X)$ for the profile $\overrightarrow{P'_X}$

|   | a | b | c |
|---|---|---|---|
| a | - | 2 | 1 |
| b | 1 | - | 2 |
| c | 2 | 1 | - |

According to the rule the alternatives $a,b,c$ will be chosen, i.e., $C(\overrightarrow{P'_X}, X) = \{a, b, c\}$.

$$C\left(\overrightarrow{P'_X}, X\right) \neq \begin{bmatrix} C(\vec{P}_X, X) \text{ or} \\ \{c\} \text{ or} \\ C(\vec{P}_X, X) \cup \{c\}. \end{bmatrix}$$

Thus, strict monotonicity condition is not satisfied.

### 27.8 The non-compensatory condition

Let $X = \{a, b, c\}$ and the profile $\vec{P}_X$ is the following

| $P_1$ | $P_2$ | $P_3$ |
|---|---|---|
| a | a | b |
| b | b | c |
| c | c | a |

Let us construct a matrix $S^-(\vec{P}_X, X)$ for the profile $\vec{P}_X$

|   | a | b | c |
|---|---|---|---|
| a | - | 2 | 2 |
| b | 1 | - | 3 |
| c | 1 | 0 | - |

According to the rule the alternative $a$ will be chosen, i.e., $C(\vec{P}_X, X) = \{a\}$.

Let us write the profile $\vec{P}_X$ in the following form

|   | $\varphi_1$ | $\varphi_2$ | $\varphi_3$ |
|---|---|---|---|
| a | 3 | 3 | 1 |
| b | 2 | 2 | 3 |
| c | 1 | 1 | 2 |

According to the non-compensatory condition the alternative $b$ is better than the alternative $a$ and the alternative $a$ is better than the alternative $c$. Thus, the non-compensatory condition is not satisfied as $\{b\} \neq C(\vec{P}_X, X)$.

## 28. The Simpson procedure

The Simpson procedure satisfies the same conditions as the minimax procedure. To prove it examples from paragraphs 27.1-27.8 of this section can be used.

# References

1. Aizerman M., Aleskerov F. Theory of Choice. Elsevier, North-Holland, 1995.
2. Aleskerov F.T., Khabina E.L., Shvarts D.A. Binary relations, graphs and collective decisions. Moscow: HSE Publishing House, 2006 (in Russian).
3. Aleskerov F.T., Kurbanov E. On the degree of manipulability of group choice rules // Automation and Remote Control (in Russian). 1998. No. 10. P. 134–146.
4. Volsky V.I. Voting rules in small groups from ancient times to the XX century (in Russian). Working paper WP7/2014/02. Moscow: HSE Publishing House, 2014.
5. Aleskerov F. Procedures of multicriterial choice. Preprints of the IFAC/IFORS Conference on Control Science and Technology for Development, Beijing, China, 1985. P. 858–869.
6. Aleskerov F. Multicriterial interval choice models, Information Sciences. 1994. No. 1. P. 14–26.
7. Aleskerov F.T., Subochev A. Matrix-vector representation of various solution concepts. Working papers WP7/2009/03. Moscow: HSE Publishing House, 2009.
8. Subochev A. Dominant, Weakly Stable, Uncovered Sets: Properties and Extensions. Working papers WP7/2008/03. Moscow: HSE Publishing House, 2008.
9. Aleskerov F.T., Yuzbashev D.V., Yakuba V.I. Threshold aggregation for three-graded rankings // Automation and Remote Control (in Russian). 2007. No. 1. P. 147–152.
10. Aleskerov F. Threshold Utility, Choice and Binary Relations // Automation and Remote Control. 2003. No. 3. P. 8–27.
11. Aleskerov F. Categories of Arrovian Voting Schemes // Handbook of Economics 19, Handbook of Social Choice and Welfare. Vol. 1 / K. Arrow, A. Sen, K. Suzumura (eds). Elsevier Science B.V., 2002. P. 95–129.
12. Moulin H. Axioms of cooperative decision-making. Cambridge, Cambridge University Press, 1987.
13. Aleskerov F. Arrovian Aggregation Models, Kluwer Academic Publishers, Dordrecht, Boston, London, 1999.
14. Aleskerov F., Chistyakov V., Kalyagin V. Social threshold aggregations // Social Choice and Welfare. 2010. Vol. 35. No. 4. P. 627–646.

*Швыдун С.В.* – НИУ ВШЭ, ИПУ РАН.

Швыдун С. В.

**Нормативные свойства процедур
многокритериального выбора и их суперпозиции: I**

(*на английском языке*)